\documentclass[12pt]{article}
\usepackage{amsmath}
\usepackage{amsfonts}
\usepackage{amssymb}
\usepackage{latexsym}

\newcommand{\comp}{{\operatorname{comp}}}
\newcommand{\arccomp}{{\operatorname{arccomp}}}
\newcommand{\dist}{{\operatorname{dist}}}
\newcommand{\id}{{\operatorname{id}}}

\newcommand{\qed}{\hfill$\Box$}
\newcommand{\R}{{\mathbb R}}

\newtheorem{theorem}{Theorem}[section]
\newtheorem{prop}[theorem]{Proposition}
\newtheorem{lemma}[theorem]{Lemma}
\newtheorem{coro}[theorem]{Corollary}

\newtheorem{remark}[theorem]{Remark} 
\newtheorem{remarks}[theorem]{Remarks}

\numberwithin{equation}{section}
\DeclareMathSymbol\leqslant{\mathrel}{AMSa}{"36}
\DeclareMathSymbol\geqslant{\mathrel}{AMSa}{"3E}

\newenvironment{proof}{\par\noindent{\bf Proof.} }{\par\noindent}
\newenvironment{example}{\par\noindent{\bf Example.} }{\par\noindent}
\newenvironment{defi}{\par\noindent{\bf Definition.} }{\par\noindent}

\newcommand{\qu}{\quad}
\newcommand{\en}{\enspace}
\newcommand{\bi}{\bigskip}
\newcommand{\me}{\medskip}
\newcommand{\no}{\noindent}    
\newcommand{\non}{\nonumber}                

\newcommand{\bu}{\bigtriangleup}                                               
\newcommand{\be}{\begin{equation}}                                              
\newcommand{\ee}{\end{equation}}                                                 
\newcommand{\bea}{\begin{eqnarray}}                                              
\newcommand{\eea}{\end{eqnarray}} 
\newcommand{\lra}{\longrightarrow}

\newcommand{\ra}{\rightarrow}       

\newcommand{\N}{{\rm I}\!{\rm N}}  
\newcommand{\Z}{\makebox[0.06cm][l]{\sf Z}{\sf Z}}

\newcommand{\C}{{\mbox{\rm $\scriptscriptstyle ^\mid$\hspace{-0.40em}C}}} 

\renewcommand{\epsilon}{\varepsilon}
\renewcommand{\phi}{\varphi}

\def\theequation{\thesection.\arabic{equation}}

\newcommand{\dil}{{\rm dil \,\,}}

\newcommand{\geom}{{\rm Geom}}
\newcommand{\gencomp}{{\rm gen \, comp}}
\newcommand{\homm}{{\rm Hom \,\,}}
\newcommand{\mult}{{\rm Mult}}
\newcommand{\clm}{{\rm CLM}}
\newcommand{\cl}{{\rm CL}}
\newcommand{\grad}{{\rm \,\, grad \,\,}}
\newcommand{\im}{{\rm im \,\,}}
\newcommand{\supp}{{\rm supp \,\,}}
\newcommand{\ind}{{\rm ind}}
\newcommand{\dimker}{{\rm dim \,\, ker \,\,}}
\newcommand{\tr}{{\rm tr \,\,}}
\newcommand{\re}{{\rm Re \,}}

\newcommand{\sk}[2]{\mbox{$\langle#1,#2\rangle$}}

\newfont{\klein}{cmmi5}

\title{Relative Zeta Functions, Determinants, Torsion, Index Theorems and Invariants
for Open Manifolds}
\author{J\"urgen Eichhorn, Greifswald}

\begin{document}

\sloppy

\maketitle

\begin{abstract}
The set of Clifford bundles of bounded geometry over open mani\-folds can be endowed 
with a metrizable uniform structure. For one fixed bundle $E$ we define the generalized
component $\gencomp (E)$ as the set of Clifford bundles $E'$ which have finite distance
to $E$. If $D$, $D'$ are the associated generalized Dirac operators, we prove for the
pair $(D,D')$ relative index theorems, define relative $\zeta$-- and $\eta$--functions, 
relative determinants and in the case of $D=\Delta$ relative analytic torsion. To
define relative $\zeta$-- and $\eta$--functions, we assume additionally that the
essential spectrum of $D^2$ has a gap above zero.
\end{abstract}

\section{Introduction}

\setcounter{equation}{0}

The succeeding paper has 5 main goals, 1. a contribution to classification theory
for open manifolds, 2. very general relative index theory, 3. scattering theory,
4. spectral invariants like the absolutely continuous spectrum and relative
analytic torsion, 5. relative determinants, important in QFT on open manifolds.
All these goals are achieved for pairs of Clifford bundles
$(E,h,\nabla^h,\cdot) \lra (M^n,g)$, 
$(E',h',\nabla^{h'},\cdot') \lra ({M'}^n,g')$
of bounded geometry and associated generalized Dirac operators
$D=D(E,h,\nabla^h,\cdot,g)$, 
$D'=D(E',h',\nabla^{h'},\cdot',g')$
if the following holds. For certain
$K \subset M$, $K' \subset M'$,
\be
(E|_{M \setminus K} \lra M \setminus K), 
(E'|_{M' \setminus K'} \lra M' \setminus K')
\en \mbox{ are vector bundle isomorphic, }
\end{equation}
and (after pulling back)
$g-g'$, $h-h'$, $\nabla^h - \nabla^{h'}$, $\cdot - \cdot'$
have finite $(p=1,r+1)$--Sobolev norm, 
$r+1>n+2$ or $r+1>n+3$, i.e.
\bea
|g-g'|_{g,1,r+1} < \infty, \quad |h-h'|_{g,h,\nabla^h,1,r+1} < \infty, \non \\
|\nabla^h-\nabla^{h'}|_{g,h,\nabla^h,1,r+1} < \infty, \quad 
|\cdot-\cdot'|_{g,h,\nabla^h,1,r+1} < \infty.
\eea
For relative determinants and analytic torsion we assume additionally
\be
\inf \sigma_e (\Delta(g)|_{(\ker \Delta(g))^\perp}) > 0,
\end{equation}
where $\sigma_e(\Delta(g))$ denotes the essential spectrum of the graded
Laplace operator $\Delta(g)$ and we do not admit compact topological perturbations.

We rewrite these conditions in a more convenient language
$E' \in \gencomp^{1,r+1}_{L,diff,rel}(E)$
where $\gencomp(\cdot)$ are generalized components in a uniform
structure on the set of all Clifford bundles of bounded geometry.
The description of (1.1), (1.2) as $E' \in \comp(E)$ is more
convenient, contentful and useful since this expresses additionally
the symmetry and transitivity of (1.1), (1.2) which is very important
in performing proofs. Sobolev estimates e. g. play an absolutely 
decisive role. These are always performed for concrete Sobolev norms
but are needed for others too. 
$E', E'' \in \comp^r(E)$ 
implies the equivalence of 
$|\en|_{g,h,\nabla^h,1,r}$, 
$|\en|_{g',h',\nabla^{h'},p,r}$, 
$|\en|_{g'',h'',\nabla^{h''},p,r}$.
Moreover, 
$\comp(E)$ and $\gencomp(E)$
are endowed with a topology such that continuity and differentiability
considerations make sense.

Special cases of our results presented in sections 9 and 10 have been discussed
e. g. in [1], [17], [2], [3].

Our approach in this paper is organized as follows. We repeat in section
2 very briefly the main concepts concerning generalized Dirac operators
and Sobolev spaces. A key role in the sequel is played by the module
structure theorem 2.6. We sketch in sections 3 and 4 uniform structures
of open manifolds, metrics, vector bundles, connections and Clifford
bundles. This is a way to fit these objects into natural intrinsic
equivalence classes, the arc components of the uniform structure under
consideration. The general approach is to define metrizable uniform
structures by local radial metrics which are locally ''symmetric and
transitive''. If the corresponding uniform topology is locally arcwise
connected then the components coincide with arc components and the 
whole space under consideration (consisting of proper metric spaces or
open, complete Riemannian manifolds) is the topological sum of its (arc)
components. In some important cases the uniform struture is not locally
arcwise connected. Here we had to make a choice, to restrict to arc
components or to generalized components, where
$\gencomp (X,d_X) = \{ (Y,d_Y) | d_{\mathfrak U}(X,Y) < \infty \}$.
We decided to work with (the much bigger) $\gencomp$ since this is for
the later applications the most important case. The generalized 
components $\gencomp(E)$ are the biggest equivalence classes such that
we can define for any pair of operators in $\gencomp (E)$ our invariants 
and prove relative index theorems.

Our starting basic uniform structure is the Lipschitz uniform structure
defined on the set (of isometry classes) of proper metric spaces. Thereafter
we restrict to complete Riemannian manifolds and smooth maps, to vector
bundles and Clifford bundles of bounded geometry. Finally we admit compact
topological perturbations in the definition of $d(E,E') < \infty$. Then
the arising uniform structures are not locally arcwise connected. The
topological background for this is that in general one cannot connect
two different manifolds by a reasonable arc consisting of manifolds only.
We omit mostly in these sections the proofs since they are either
performed in [6] or modeled by proofs in this paper.

Section 5 sums up some heat kernel estimates which are permanently used
in the sequel.

Section 6 has model character for our approach. We fix a Clifford bundle
$(E,h,\nabla^h,\cdot) \lra (M^n,g)$
of bounded geometry, the corresponding genera\-li\-zed Dirac operator
$D=D(g,h,\nabla^h,\cdot)$
and permit variation of the (metric) connection $\nabla^h$, such that
$E'=(E,h,{\nabla'}^h,\cdot) \in \comp^{1,r}(E,h,\nabla^h,\cdot)$, 
in particular 
$|\nabla-\nabla'|_{g,h,\nabla^h,1,r} < \infty$.
Then $D$ and $D'=D(g,h,{\nabla'}^h,\cdot)$ are self adjoint in the
same Hilbert space $L_2((M,E),g,h)$.
$e^{-tD^2}-e^{-t{D'}^2}$
is defined and we want to prove the trace class property for 
$t>0$. By Duhamel's principle
\be
e^{-tD^2}-e^{-t{D'}^2} = \int\limits^t_0 e^{-sD^2} ({D'}^2-D^2) 
e^{-(t-s){D'}^2} \,\, ds
\end{equation}
If we write 
$D^2-{D'}^2 = D(D-D')+(D-D')D'$,
$D'-D = \eta = \eta^{op}$, 
where
$\eta^{op} (\psi) = \sum\limits^n_{i=1} e_i \cdot \eta_{e_i} (\psi)$, 
then
\[ e^{-tD^2}-e^{-t{D'}^2} = \int\limits^t_0 e^{-sD^2} D \eta e^{-(t-s){D'}^2} \,\, ds +
\int\limits^t_0 e^{-sD^2} \eta D' e^{-(t-s){D'}^2}\,\, ds .  \]
We split 
$\int\limits^t_0 = \int\limits^{\frac{t}{2}}_0 + \int\limits^t_{\frac{t}{2}}$
and obtain
\def\theequation{$I_1$}   
\be
e^{-tD^2}-e^{-t{D'}^2} 
= \int\limits^{\frac{t}{2}}_0 e^{-sD^2} D \eta e^{-(t-s){D'}^2} \,\, ds 
\end{equation}
\def\theequation{$I_2$}   
\be
+ \int\limits^{\frac{t}{2}}_0 e^{-sD^2} \eta D' e^{-(t-s){D'}^2} \,\, ds 
\end{equation}
\def\theequation{$I_3$}   
\be 
+ \int\limits^t_{\frac{t}{2}} e^{-sD^2} D \eta e^{-(t-s){D'}^2} \,\, ds 
\end{equation}
\def\theequation{$I_4$}   
\be
+ \int\limits^t_{\frac{t}{2}} e^{-sD^2} \eta D' e^{-(t-s){D'}^2} \,\, ds.
\end{equation}
\def\theequation{\thesection.\arabic{equation}}  
\setcounter{equation}{4}

We want to show that each integral ${\rm (I_1) - (I_4)}$ is a product of
Hilbert--Schmidt operators and to estimate their Hilbert--Schmidt norm.
Consider the integrand of ${\rm (I_4)}$,
\[ (e^{-sD^2} \eta)(D' e^{-(t-s){D'}^2}). \]

One has 
$|e^{-(t-s){D'}^2}|_{L_2 \ra H^1} \le C \cdot (t-s)^{-\frac{1}{2}}$.
If one writes
\be
(e^{-sD^2} \eta)(D' e^{-(t-s){D'}^2}) = (e^{-\frac{s}{2}D^2} f)
(f^{-1} e^{-\frac{s}{2}D^2} \eta) (D'e^{-(t-s){D'}^2})
\end{equation}
where $f$ shall be scalar function acting by multiplication then we would
be done if 
$e^{-\frac{s}{2}D^2} f$,
$f^{-1} e^{-\frac{s}{2}D^2} \eta$
would be Hilbert--Schmidt with bounded HS--norm on compact $t$--intervals
$[a_0,a_1]$, $a_0>0$.
For 
$e^{-sD^2}f$
this holds if $f \in L_2$.
The absolutely decisive question is whether $f \in L_2$ can be
additionally chosen in such a manner that
\be
f^{-1} e^{-\frac{s}{2}D^2} \eta
\end{equation}
is Hilbert--Schmidt too. This is in fact possible via rather delicate
estimates. Here we use
$|\eta|_{1,r} < \infty$.

$(I_1)$ -- $(I_3)$ can be discussed quite parallel. By a similar
decomposition we show that 
$De^{-tD^2}-D'e^{-t{D'}^2}$
is for $t>0$ also of trace class and its trace norm is uniformly bounded
on compact $t$--intervalls $[a_0,a_1]$, $a_0>0$. 
In section 7, we admit complete perturbation of the Clifford structure, 
i.e. of $g,h,\nabla^h,\cdot$.
This has several unconvenient consequences, e. g. variation of $g,h$
implies that $D=D(g,h,\nabla^h,\cdot)$
and $D'=D(g',h',\nabla^{h'},\cdot')$
are self adjoint only in distinct Hilbert spaces, hence
$e^{-tD^2}-e^{-t{D'}^2}$
is not defined. We transform $D'$ from
$L'_2=L_2((M,E),g',h')$ to $L_2=L_2((M,E),g,h)$,
thus getting an operator
$D^{'2}_{L_2} = (U^*i^*D'iU)^2$.
But for later applications in section 9 we do this in two steps.
In the first step we admit perturbation of
$g,\nabla^h,\cdot$ to $g',{\nabla'}^h,\cdot'$, 
i.e. the fibre metric $h$ remains fixed. For this reason we write
$E' \in \comp^{1,r+2}_{L,diff,F}(E)$
($F$ stands for fixing of the fibre metric). Then
$D^{'2}_{L_2} = (U^*D'U)^2$.
We apply an adapted version of Duhamel's principle thus getting
\be
e^{-tD^2} - e^{-tD^{'2}_{L_2}} = - e^{-tD^2} (\alpha-1) - \int\limits^t_0
e^{-sD^2} (D^2-D^{'2}_{L_2}) e^{-(t-s)D^{'2}_{L_2}} \,\, ds,
\end{equation}
where $\alpha-1=\frac{\sqrt{\det g}}{\sqrt{\det g'}}-1$.
The term 
$e^{-tD^2}(\alpha-1)$
can be easily be settled since $\alpha-1$ is $(1,r+1)$--Sobolev.
Then with 
$-(D-D')=\eta$ we have to estimate
\bea
&& \int\limits^{\frac{t}{2}}_0 e^{-sD^2} D (\eta - \frac{\grad' \alpha}{2 \alpha}) 
e^{-(t-s){D'}^2_{L_2}} \,\, ds + \non \\
&& \int\limits^{\frac{t}{2}}_0 e^{-sD^2} (\eta - \frac{\grad' \alpha}{2 \alpha}) 
D'_{L_2} e^{-(t-s){D'}^2_{L_2}} \,\, ds + \non \\
&& \int\limits^t_{\frac{t}{2}} e^{-sD^2} D (\eta - \frac{\grad' \alpha}{2 \alpha}) 
e^{-(t-s){D'}^2_{L_2}} \,\, ds + \non \\
&& \int\limits^t_{\frac{t}{2}} e^{-sD^2} (\eta - \frac{\grad' \alpha}{2 \alpha}) 
D'_{L_2} e^{-(t-s){D'}^2_{L_2}} \,\, ds . \non 
\eea
We can write 
$-(D-D') \Phi = (\eta^{op}_1+\eta^{op}_2+\eta^{op}_3) \Phi$.
If we still set 
$\eta_0 = - \frac{\grad' \alpha \cdot'}{\alpha}$
then we must study the 16 integrals
\def\theequation{$I_{\nu,1}$}   
\be
   \int\limits^{\frac{t}{2}}_0 e^{-sD^2} D \eta_\nu e^{-(t-s){D'}^2_{L_2}} \,\, ds ,
\end{equation}
\def\theequation{$I_{\nu,2}$}  
\be
   \int\limits^{\frac{t}{2}}_0 e^{-sD^2} \eta_\nu D'_{L_2} e^{-(t-s){D'}^2_{L_2}} \,\, ds ,
\end{equation}
\def\theequation{$I_{\nu,3}$}  
\be
   \int\limits^t_{\frac{t}{2}} e^{-sD^2} D \eta_\nu  e^{-(t-s){D'}^2_{L_2}} \,\, ds ,
\end{equation}
\def\theequation{$I_{\nu,4}$}  
\be
   \int\limits^t_{\frac{t}{2}} e^{-sD^2} \eta_\nu D'_{L_2} e^{-(t-s){D'}^2_{L_2}} \,\, ds .
\end{equation}                         
\def\theequation{\thesection.\arabic{equation}}  
\setcounter{equation}{7} 

Unfortunately $\eta^{op}_1$ and $\eta^{op}_3$ are first order differential operators
but using estimates for the first derivatives of the heat kernels we can show that
for $t>0$
\be
e^{-tD^2} - e^{-tD^{'2}_{L_2}} \quad \mbox{and} \quad
e^{-tD^2} D - e^{-tD^{'2}_{L_2}} D'_{L_2}
\end{equation}
are of trace class and their trace norm is uniformly bounded on compact
$t$--intervalls $[a_0,a_1]$, $a_0>0$.

There is a very important class of Clifford structures where perturbation of $g$
automatically induces a perturbation of the fibre metric $h$. This is
$E = (\Lambda^* T^* M \otimes \C, g_{\Lambda^*}, \nabla^{g_{\Lambda^*}})
\lra (M^n,g)$ 
with Clifford multiplication
\[ X \otimes \omega \in T_m M \otimes \Lambda^* T^* M \otimes \C \lra X \cdot \omega
= \omega_X \wedge \omega - i_X \omega, \]
where $\omega_X := g(, X)$. In this case $E$ as a vector bundle remains fixed but the
Clifford module structure varies smoothly with $g,g' \in \comp(g)$. It is well known
that in this case $D=d+d^*$, $D^2=(d+d^*)^2$ Laplace operator $\Delta$.
This example later produces the relative analytic torsion.
Therefore we must admit perturbation of the fibre metric too but we can prove that
if 
$E' \in \gencomp^{1,r+1}_{L,diff}(E)$
then for $t>0$
\be
e^{-tD^2} - e^{-t(U^*i^*D'iU)^2} \quad \mbox{and} \quad
e^{-tD^2} D - e^{-t(U^*i^*D'iU)^2} (U^*i^*D'iU)
\end{equation}
are of trace class and their trace norm is uniformly bounded on compact
$t$--intervalls $[a_0,a_1]$, $a_0>0$.

In section 8 finally, we additionally admit compact topological perturbations, 
$E|_{M \setminus K} \cong E'|_{M' \setminus K'}$
as vector bundles together with (1.2).
Consider
${\cal H} = L_2 ((K,E|_K),g,h) \oplus L_2 ((K',E'|_{K'}),g'h') 
\oplus L_2 ((M \setminus K),g,h)$.
Then $D$ acts in ${\cal H}$ by defining it to be zero on 
$L_2 ((K',E'|_{K'}),g'h')$. 
Similarly, a suitable transformation $\tilde{D'}$ of $D'$ acts in 
${\cal H}$ by defining it additionally to be zero on 
$L_2 ((K,E|_K),g,h)$. 
Let $P, P'$ be the orthogonal projection onto 
$L_2 ((K',E'|_{K'}),g'h')^\perp$, $L_2 ((K,E|_K),g,h)^\perp$,
respectively.
Then we have to prove
(if $E' \in \gencomp^{1,r+1}_{L,diff,rel}(E)$)
that for $t>0$
\be
e^{-tD^2} P - e ^{-t{\tilde{D'}}^2} P' \quad \mbox{and} \quad
e^{-tD^2} D - e^{-t{\tilde{D'}}^2} {\tilde{D'}}^2
\end{equation}
are of trace class and their trace norm is uniformly bounded on compact
$t$--intervalls $[a_0,a_1]$, $a_0>0$.

In principle (1.10) contains the other theorems (as special case
$K=K'=\Phi$) such that we could omit them. But there is no doubling
of proofs. The proof for the more general perturbation adds to the
proof of the foregoing case (with smaller admitted perturbations) only
those features which come specifically from the more general perturbation.

In the special case of the graded Laplace operator we obtain that for
$t>0$
\be
e^{-t\Delta} P - e^{-t(U^*i^*\Delta'iU)} P'
\end{equation}
is of trace class and its trace norm is uniformly bounded on compact
$t$--intervalls $[a_0,a_1]$, $a_0>0$.

After a long preparatory material, the sections 9 and 10 contain the 
applications to relative index theory, scattering theory, $\zeta$--functions,
determinants and relative torsion. Here we assume as usual $E,E'$ endowed
with an involution = $\Z_2$--grading
$\tau$, $\tau^2=1$, $\tau^*=\tau$, $[\tau,X\cdot]_+=0$
for $X \in TM$, $[\nabla,\tau]=0$.
Then 
$L_2((M,E),g,h) = L_2(M,E^+) \oplus L_2(M,E^-)$
and
$D = \left( \begin{array}{cc} 0 & D^- \\ D^+ & 0  \end{array} \right)$,
$D^- = (D^+)^*$.
If $M^n$ is compact then as usual
\[ \ind D := \ind D^+ := \dimker D^+ - \dimker D^- \equiv \tr (\tau e^{-tD^2}), \]
where we understand $\tau$ as
\[ \tau = \left( \begin{array}{cc} I & 0 \\ 0 & -I \end{array} \right) . \]
For open $M^n$ $\ind D$ in general is not defined since $\tau e^{-tD^2}$ is not
of trace class. The appropriate approach on open manifolds is relative index 
theory for pairs of operators $D,D'$. If $D,D'$ are selfadjoint in the same
Hilbert space and $e^{tD^2}-e^{-t{D'}^2}$ would be of trace class then
\be
\ind (D,D') := \tr (\tau (e^{-tD^2} - e^{-t{D'}^2}))
\end{equation}
makes sense, but at the first glance (1.12) should depend on $t$.

The admitted perturbations then must be $\tau$--compatible, i.e.
$[\nabla',\tau]=0$, $[\tau,X\cdot']_+=0$.
Then we prove in theorem 9.1 that if 
$E' \in \comp^{1,r+1}_{L,diff,rel} (E)$
then
\be
\tr \tau (e^{-tD^2} - e^{-t(U^*i^*D'iU)^2})
\end{equation}
is independent of $t$. A similar result has been proved in [2]
under much more restrictive assumptions.

We conclude in 9.4 that under the hypothesis of 9.1 the pair
$D$, $U^*i^*D'iU$ 
forms a supersymmetric scattering system. In particular,
the restriction of 
$D$, $U^*i^*D'iU$ 
to their absolutely continious spectral
subspaces are unitarily equivalent.

If moreover, 
$\inf \sigma_e (D^2) > 0$ 
then
$\inf \sigma_e ((U^*i^*D'iU)^2) > 0$ 
and
\[ \tr \tau (e^{-tD^2} - e^{-t(U^*i^*D'iU)^2}) = \ind D - \ind D' , \]
where as usual 
$\ind D := \ind D^+$. 
This is the content of theorem 9.6.
In theorem 9.7 we admit more general perturbations for the relative index
theorem. Thereafter we introduce the spectral shift function and show that
under certain assumptions the scattering index 
$n^c (\lambda,D,\tilde{D'})$
is constant and if additionally the essential spectrum has a gap above zero
then even 
$n^c (D,\tilde{D'}) = 0$.

The concluding section 10 is devoted to relative $\zeta$--functions, 
$\eta$--functions, determinants and analytic torsion. Here we prove in theorem
10.4 that for 
$E' \in \gencomp^{1,r+1}_{L,diff,F} (E)$, 
$\inf \sigma_e (D^2|_{(\ker D^2)^\perp}) > 0$
a relative zeta--function
$\zeta (s,D^2,\tilde{D'}^2)$
is defined for $\re s > -1$
and is holomorphic in $s=0$.
This allows to define relative determinants which are very important in QFT
and which satisfy the usual rules. This is the content of theorem 10.5.

Unfortunately
$g' \in \comp^{1,r+1}(g)$
does not imply
$E' \in \gencomp^{1,r+1}_{L,diff,F} (E)$
for
$E = ( \Lambda^* \otimes \C, g_{\Lambda^*})$, 
$E' = ( \Lambda^* \otimes \C, g'_{\Lambda^*})$,
since the fibre metric also changes, 
$g_{\Lambda^*} \lra g'_{\Lambda^*}$.
Here we have to refine our considerations. This is performed in 10.6. 
We prove in theorem 10.7 that for
$g' \in \comp^{1,r+1}(g)$
and
$\inf \sigma_e (\Delta |_{(\ker \Delta)^\perp}) > 0$
the relative analytic torsion $\tau_a$,
$\log \tau_a (M^n,g,g') := \sum\limits^n_{q=0} (-1)^q q \cdot
\zeta'_q (0,\Delta,\Delta')$
is well defined. Thereafter we present classes of examples where 
all assumptions are satisfied.

The final calculations are devoted to the relative $\eta$--function and
$\eta$--invariant which is the content of theorem 10.11.

Part of our future attention will turn to combinatorial approaches of the 
topic of this paper.

\section{Clifford bundles, generalized Dirac operators and Sobolev spaces}

\setcounter{equation}{0}

We recall for completeness very briefly the basic properties of 
generalized Dirac operators on open manifolds. Let $(M^n,g)$ be a
Riemannian manifold, $m \in M, Cl(T_mM,g_m)$ the corresponding Clifford
algebra at $m$. $Cl(T_m,g_m)$ shall be complexified or not, depending
on the other bundles and structure under consideration. A hermitian
vector bundle $E \rightarrow M$ is called a bundle of Clifford modules
if each fibre $E_m$ is a Clifford module over $Cl(T_m,g_m)$ with
skew symmetric Clifford multiplication. We assume $E$ to be endowed
with a compatible connection $\nabla^E$, i.e. $\nabla^E$ is metric and 
\[ \nabla^E_X (Y \cdot \Phi) = (\nabla^g_X Y) \cdot \Phi + Y \cdot  
(\nabla^E_X \Phi) \]
$X,Y \in \Gamma (TM), \Phi \in \Gamma(E)$. Then we call the pair 
$(E,\nabla^E)$ a Clifford bundle. The composition 
\[ \Gamma(E) \stackrel{\nabla}{\longrightarrow} \Gamma (T^* M \otimes E) 
   \stackrel{g}{\longrightarrow} \Gamma ( T M \otimes E) 
   \stackrel{\cdot}{\longrightarrow} \Gamma (E) \]        
shall be called the generalized Dirac operator $D$. We have
$D=D(g,E,\nabla)$. If $X_1, \dots X_n$ is an orthonormal basis in $T_mM$ then
\[ D = \sum^n_{i=1} X_i \cdot \nabla^E_{X_i} . \]
$D$ is of first order elliptic, formally self-adjoint and 
\[ D^2 = \bu^E + R , \]
where $\bu^E = (\nabla^E)^* \nabla^E$ and $ R \in \Gamma(End(E))$ is the 
bundle endomorphism
\[ R \Phi = \frac{1}{2} \sum^n_{i,j=1} X_i X_j R^E (X_i,X_j) \Phi . \]
Next we recall some associated functional spaces and their properties
if we assume bounded geometry. These facts are contained in [10], [11] and
partially in [13].

Let $E \rightarrow M$ be a Clifford bundle, $\nabla = \nabla^E$, $D$ the 
generalized Dirac operator. Then we define for $\Phi \in \Gamma(E),
p \ge 1, r \in {\bf Z}, r \ge 0$,  
\begin{eqnarray} 
   |\Phi|_{W^{p,r}} &:=&  \left( 
   \int \sum^r_{i=0} | \nabla^i \Phi |^p_x \,\, dvol_x(g)
   \right)^{\frac{1}{p}} , \nonumber \\
   |\Phi|_{H^{p,r}} &:=&  \left( 
   \int \sum^r_{i=0} | D^i \Phi |^p_x \,\, dvol_x(g)
   \right)^{\frac{1}{p}} , \nonumber \\  
   W^p_r(E) &:=& \left\{ 
   \Phi \in \Gamma(E) \big| |\Phi|_{W^{p,r}} < \infty
   \right\} , \nonumber \\
   W^{p,r}(E) &:=& \en \mbox{completion of} 
   \en W^p_r \en \mbox{w. r. t.}
   \en | \en |_{W^{p,r}} , \nonumber \\
   H^p_r(E) &:=& \left\{ 
   \Phi \in \Gamma(E) \big| |\Phi|_{H^{p,r}} < \infty
   \right\} , \nonumber \\
   H^{p,r}(E) &:=& \en \mbox{completion of} 
   \en H^p_r \en \mbox{w. r. t.}
   \en | \en |_{H^{p,r}} . \nonumber 
\end{eqnarray}
In a great part of our consideration we restrict to $p=1,2$. In the case $p=2$
we write $W^{2,r} \equiv W^r, H^{2,r} \equiv H^r $ etc.. If $r<0$ then
we set 
\begin{eqnarray} 
W^r(E) &:=& \Big( W^{-r}(E) \Big)^* , \nonumber \\    
H^r(E) &:=& \Big( H^{-r}(E) \Big)^* . \nonumber 
\end{eqnarray}
Assume $(M^n,g)$ complete. Then $C^\infty_c(E)$ is a dense subspace of 
$W^{p,1}(E)$ and $H^{p,1}(E)$. This follows from proposition 1.4 in [11].
If we use this density and the fact
\[ |D \Phi(m)| \le C \cdot |\nabla \Phi(m)| , \]
we obtain $|\Phi|_{H^{p,1}} \le C' \cdot |\Phi|_{W^{p,1}}$ and a
continious embedding
\[ W^{p,1}(E) \hookrightarrow H^{p,1}(E) . \]
For $r>1$ this cannot be established, and we need further assumptions.
Consider as in the introduction the following conditions 
\be 
\begin{array}{lll}
(I) & \qquad & 
\qu r_{inj} (M,g) = \inf_{x \in M} r_{inj}(x) > 0 , \\
(B_k(M,g)) & & 
\qu |(\nabla^g)^i R^g| \le C_i, \en 0 \le i \le k , \\
(B_k(E,\nabla^E)) & & 
\qu |(\nabla^g)^i R^E| \le C_i, \en 0 \le i \le k . 
\end{array} \non
\end{equation}
It is a well known fact that for any open manifold and given $k, 
0 \le k \le \infty$, there exists a metric $g$ satisfying $(I)$ and 
$(B_k(M,g))$. Moreover, $(I)$ implies completeness of $g$.

\begin{lemma}  
Assume $(M^n,g)$ with (I) and $(B_k)$. Then $C^\infty_c(E)$ is a dense
subset of $W^{p,r}(E)$ and $H^{p,r}(E)$ for $0 \le r \le k+2$.
\end{lemma}

See [11], proposition 1.6 for a proof. \qed

\begin{lemma}
Assume $(M^n,g)$ with (I) and $(B_k)$. Then there exists a
conti\-nuous embedding
\[ W^{p,r}(E) \longrightarrow H^{p,r}(E), \quad 0 \le r \le k+1 . \]
\end{lemma}

\begin{proof} 
According to 2.1, we are done if we can prove
\[ |\Phi|_{H^{p,r}} \le C \cdot |\Phi|_{W^{p,r}} \]
for $0 \le r \le k+1$ and $\Phi \in C^\infty_c(E)$. Perform induction.
For $r=0$ $|\Phi|_{H^{p,0}} = |\Phi|_{W^{p,0}}$. Assume 
$|\Phi|_{H^{p,r}} \le C \cdot |\Phi|_{W^{p,r}}$. Then
\begin{eqnarray} 
  |\Phi|_{H^{p,r+1}} & \le & C \cdot ( |\Phi|_{H^{p,r}} +
  |D^r D \Phi|_{H^{p,r}})  \nonumber \\  
  & \le & C \cdot ( |\Phi|_{W^{p,r}} + |D \Phi|_{W^{p,r}}) . \nonumber
\end{eqnarray}
Let 
$\frac{\partial}{\partial x^i}, i=1, \dots, n$ 
coordinate vectors fields which are orthonormal in $m \in M$. 
Then with 
$\nabla_i = \nabla_{\frac{\partial}{\partial x^i}}$
\[ |\nabla^s D \Phi|^p_m \le C \cdot \sum_{i_1, \dots ,i_s, j} 
   |\nabla_{i_1} \dots \nabla_{i_s} \frac{\partial}{\partial x^j} \cdot
   \nabla_j \Phi|^p . \]   
Now we apply the Leibniz rule and use the fact that in an atlas of 
normal charts the Christoffel symbols have bounded euclidean 
derivatives up to order $k-1$. This yields
\[ |\nabla^r D \Phi|^p_m \le C \cdot \sum_{i_1, \dots ,i_{r+1}} 
   |\nabla_{i_1} \dots \nabla_{i_{r+1}} \Phi|^p_m \en \mbox{for} \en
   r \le k , \]   
i. e.
\[ |D\Phi|_{W^{p,r}} \le C \cdot |\Phi|_{W^{p,r+1}} \]    
altogether
\[ |\Phi|_{H^{p,r+1}} \le C \cdot |\Phi|_{W^{p,r+1}}. \]
\qed 
\end{proof}

\me
\no
{\bf Remark.} For $p=2$ this proof is contained in [3]. \hfill $\Box$

\begin{theorem} 
Assume $(M^n,g)$ with (I) and $(B_k)$ and
$(E,\nabla)$ with $(B_k)$ and $p=2$. Then for $r \le k$
\[ H^{2,r}(E) \equiv H^r(E) \cong W^r(E) \equiv W^{2,r}(E) \]
as equivalent Hilbert spaces.
\end{theorem}

\begin{proof} 
According to 2.2., $W^r(E) \subseteq H^r (E)$
continuously. Hence we have to show $H^r(E) \subseteq W^r(E)$
continuously. The latter follows from the local elliptic 
inequality, a uniformly locally finite cover by normal charts of
fixed radius, uniform trivializations and the existence of 
uniform elliptic constants. The proof is performed in [3].
\qed
\end{proof}

\me
\no
{\bf Remark.} 2.3 holds for $1 < p < \infty$ (cf. [13]). 
\qed

\me
As it is clear from the definition that the spaces $W^{p,k}(E)$ can be
defined for any Riemannian vector bundle $(E, h_E, \nabla^E)$. We
assume this more general case and define additionally
\[ {}^{b,s}W(E) := \left\{ \varrho \in C^s(E) \en \Big| \en 
   {}^{b,s}|\varrho| := \sum^s_{i=0} \sup_{x \in M}  
   |\nabla^i \varrho|_x < \infty   \right\} \]          
and in the case of a Clifford bundle
\[ {}^{b,s}H(E) := \left\{ \varrho \in C^s(E) \en \Big| \en 
   {}^{b,s,D}|\varrho| := \sum^s_{i=0} \sup_{x \in M}  
   |D^i \varrho|_x < \infty   \right\} . \]          
${}^{b,s}W(E)$ is a Banach space and coincides with the completion
of the space of all $\varrho \in \Gamma(E)$ with 
${}^{b,s}|\varrho| < \infty$ with respect to ${}^{b,s}| |$.

\begin{theorem} 
Let $(E,h,\nabla^E)$ be a Riemannian vector
bundle satisfying (I), $(B_k(M^n,g))$, $B_k(E;\nabla))$.

{\bf a)} Assume $k \ge r, k \ge 1, r-\frac{n}{p} \ge s-\frac{n}{q},
r \ge s, q \ge p $, then
\be
W^{p,r}(E) \hookrightarrow W^{q,s}(E) 
\end{equation}
continuously.

{\bf b)} If $k \ge r > \frac{n}{p} + s$ then
\be 
W^{p,r}(E) \hookrightarrow {}^{b,s}W(E) 
\end{equation}
continuously.
\end{theorem}

\me
We refer to [13] for the proof. 
\qed

\begin{coro} 
Let $E \rightarrow M$ be a Clifford bundle
satisfying (I), $(B_k(M))$, $(B_k(E))$, $k \ge r $, $r > \frac{n}{2}+s$. Then
\be 
H^r(E) \hookrightarrow {}^{b,s}H(E) 
\end{equation}
continuously.
\end{coro}

\begin{proof} 
We apply 2.3, (2.2) and obtain
\be
H^r(E) \hookrightarrow {}^{b,s}W(E). 
\end{equation}    
Quite similar as in the proof of 2.2., 
\be
H^r(E) \hookrightarrow {}^{b,s}W(E). 
\end{equation}    
continuously. 
\qed
\end{proof}

A key role for everything in the sequel plays the module structure theorem
for Sobolev spaces.

\begin{theorem}
Let $(E_i,h_i,D_i) \rightarrow (M^n,g)$
  be vector bundles with $(I)$, $(B_k(M^n,g))$, $(B_k(E_i,\nabla_i))$,
  $i=1,2$. Assume $0\le r\le r_1,r_2\le k$. If $r=0$ assume
  \[
  \left\{ 
    \begin{array}{rcl}
      r-\frac{n}{p} & < & r_1-\frac{n}{p_1} \\
      r-\frac{n}{p} & < & r_2-\frac{n}{p_2} \\
      r-\frac{n}{p} & \le & r_1-\frac{n}{p_1} + r_2-\frac{n}{p_2} \\
      \frac{1}{p} & \le & \frac{1}{p_1} +\frac{1}{p_2}
    \end{array}
  \right\}
  \en \mbox{or} \]
\be
  \left\{ 
    \begin{array}{rcl}
      r-\frac{n}{p} &  \le  & r_1-\frac{n}{p_1} \\
      0  & <  & r_2-\frac{n}{p_2} \\
      \frac{1}{p} & \le  & \frac{1}{p_1}
    \end{array}
  \right\}
  \en \mbox{or} \en
  \left\{ 
    \begin{array}{rcl}
      0  & <  & r_1-\frac{n}{p_1} \\
      r-\frac{n}{p} &  \le  & r_2-\frac{n}{p_2} \\
      \frac{1}{p} & \le  & \frac{1}{p_2}
    \end{array}
  \right\}.
\end{equation}
  If $r>0$ assume $\frac{1}{p}\le \frac{1}{p_1} + \frac{1}{p_2}$ and
\be
  \left\{ 
    \begin{array}{rcl}
      r-\frac{n}{p} & < & r_1-\frac{n}{p_1} \\
      r-\frac{n}{p} & < & r_2-\frac{n}{p_2} \\
      r-\frac{n}{p} & \le & r_1-\frac{n}{p_1} + r_2-\frac{n}{p_2}
    \end{array}
  \right\}
  \en \mbox{or} \en
  \left\{ 
    \begin{array}{rcl}
      r-\frac{n}{p} & \le & r_1-\frac{n}{p_1} \\
      r-\frac{n}{p} & \le & r_2-\frac{n}{p_2} \\
      r-\frac{n}{p} & < & r_1-\frac{n}{p_1} + r_2-\frac{n}{p_2}
    \end{array}
  \right\}.
\end{equation}
Then the tensor product of sections defines a continuous bilinear map
\[ W^{p_1,r_1}(E_1,\nabla_1) \times W^{p_2,r_2}(E_2,\nabla_2) \longrightarrow
   W^{p,r} (E_1 \otimes E_2, \nabla_1 \otimes \nabla_2) . \]
\end{theorem}

We refer to [13] for the proof. \qed

Define for $u \in C^0(M), c>0$ 
\[ \overline{u}_c(x) := \frac{1}{vol B_c(x)} 
   \int\limits_{B_c(x)} u(y) \,\, d vol_y (g) . \]

\begin{lemma}
Let $(M^n,g)$ be complete, $Ric(g) \ge k, k \in {\bf R}, k>0$.
Then there exists a positive constant $C=C(n,k,R)$, depending only on 
$h,k,R$ such that for any $c \in  ]0,R[$ and any 
$u \in W^{1,1}(M) \cap C^\infty (M)$
\[ \int\limits_M |u-\overline{u}_c| d vol_x (g) \le  
   C \cdot c \cdot \int\limits_M |\nabla u|  \,\, d vol_x (g) .  \]
\end{lemma}

\begin{proof} 
For $u \in C_c^\infty (M)$ the proof is performed in [18],
p. 31--33. But what is only needed in the proof is 
$\int |u| dx, \int |\nabla u| \, dx < \infty $
(even only 
$\int |\nabla u| \, dx < \infty $).
   
The key is a the lemma of Buser, 
\[ \int\limits_{B_c(x)} |u-\overline{u}_c| \,\, dy \le
   C \cdot c \cdot \int\limits_{B_c(x)} |\nabla u| \,\, dy . \]
\qed
\end{proof}

\me
\no
{\bf Remark.} The conditions $u, \nabla u \in C^\infty$ are not necessary, 
$u \in C^1$ is sufficient. 

\qed

\begin{prop}
Let $(E,h,\nabla) \rightarrow (M^n,g)$ be a Riemannian vector bundle,
$(M^n,g)$ with (I), $(B_0), r>n+1, 0<c<r_{inj}$ and 
$\eta \in W^{1,r}(E)$. Then 
$\overline{|\eta|}_c \in W^{1,0}(M) \equiv L_1(M)$, where
\[ \overline{|\eta|}_c(x) := \frac{1}{vol B_c(x)}
   \int\limits_{B_c(x)} |\eta(y)| \,\, dy . \]
\end{prop}

\begin{proof} 
Set $u(x)=|\eta(x)|$. Then   
$\overline{u(x)}=\overline{|\eta|}_c(x)$ and, according to Kato's
inequality, 
\[ \int |\nabla u| dx = \int |\nabla |\eta|| dx \le 
   \int |\nabla \eta| \,\, dx < \infty . \]
Hence we obtain from 2.7, $|u| = |\eta| \in L_1, 
|\eta| - \overline{|\eta|_c} \in L_1$,
\be
\overline{|\eta|}_c \in L_1 .
\end{equation}
\qed
\end{proof}

\me
{\bf Remark.} For (2.2) is the assumption $(B_0(E))$ superfluous.
Nevertheless, in our applications we even have $(B_k(E))$.
\qed

\me
Finally we recall for clarity and distinctness a fact which will be
very important later. Let $(E,h,\nabla) \rightarrow (M^n,g)$ be a
Riemannian vector bundle with $(I), (B_k(M)), (B_k(E)), h \ge r+1,
r > \frac{n}{p}+1, 0 < c < r_{inj}$. Then the spaces 
$W^{p,r}(E |_{B_c(x)}) = \{ \varrho \Big| \varrho$ distributional section
of $E |_{B_c(x)}$ s. t. $|\varrho|_{p,r} < \infty \}$ are well defined,
$x \in M$ arbitrary. Radial parallel translation of an orthonormal basis
defines an isomorphism
\be
A_x : W^{p,r} (E |_{B_c(x)}) \stackrel{\cong}{\longrightarrow}
W^{p,r}(B_c(0),V^N),
\end{equation}                 
$B_c(0) \subset {\bf R}^n, V^N = {\bf R}^N$ or ${\bf C}^N, N = \mbox{rk} E$.
We conclude from $(B_k(M)), (B_k(E))$, $k \ge r+1$ and [15] that there
exists constants $c_1, C_1$ s.t.
\be
c_1 \cdot |\varrho|_{p,r,B_c(x)} \le  |A_x \varrho|_{p,r,B_c(0)}  
\le C_1 \cdot |\varrho|_{p,r,B_c(x)} ,  
\end{equation}                                   
$c_1, C_1$ independent of $x$. Moreover, if $\varrho \in W^{p,r}(E)$
then $\varrho |_{B_c(x)} \in W^{p,r} (E |_{B_c(x)})$. Similarly, 
\be
c_2 {}^{b,s} |\varrho|_{B_c(x)} \le {}^{b,s} |A_x \varrho|_{B_c(0)}  
\le C_2 \cdot {}^{b,s}|\varrho|_{B_c(x)} ,  
\end{equation}
$c_2, C_2$ independent of $x$. $(B_k(M)), (B_k(E)),0 < c < r_{inj}$
imply that $B_c(x)$ satisfies all required smoothness conditions and we
obtain from the Sobolev embedding theorem, (2.10), (2.11)
\be
\hspace*{0.7cm} 
W^{p,r} (E |_{B_c(x)}) \hookrightarrow {}^{b,1} W(E |_{B_c(x)}) ,
\end{equation}                                 
\be
{}^{b,1}|\varrho|_{B_c(x)} \le C \cdot |\varrho|_{p,r,B_c(x)} ,
\end{equation}
$C$ independent of $x$.

\section{Uniform structures of open manifolds, metric, vector bundles, connections and
Clifford structures}

\setcounter{equation}{0}

As we pointed out in [6], there exist in any dimension $n \ge 2$ for closed manifolds
only countably many diffeomorphism types but for open manifolds uncountably many homotopy
types. Moreover, there exists no nontrivial additive number valued invariant which is 
defined for all open (oriented) manifolds. We attack these two problems as follows. We 
developed in [6], [7], [8] a classification approach for open manifolds consisting
of 4 main steps. 

1. Definition of uniform structures of open manifolds.

2. Characterization of the (arc) components which serve as rough equivalence classes.

3. Classification of this components by invariants.

4. Classification of the elements inside a component by invariants.

The components in the uniform structure under consideration are the rough equivalence
classes. In the for us essential cases they contain only countably many diffeomorphism
types. In this paper, we are concerned essentially with step 4 above. We define number
valued invariants for the manifolds inside a component. But according to the remark
above, this will be relative invariants. This means we fix one manifold $M_0$ inside
a component and construct for any other manifold $M$ an invariant defined by the pair
$(M_0,M)$. At the purely Riemannian level this will be e. g. the relative torsion. But
we extend this approach to the more general case of Clifford structures, genera\-lized
Dirac operators and associated relative zeta functions, in particular we define relative
determinants which play, as well known, a big role in partition functions in QFT.

We give in this and the next section a brief outline of the uniform structure approach
necessary for the sequel and start with uniform structures of open manifolds. A complete
Riemannian manifold $(M^n,g)$ defines a proper metric space. Therefore we started in 
[6] with uniform structures of proper metric space, i. g. the Gromov--Hausdorff uniform
structure or series of Lipschitz uniform structures. The general image for Lipschitz
uniform structures is as follows. Let ${\mathfrak M}$ be the set of (isometry classes
of) proper metric spaces. Denote for a Lipschitz map
$\Phi: X = (X,d_X) \lra (Y,d_Y) = Y$, \\
$\dil \Phi = \sup\limits_{x_1 \neq x_2} d(\Phi x_1, \Phi x_2) / d(x_1,x_2)$
and define
\[
d_L(X,Y) = \left\{ \begin{array}{l} 
\inf \{ \max \{ 0, \log \dil \Phi \} + \max \{ 0, \log \dil \Phi \}
+ \sup\limits_x d(\Psi \Phi x,x) \\
+ \sup\limits_y d(\Phi \Psi y, y) \en | \en
\Phi : X \lra Y, \Psi: Y \lra X \mbox{ Lipschitz maps } \} \\
\mbox{ if }
\{ \dots \} \neq \emptyset \mbox{ and } \inf \{ \dots \} < \infty
\\[2ex]
\infty \mbox{ in the other case.} 
\end{array} \right.
\]
$X \begin{array}{c} {} \\[-1ex] \sim \\[-1ex] L \end{array} Y$
if $d_L (X,Y) = 0$ is an equivalence relation, 
${\mathfrak M}_L := {\mathfrak M}/\begin{array}{c} {} \\[-1ex] \sim \\[-1ex] L \end{array}$.
Set for $\delta > 0$
\[ V_\delta = \{ (X,Y) \in {\frak M}^2_L | d_L(X,Y) < \delta \} . \]

\begin{lemma}
${\mathfrak L} = \{ V_\delta \}_{\delta > 0}$ is a basis for a metrizable uniform structure
${\mathfrak U}_L ({\mathfrak M}_L)$.
\qed
\end{lemma}

Let $\overline{{\mathfrak M}_L}$ be the corresponding completion.

\begin{theorem}
{\bf a)} $\overline{{\mathfrak M}_L} = {\mathfrak M}_L$.

{\bf b)} ${\mathfrak M}_L$ is locally arcwise connected, hence components coincide with
arc components.

{\bf c)} ${\mathfrak M}_L = \sum\limits_{i \in I} \comp (X_i)$ as topological sum.

{\bf d)} $\comp (X) = \{ Y \in {\mathfrak M}_L | d_L(X,Y) < \infty \}$.
\qed
\end{theorem}

These definitions and considerations do not take into account the smooth structure and the
Riemannian metric if $(X,d_X) = (X^n,g,d_g)$ is a Riemannian manifold. If we restrict to
Riemannian manifolds we restrict the class of admitted maps and replace $d_L$ by sharper
distances. First we start with $\id_m$ as admitted maps and define the Sobolev uniform
structures for the Riemannian metrics on a fixed manifold $M$ and for the Clifford
connections on a Clifford bundle $(E,h) \lra (M,g)$ without fixed connection $\nabla^h$.

Denote by ${\cal M} (I,B_k)$ the set of all metrics $g$ on $M$ satisfying the
conditions $(I)$ and $(B_k)$.

Let $1 \le p \le \infty, k \ge r \ge \frac{n}{p}+2, \delta > 0$ and set
with $C(n,\delta) = 1 + \delta + \delta \sqrt{2n(n-1)}$
\begin {eqnarray} 
   V_\delta &=& \Bigg\{ (g,g') \in M(I,B_k)^2 \Big| \en  
   C(n,\delta)^{-1} g \le g' \le C(n,\delta) g
   \en \mbox{and} \en
   |g-g'|_{g,p,r} \nonumber \\
   &:=& \left( \int \left( |g-g'|^p_{g,x} +
   \sum^{r-1}_{i=0} | (\nabla^g)^i (\nabla^g-\nabla^{g'}) |^p_{g,x}
   \right) \,\, dvol_x(g)
   \right)^\frac{1}{p} < \delta
   \Bigg\} . \nonumber
\end{eqnarray}

Here $g,g'$ quasi isometric means $C_1 \cdot g \le g' \le C_2 \cdot g$
in the sense of quadratic forms. This is equivalent to 
$ {}^b|g-g'|_g < \infty $ and ${}^b|g-g'|_{g'} < \infty $, 
where for a tensor 
$t$ ${}^b|t|_g = \sup_{x \in M} |t|_{g,x}$.

\begin{prop}
Assume $p,k,r$ as above. Then ${\cal B} = \{ V_\delta \}_{\delta>0}$
is a basis for a metrizable uniform structure 
${\mathfrak U}^{p,r}(M(I,B_k))$.
\end{prop}

We refer to [12] for a proof. The key to the proof is the module
structure theorem. 
\qed

If we would replace in the definition of $V_\delta$
$|g-g'|_{g,p,r}$ above by
$|g-g'|_{g,p,r}= \left( \int \left( \sum^{r}_{i=0} | (\nabla^g)^i 
(g-g') |^p_{g,x} \right) \, dvol_x(g) \right)^\frac{1}{p}$ 
then we would get an equivalent uniform structure.

Let ${\cal M}^p_r(I,B_k) = {\cal M}(I,B_k)$ endowed with the uniform topology.
${\cal M}^{p,r} := \overline{{\cal M}^p_r}$ the completion. If 
$k \ge r > \frac{n}{p}+1$ then ${\cal M}^{p,r}$ still consists of $C^1$--metrics,
i.e. does not contain semi definite elements. This has been proved by
Salomonsen in [20].

\begin{theorem}
Let $k \ge r > \frac{n}{p}+2, g \in {\cal M}(I,B_k)$, 
$U^{p,r}(g) = \Big\{ g' \in {\cal M}^{p,r}(I,B_k) \big|$ ${}^b|g-g'|_g < \infty, \en
{}^b|g-g'|'_g < \infty$ and $|g-g'|_{g,p,r} < \infty \Big\}$
and denote by $comp(g) \subset {\cal M}^{p,r}(I,B_k)$ the component of 
$g$ in ${\cal M}^{p,r}(I,B_k)$. Then 
\be
comp(g) = U^{p,r}(g) ,
\end{equation}
$comp(g)$ is a Banach manifold, for $p=2$ a Hilbert manifold and
${\cal M}^{p,r}(I,B_k)$ has a representation as topological sum
\be
{\cal M}^{p,r}(I,B_k) = \sum_{j \in J} comp(g_j),
\end{equation}
J an uncountable set.
\end{theorem}

The proof is performed in [12]. \hfill $\Box$

\begin{remarks}

{\rm
{\bf 1.} If $M^n$ is compact then $J$ consists of one element.

{\bf 2.} All metrics in the completed space are at least of class 
$C^2$. Hence curvature is well defined. \hfill $\Box$
}
\end{remarks}

Let $(E,h) \rightarrow (M^n,g)$ be a Clifford bundle without a fixed
connection, $(M^n,g)$ with $(I)$ and $(B_k)$. 

Set

$C_E(B_k) = \Big\{ \nabla \Big| \nabla$ is Clifford connection, metric with
respect to $h$ and satisfies 

\hspace{2cm}$(B_k(E,\nabla)) \Big\}$

Assume $(E,h) \rightarrow (M^n,g)$ as above, $k \ge r > \frac{n}{p}+2,
\delta >0$ and set 
\begin{eqnarray} 
  V_\delta &=& \Big\{ (\nabla, \nabla') \in C_E(B_k)^2 \Big| 
  |\nabla-\nabla'|_{\nabla,p,r} \nonumber \\
  &:=& \Big( \int \sum^r_{i=0} |\nabla (\nabla - \nabla')|^p_x \,\,
  dvol_x(g) \Big)^{\frac{1}{p}} < \delta \Big\} . \nonumber
\end{eqnarray}

\begin{prop}
Assume $p,k,r$ as above. Then ${\cal B} = \{V_\delta \}_{\delta>0}$
is a basis for a metrizable uniform structure ${\mathfrak U}^{p,r}(C_E(B_k))$.
\end{prop}

We refer to [10] for a proof. \hfill $\Box$

Let $C^{p,r}_E(B_k)$ be the completion of $(C_E(B_k), 
{\mathfrak U}^{p,r}(C_E(B_k)))$. If $\nabla, \nabla' \in C^{p,r}_E(B_k)$ then
$\nabla-\nabla'$ is a 1--form $-\eta$ with values in ${\cal G}_E$
satisfying
\be
\eta_x(Y \cdot \Phi) = Y \cdot \eta_x(\Phi) .
\end{equation}

As well known, a metric connection $\nabla$ in $E$ induces a connection
$\nabla$ in ${\cal G}_E$. Denote
\begin{eqnarray}
   \Omega^1({\cal G}^{\cl}_E) &:=& \big\{ \eta \in \Omega^1 ({\cal G}_E)
   \Big| \eta \en \mbox{satisfies (3.3)} \en \big\} , \nonumber \\
   \Omega^{1,p}_r({\cal G}^{\cl}_E,\nabla) &:=&
   \Big\{ \eta \in \Omega^1 ({\cal G}^{\cl}_E)
   \Big| \nonumber \\
   &{}& |\eta|_{\nabla,p,r} 
   := \Big( \int \sum^r_{i=0} |\nabla^i \eta|^p_x \,\,
   dvol_x(g) \Big)^{\frac{1}{p}} < \infty \Big\} , \nonumber \\
   \Omega^{1,p,r}({\cal G}^{\cl}_E,\nabla) &:=&
   {\overline{\Omega^{1,p}({\cal G}^{\cl}_E,\nabla)}}^{| |_{\nabla,p,r}} . \nonumber
\end{eqnarray}

If $(M^n,g)$ satisfies $(I), (B_k)$ then
\be     
\hspace*{-0.2cm}
\Omega^{1,p,r}({\cal G}^{\cl}_E,\nabla) :=     
{\overline{C^\infty_c ({\cal G}_E)}}^{| |_{\nabla,p,r}} =
\Big\{ \eta \en \mbox{distributional} \en  \Big| |\eta|_{\nabla,p,r} 
< \infty \Big\} ,
\end{equation}     

\hspace*{-0.2cm} $ r \le k+2 $.

\begin{theorem}
Assume $(E,h) \rightarrow (M^n,g), p, k, r$ as above. Denote for
$\nabla \in C_E(B_k)$ by $comp(\nabla) \subset C^{p,r}_E(B_k)$ the component
of $\nabla$ in $C^{p,r}_E(B_k)$. Then
\be
comp(\nabla) = \nabla + \Omega^{1,p,r,} ({\cal G}^{\cl}_E, \nabla)
\end{equation}     
and $C^{p,r}_E(B_k)$ has a representation as topological sum 
\be     
C^{p,r}_E(B_k) = \sum_{j \in J} comp(\nabla_j) .
\end{equation}
\end{theorem}

The proof is performed in [10]. \hfill $\Box$

\begin{remarks}
    
{\rm     
{\bf 1.} If $M^n$ is compact then $C^{p,r}_E(B_k) = C^{p,r}_E$ consists
of one component
         
{\bf 2.} If $\nabla$ is non smooth then one sets $\nabla^i=(\nabla_0+(\nabla-\nabla_0))^i$,
$\nabla_0 \in comp{\nabla} \cap C_E(B_k)$, and the right hand side makes sense.
         
{\bf 3.} All connections in the complete space are at least of class $C^2$.
Hence curvature is well defined.
\qed
}
\end{remarks}

For the sequel, we must sharpen our considerations concerning Sobolev 
spaces. Let $(E,h,\nabla) \rightarrow (M^n,g)$ be a Riemannian vector
bundle. The connection $\nabla$ enters into the definition of the Sobolev
spaces $W^{p,r}$. Hence we should write $W^{p,r}(E,\nabla)$. Now there 
arises the natural question, how do the spaces $W^{p,r}(E,\nabla)$ depend on 
$\nabla$? We present here one answer. Other considerations are performed
in [9], [10].

\begin{prop}                                        
Let $(E,h,\nabla) \rightarrow (M^n,g)$ be a Riemannian vector
bundle with (I), $(B_k)$, $(B_k(E,\nabla))$, $k \ge r > \frac{n}{p}+1,
1 \le p < \infty$. Suppose $\nabla' \in comp(\nabla) \subset C^{p,r}_E(B_k)$,
$\nabla'$ smooth, i. e. $\nabla'=\nabla+\eta, \eta \in 
\Omega^{1,p,r}({\cal G}_E,\nabla) \cap C^\infty $.
Then              
\be               
W^{p,i}(E,\nabla) = W^{p,i}(E,\nabla'), \en 0 \le i \le r,
\end{equation}              
as equivalent Banach spaces.
\end{prop}

For the proof we refer to [9], [10]. The proof includes some combinatorial
considerations and essentially uses the module structure theorem. This
is the reason why we assumed $k \ge r > \frac{n}{p}+1$. But this
assumption can be weakened. We only need the validity of the module
structure theorem. \hfill $\Box$

\begin{remark} 
{\rm
The assumption $\eta$ smooth in superfluous. As we 
mentioned already several times, we can define $W^{p,r}(E,\nabla')$
and prove (3.7) for $\nabla' \in comp(\nabla)$ only. \hfill $\Box$
}
\end{remark}

\begin{coro}
Suppose $(E,h,\nabla) \rightarrow (M^n,g)$ as above, 
$k \ge r >\frac{n}{p}+2$, $\nabla' = \nabla + \eta$, 
$\eta \in \Omega^{1,p,r}({\cal G}_E,\nabla)$. Then
\be                                    
W^{2p,i}(E,\nabla) = W^{2p,i}(E,\nabla'), \en
0 \le i \le \frac{r}{2} .
\end{equation}
\end{coro}

\begin{proof}
$r>\frac{n}{p}+2$ implies 
$r-\frac{n}{p} \le \frac{r}{2} - \frac{n}{2p} , 2p \ge p , 
r \ge \frac{r}{2}$, i. e. 
\[ \Omega^{1,p,r}({\cal G}_E,\nabla) \subseteq 
   \Omega^{1,2p,\frac{r}{2}}({\cal G}_E,\nabla). \]
Now we apply the proof of 3.9 replacing $p \rightarrow 2 p, r \rightarrow
\frac{r}{2}$. \hfill $\Box$
\end{proof}

\begin{coro} 
Suppose $(E,h\nabla) \rightarrow (M^n,g)$ a Clifford bundle with the
conditions above for $p=1$, i. e. $k \ge r > n+2$ , 
$\nabla' \in comp(\nabla) \subset C^{1,r}_E(B_k)$, $\nabla'$ smooth. Then
\be
W^{2,i}(E,\nabla) \equiv W^i(E,\nabla) = W^i(E,\nabla') \equiv W^{2,i}(E,\nabla'),
\en 0 \le i \le \frac{r}{2} .
\end{equation}
\end{coro}

\begin{coro} 
Assume the hypothesises of 3.7. and write $D=D(\nabla,g), D'=D(\nabla',g)$.
Then
\be
H^i(E,D) = H^i(E,D'), \en 0 \le i \le \frac{r}{2} .
\end{equation}
In particular,
\be
{\cal D}_{D^i} = {\cal D}_{{D'}^i}, \en 0 \le i \le \frac{r}{2} ,
\end{equation}           
where ${\cal D}_{D^i}$ denotes the domain of definition of
$\overline{D^i}$ .
\end{coro}

\begin{proof}
(3.11) follows from the result of Chernoff that
$D^i$ is essentially self adjoint on $C^\infty_c(E)$ and 
${\cal D}_{D^i} = H^i(E,D)$. (3.10) follows from (3.9) and 2.3.

\qed
\end{proof}

Finally we make some remarks concerning the essential spectrum of
$D$ and $D^2$. More precisely, we prove that it is an invariant of 
$comp(\nabla)$. We have several distinct proofs for this and present
here a particularly simple one.

We consider Weyl sequences and restrict to orthonormal ones. Denote
by $\sigma_e(D)$ the essential spectrum of $D$. 
$\lambda \in \sigma_e(D)$ if and only if there exists a Weyl 
sequence for $\lambda$, i. e. an orthonormal sequence 
$(\Phi_\nu)_\nu, \Phi_\nu \in {\cal D}_D$; s. t. 
\be
\lim_{\nu \rightarrow \infty} (D-\lambda) \Phi_\nu = 0 .
\end{equation}

\begin{lemma} 
Suppose $\lambda \in \sigma_e(D)$. Then there exists a Weyl sequence
$(\Phi_\nu)_\nu$ for $\lambda$ s. t. for any compact subset
$K \subset M$
\be
\lim_{\nu \rightarrow \infty} |\Phi_\nu|_{L_2(K,E)} = 0 .
\end{equation}
\end{lemma}

This is Lemma 4.29 of [3]. One simply chooses an exhaustion 
$K_1 \subset K_2 \subset \dots , \bigcup K_i = M$, starts with an
arbitrary Weyl sequence $(\Psi_\nu)_\nu$, produces by the Rellich
lemma and a diagonal choice a subsequence $\raisebox{0.3ex}{$\chi$}_\nu$ such that
$(\raisebox{0.3ex}{$\chi$}_\nu)_\nu$ converges on any $K_i$ in the $L_2$-sense and defines
$\Phi_\nu := (\raisebox{0.3ex}{$\chi$}_{2\nu+1}-\raisebox{0.3ex}{$\chi$}_{2\nu}) | \sqrt{2}$. 
$\Phi_\nu)_\nu$ has the desired properties. \hfill $\Box$

\begin{prop} 
Suppose $(E,h,\nabla) \rightarrow (M^n,g)$ a Clifford bundle with (I),
$(B_k(M))$, $(B_k(E,\nabla))$, $k \ge r > n+2, n \ge 2$,
$\nabla' \in comp(\nabla) \subset C^{1,r}_E(B_k)$,
$D=D(\nabla,g), D'=D(\nabla',g)$. 
Then
\be
\sigma_e(D) = \sigma_e(D') .
\end{equation}
\end{prop}

\begin{proof}
\[ D'= \sum_i e_i \nabla'_{e_i} = \sum_i e_i \cdot (\nabla_{e_i} + 
   \eta_{e_i}(\cdot)) 
   = D + \eta^{op} ,\]
where the operator $\eta^{op}$ acts as
\[ \eta^{op}(\Phi)|_x = \sum_i e_i \cdot \eta_{e_i}(\Phi)|_x . \]
Then, pointwise, $|\eta^{op}|_x \le C \cdot |\eta|_x$,
$C$ independent of $x$. Given $\varepsilon > 0$, there exists a compact set
$K=K(\varepsilon) \subset M$ such that
\be                                                        
\sup_{x \in M \setminus K} |\eta|_x < \frac{\varepsilon}{C}, \en 
i. e. \en \sup_{x \in M \setminus K} |\eta^{op}|_x < \varepsilon .    
\end{equation}
Assume now $\lambda \in \sigma_e(D), (\Phi_\nu)_\nu$ a Weyl sequence
as in (3.13). According to (3.11), $\Phi_\nu \in {\cal D}_{D'}$.
Then 
\[ (D'-\lambda) \Phi_\nu = (D'-D) \Phi_\nu + (D-\lambda) \Phi_\nu . \]
By assumption, $(D-\lambda)\Phi_\nu \rightarrow 0$. Moreover,
\[ |(D'-D)\Phi_\nu|_{L_2(M,E)} =
   |\eta^{op}\Phi_\nu|_{L_2(M,E)} \le
   C \cdot (|\eta\Phi_\nu|_{L_2(K,E)} + 
   |\eta\Phi_\nu|_{L_2(M \setminus K,E)}) . \]
$|\eta \Phi_\nu|_{L_2(K,E)} \rightarrow 0$ and
\[ C \cdot |\eta \Phi_\nu|_{L_2(M \setminus K,E)} \le
   C \cdot \sup_{x \in M \setminus K} |\eta|_x  \cdot
   |\Phi_\nu|_{L_2(M \setminus K,E)} < \varepsilon . \]
Hence
$(D'-\lambda)\Phi_\nu \rightarrow 0$,
$\lambda \in \sigma_e(D')$,
$\sigma_e(D) \subseteq \sigma_e(D')$.
Exchanging the role of $D,D'$, we obtain 
$\sigma_e(D') \subseteq \sigma_e(D)$.
\qed
\end{proof}

Until now we fixed the fibre metric $h$ and allowed variation of the metric
connection $\nabla^h$. The next step is to admit simultaneous variation of
$g,h$ and $\nabla^h$, i. e. we consider Riemannian vector bundles 
$(E,h,\nabla^h) \lra (M^n,g)$ $\mbox{rk} E = N$ with $(I), (B_k(g)), (B_k(\nabla^h))$, 
$\nabla^h$ a metric connection for $h$ and admit variation of $g, h, \nabla^h$
in this class. Denote by $\geom (E \ra M, k)$ this set. Let
$k \ge r > \frac{n}{p}+2$, $1 \le p < \infty$, $\delta > 0$ and set
\bea
V_\delta &=& \Bigg\{ ((g,h,\nabla^h), (g'm,h',\nabla^{h'})) \in \geom (E \ra M,k)^2 \en | \non \\
&& C \left( n,\frac{\delta}{3} \right)^{-1} g \le g' \le C(n,\delta) g, \en
C \left( N,\frac{\delta}{3} \right)^{-1} h \le h' \le C \left( N,\frac{\delta}{3} h \right), \non \\
&& |g-g'|_{g,p,r} < \frac{\delta}{3}, \en |h-h'|_{g,h,\nabla^h,p,r} < \frac{\delta}{3}, \non \\
&& \left. |\nabla^h-\nabla^{h'}|_{g,h,\nabla^h,p,r} < \frac{\delta}{3} \right\}.
\eea
Here
\bea
|h-h'|_{g,h,\nabla^h,p,r} &=& \left( \int \sum^r_{i=0} |(\nabla^h)^i(h-h')|^p_{g,h,x}
\,\, dvol_x(g) \right)^{\frac{1}{p}}, \non \\
|\nabla^h-\nabla^{h'}|_{g,h,\nabla^h,p,r} &=& \left( \int \sum^r_{i=0} |(\nabla^h)^i
(\nabla^h-\nabla^{h'})|^p_{g,h,x} \,\, dvol_x(g) \right)^{\frac{1}{p}}. \non
\eea
We remark that $\nabla^h-\nabla^{h'}$ is still tensorial although 
$\nabla^h, \nabla^{h'}$ belong to different fibre metrics.

\begin{prop}
${\mathfrak L} = \{ V_\delta \}_{\delta > 0}$ ist a basis for a metrizable uniform
structure ${\mathfrak U}^{p,r} (\geom (E \ra M, k))$.
\end{prop}

\bi
\no
{\bf Sketch of proof.}
For the symmetry we have to prove that
$|g-g'|_{g,p,r} < \frac{\delta}{3}$, 
$|h-h'|_{g,h,\nabla^h,p,r} < \frac{\delta}{3}$, 
$|\nabla^h-\nabla^{h'}|_{g,h,\nabla^h,p,r} < \frac{\delta}{3}$ 
imply the existence of $\delta'=\delta'(\delta)>0$ s. t., replacing in all three
expressions $| \en |_{g,h,\nabla^h,p,r}$ by $| \en |_{g',h',\nabla^{h'},p,r}$,
these are $< \frac{\delta'}{3}$ and $\delta'(\delta) \underset{\delta \ra 0} {\lra} 0$.
For the first expression we proved this in [12]. For the second and third expression
this is quite analogous. We have to perform a simple induction and we indicate only the
initial steps. The main point are the combinatorial formulas in [10], [12].
Nevertheless, there are some deviations from [10], [12] since we have here
simultaneous variation of $g,h,\nabla^h$. Assume
\be 
|g-g'|_{g,p,r} < \frac{\delta}{3}, \quad 
|h-h'|_{g,h,\nabla^h,p,r} < \frac{\delta}{3}, \quad
|\nabla^h-\nabla^{h'}|_{g,h,\nabla^h,p,r} < \frac{\delta}{3} 
\end{equation}
and consider
\be
|(\nabla^{h'}(h-h')|_{g',h',x} \quad \mbox{and} \quad
|(\nabla^{h'})^i(\nabla^h-\nabla^{h'})|_{g',h',x} 
\end{equation}
For $i=0$, according to (3.16) with $\nabla=\nabla^h$, $\nabla'=\nabla^{h'}$
\bea
&& C_1 |h-h'|_{h,x} \le |h-h'|_{h',x} \le C_2 |h-h'|_{h,x} \\
&& C_3 |\nabla-\nabla'|_{g,h,x} \le |\nabla-\nabla'|_{g',h',x} \le C_4 |\nabla-\nabla'|_{g,h,x}
\eea
and similarly for the $L_p$--norm. (3.20) is a consequence of the fact that quasiisometry
of the metrics $g$ and $g'$, $h$ and $h'$ induce quasiisometry of pointwise norms of tensors
and of volume elements. Let $i=1$. Then, omitting the indices at the pointwise norms, 
\bea
&& |\nabla'(h-h')| \le |(\nabla'-\nabla)(h-h')| + |\nabla (h-h')|, \\
&& |\nabla'(\nabla-\nabla')| \le |(\nabla'-\nabla)(\nabla-\nabla')| + |\nabla (\nabla-\nabla')|
\eea
The first terms on the r. h. s. are in $L_p$ according to (3.19), (3.20) and the module
structure theorem. The second term is in $L_p$ by the assumption (3.17). All $L_p$--norms
are controlled by (3.17), i. e. by a polynomial in $\frac{\delta}{3}$ with positive coefficients
and without constant terms. For $i=2$. ${\nabla'}^2 = \nabla' \circ \nabla'$.
The left hand $\nabla'$ on the right hand side is of the form 
$\nabla^{g'} \otimes \nabla^{h'}$ and this can be rewritten as
\be
(\nabla^{g'} - \nabla^g) \otimes (\nabla^{h'} - \nabla^h) + (\nabla^{g'} - \nabla^g)
\otimes \nabla^h + \nabla^g \otimes (\nabla^{h'} - \nabla^h) + \nabla^g \otimes \nabla^h.
\end{equation}
Moreover, 
\bea
{\nabla'}^2 &=& \nabla' (\nabla' - \nabla) + \nabla' \nabla \non \\
&=& (\nabla' - \nabla) (\nabla' - \nabla) + \nabla (\nabla' - \nabla)
+ (\nabla' - \nabla) \nabla + \nabla^2.
\eea
Using this and inserting into the left hand $(\nabla')'$s the expression (3.23), into
the left hand $(\nabla)'$s $\nabla = \nabla^g \otimes \nabla^h$, we obtain from the
assumption (3.17), the settled case $i=1$ and from the module structure theorem
\bea
&& |{\nabla'}^2 (h-h')|_{g',h',\nabla^{h'},p,r} < \frac{\delta_2}{3} , \\
&& |{\nabla'}^2 (\nabla-\nabla')|_{g',h',\nabla^{h'},p,r} < \frac{\delta_2}{3} ,
\eea
where $\delta_2$ is a polynomial in $\frac{\delta}{3}$ with positive coefficients
and without constant term, i. e. $\delta_2(\delta) \underset{\delta \ra 0} {\lra} 0$.
For higher $i$ we apply the combinatorial considerations performed in [10], [12]
and obtain that 
$|g-g'|_{g,p,r}$, 
$|h-h'|_{g,h,\nabla^h,p,r}$, 
$|\nabla^h-\nabla^{h'}|_{g,h,\nabla^h,p,r}$ 
are bounded by polynomials in 
$|\nabla^i(g-g')|_{g,p}$, 
$|(\nabla^h)^j(h-h')|_{g,p}$, 
$|(\nabla^h)^k(\nabla^h-\nabla^{h'})|_{g,p}$
$0 \le i,j,k \le r$
without constant term and with positive coefficients. This indicates the proof for the
symmetry of ${\mathfrak L} = \{ V_\delta \}_{\delta>0}$. The proof of transitivity
is a little more difficult and completely modeled in [10], [12]. 
\qed

\bi
Let $\geom^p_r(E \ra M,k)$ be the pair
$(\geom(E \ra M,k), {\mathfrak U}^{p,r})$
and denote by $\geom^{p,r}(E \ra M,k)$ the completion.

\begin{prop}
Denote by $\arccomp^{p,r}(g,h,\nabla^h)$ the arc component of $(g,h,\nabla^h)$
in $\geom^{p,r}(E \ra M, k)$. If 
$(g',h',\nabla^{h'}) \in \arccomp^{p,r}(g,h,\nabla^h)$
then
\[ |g-g'|_{g,p,r} < \infty, \quad 
   |h-h'|_{g,h,\nabla^h,p,r} < \infty, \quad
   |\nabla^h - \nabla^{h'}|_{g,h,\nabla^h,p,r} < \infty. \]
\qed
\end{prop}

We are not able to prove that $\geom^{p,r}(E \ra M,k)$
is locally arcwise connected. The reason for this is 
indicated after proposition 3.19.

Until now we defined our uniform structures for fixed $M$ or fixed $E \ra M$,
respectively. It is possible and for many applications in classification theory
necessary to include maps into the definition of the uniform structures. We define,
as we did in ${\mathfrak U}_L$. This is an extensive task. We did some work in this
direction in [6] and refer to this paper. For our main theorems only these uniform
structures are of interest. One we introduce now, the others in the next section. 
First we recall some facts from mapping theory.

Let $(M^n,g)$,
$(N^{n'},h)$ be open, satisfying $(I)$ and $(B_k)$ and let
$f \in C^\infty(M,N)$. Then the differential $df = f_* = Tf$ is a section
of $T^*M \otimes f^*TN$. $f^*TN$ is endowed with the induced connection
$f^* \nabla^h$. The connections $\nabla^g$ and $f^*\nabla^h$ induce
connections $\nabla$ in all tensor bundles
$T^q_s(M) \otimes f^*T^u_v(N)$. Therefore, $\nabla^mdf$ is well defined.
Assume $m \le k$. We denote by $C^{\infty,m}(M,N)$ the set of all
$f \in C^{\infty}(M,N)$ satisfying
\[ ^{b,m} |df| = \sum\limits^{m-1}_{i=0} \sup\limits_{x \in M}
    |\nabla^i df|_x < \infty . \]

Consider complete manifolds $(M^n,g)$ $({M'}^n,g')$ and 
$C^{\infty,m}(M,M')$. A diffeomorphism $f : M \lra M'$ will be
called $m$--bibounded if $f \in C^{\infty,m}(M,M')$ and
$f^{-1} \in C^{\infty,m}(M',M)$. Sufficient for this is
1. $f$ is a diffeomorphism, 2. $f \in C^{\infty,m}(M,M')$,
3. $\inf\limits_x | \lambda |_{\min} (df)_x > 0$.

Starting with bibounded diffeomorphisms of $C^{\infty,m+1}$, $m+1 \ge r+1$,
one can define a Sobolev uniform structure ${\mathfrak U}^{p,r}({\cal D}^{m+1}(M))$,
thus getting a completed Sobolev diffeomorphism group ${\cal D}^{p,r}(M,g)$.
This procedure is really very complicated and the content of [9].
A bibounded $C^1$ diffeomorphism $f$ is in ${\cal D}^{p,r}(M,g)$ 
if and only if $f=\exp X \circ \tilde{f}$, where $\tilde{f}$
is smooth and $(r+1)$--bibounded, $X$ a vector field anlong $\tilde{f}$,
$X \in \Omega^{p,r}(\tilde{f}^*TM)$ small. Since we mostly do not
consider completed versions in our manifold classification approach,
we restrict to a subgroup $\tilde{{\cal D}}^{p,r}(M,g)$. Here
$f \in \tilde{{\cal D}}^{p,r}$ if $f=\exp X \circ \tilde{f}$, $\tilde{f}$
$(r+1)$--bibounded and $X$ a smooth vector field along $\tilde{f}$
satisfying $|X|_{p,r} < \infty$. Recall
$|X|_{p,r} = \left( \int \sum\limits^r_{i=0} |(\tilde{f}^* 
\nabla)^i X|^p_{g,x} \,\, dvol_x(g) \right)^{\frac{1}{p}}$.

We define as generalization of $\tilde{{\cal D}}^{p,r}(M,g)$ the spaces
$\tilde{{\cal D}}^{p,r}((M_1,g_1),(M_2,g_2))$ 
and
$C^{\infty,p,r}((M_1,g_1),(M_2,g_2))$. Here 
$f \in \tilde{{\cal D}}^{p,r}(M_1,M_2)$ if
$f = \exp X \circ \tilde{f}$. 
$\tilde{f} : M_1 \begin{array}{c} \cong \\[-2ex] \lra \\[-2ex] {} \end{array} M_2$
an $(r+1)$--bibounded smooth diffeomorphism, $X$ a smooth vector field along
$\tilde{f}$ with $|X|_{p,r}$ small; $f \in C^{\infty,p,r}(M_1,M_2)$
if $f = \exp X \circ \tilde{f}$. $f \in X^{\infty,r+1}(M_1,M_2)$ and $X$ a
smooth vector field along $f$ with $|X|_{p,r}$ small. Small means in all cases
$\sup\limits_x |X|_{x} < r_{inj}$, which can be assured by a sufficiently small
Sobolev norm.

Consider pairs $(M^n_1,g_1), (M^n_2,g_2) \in {\mathfrak M}^n (mf, I, B_k)$ with
the following pro\-per\-ty.

There exist compact submanifolds $K^n_1 \subset M^n_1$, $K^n_2 \subset M^n_2$ and an 
$f \in \tilde{{\cal D}}^{p,r} (M^n_1 \setminus K_1$, $M^n_2 \setminus K_2)$.

For such pairs define
\bea
 d^{p,r}_{L,diff,rel} ((M_1,g_1), (M_2,g_2)) 
 & := & \inf \Big\{ \max \{ 0, \log {}^b |df| \}
 + \max \{ 0, \log {}^b |dh| \} \non \\
 & + & \sup\limits_{x \in M_1} \dist (x, hfx)
 + \sup\limits_{y \in M_2} \dist (y, fhy)  \non \\
 & + & |(f|_{M_1 \setminus K_1})^* g_2-g_1|_{M_1 \setminus K_1}|_{g,p,r} \non \\
 && \Big| \en f \in C^{\infty,p,r} (M_1,M_2), h \in C^{\infty,p,r} (M_1,M_2) \non \\
 && \mbox{ and for some } K_1 \subset M \mbox{ holds } \non \\
 && f|_{M_1 \setminus K_1} \in \tilde{{\cal D}}^{p,r} (M_1 \setminus K_1, f(M_1 \setminus K_1))
 \non \\
 && \mbox{ and } h|_{f(M_1 \setminus K_1)} = (f|_{M_1 \setminus K_1})^{-1}
 \Big\}. \non
\eea
if $\{ \dots \} \neq \emptyset$ and $\inf \{ \dots \} < \infty$. In the other case
set $d^{p,r}_{L,diff,rel} ((M_1,g_1), (M_2,g_2)) = \infty$.

Set
\bea
 V_\delta & = & \{ ((M_1,g_1),(M_2,g_2)) \in ({\mathfrak M}^n (mf,I,B_k))^2
 \non \\
 && | \en d^{p,r}_{L,diff,rel} ((M_1,g_1),(M_2,g_2)) < \delta \} .
 \non
\eea

\begin{prop}
${\mathfrak B} = \{ V_\delta \}_{\delta > 0}$ is a basis for a metrizable uniform
structure on ${\mathfrak M}^n(mf,I,B_k)/_\sim$, where $(M_1,g_1) \sim (M_2,g_2)$
if $d^{p,r}_{L,diff,rel} ((M_1,g_1),(M_2,g_2)) = 0$.
\qed
\end{prop}

Denote the corresponding uniform structure with ${\mathfrak U}^{p,r}_{L,diff,rel}$
and ${\mathfrak M}^{n,p,r}_{L,diff,rel}$ for ${\mathfrak M}^n(mf,I,B_k)$ endowed 
with this uniform structure.

It follows again from the definition that 
$d^{p,r}_{L,diff,rel}((M_1,g_1),(M_2,g_2)) < \infty$ implies 
$d_L((M_1,g_1),(M_2,g_2)) < \infty$. Hence $(M_2,g_2) \in \comp_L(M_1,g_1)$,
i. e. 
\[ \left\{ (M_2,g_2) \in {\mathfrak M}^{n,p,r}_{L,diff,rel} \en | \en 
   d^{p,r}_{L,diff,rel} (M_1,M_2) < \infty \right\} \subseteq \comp_L (M_1,g_1). \]
For this reason we denote the left hand side $\{ \dots \}$ by
$\comp^{p,r}_{L,diff,rel} (M_1,g_1) = \{ \dots \} = \{ \dots \} \cap \comp_L (M_1,g_1)$
keeping in mind that this is not an arc component but a subset (of manifolds) of
a Lipschitz arc component, endowed with the induced topology.

We extend all this to Riemannian vector bundles $(E,h,\nabla^h) \lra (M^n,g)$ of 
bounded geometry. First we have to define ${\cal D}^{p,r}(E \ra M)$. For this we
consider the total space $E$ as open Riemannian manifold of bounded geometry with
respect to the Kaluza--Klein metric and restrict the uniform structure of [9]
to bundle maps $f=(f_E,f_M)$. Quite similar we define for 
$E_i = ((E_i, h_i, \nabla^{h_i}) \lra (M^n_i, g_i))$, $i=1,2$,
${\cal D}^{p,r}(E_1,E_2)$ by corresponding bundle isomorphisms and
$\tilde{{\cal D}}^{p,r}(E_1,E_2)$ as the smooth generating elements,
i. e. $f=(f_E,f_M) \in \tilde{{\cal D}}^{p,r}(E_1,E_2)$
if and only if $f_E = \exp_{\tilde{f}_E} X \circ \tilde{f}_E$,
$\tilde{f}_E \in C^{\infty,r} (E_1,E_2)$ a diffeomorphism,
$\tilde{f}^{-1}_E \in C^{\infty,r} (E_1,E_2)$, 
$X \in \Omega^p_r (f^*TE_2)$, $f_E$ a bundle isomorphism, similarly for
$f_M : M_1 \lra M_2$. Quite analogously to ${\mathfrak M}^n (mf,I,B_k)$
above we denote the bundle isometry classes of Riemannian vector bundles
$(E,h,\nabla) \lra (M^n,g)$ with $(I), (B_k(g)), (B_k(\nabla))$ of
$rk N$ over $n$--manifolds by ${\cal B}^{N,n}(I,B_k)$. Set for
$k \ge r > \frac{n}{p}+2$, 
$E_i=((E_i,h_i,\nabla^{h_i}) \lra (M^n_i,g_i)) \in {\cal B}^{N,n}(I,B_k)$,
$i=1,2$
\bea
d^{p,r}_{L,diff} (E_1,E_2) &=& \inf \{ \max \{ 0, \log {}^b|df_E| \}
+ \max \{ 0, \log {}^b|df^{-1}_E| \} 
+ \non \\
&& + \max \{ 0, \log {}^b|df_M| \} + \max \{ 0, \log {}^b|df^{-1}_M| \} 
+ \non \\
&& + |g_1 - f^*_M g_2|_{g_1,p,r} + |h_1 - f^*_E h_2|_{g_1,h_1,\nabla^{h_1},p,r}
+ \non \\
&& + |\nabla^{h_1} - f^*_E \nabla^{h_2}|_{g_1,h_1,\nabla^{h_1},p,r} \en | 
\non \\
&& f=(f_E,f_M) \in {\cal D}^{p,r}(E_1,E_2) \} \non
\eea
if $\{ \dots \} \neq \emptyset$ and $\inf \{ \dots \} < \infty$.
In the other case set $d^{p,r}_{L,diff} (E_1,E_2) = \infty$.
Here we remark that ${}^b|df_E|, {}^b|df^{-1}_E|, {}^b|df_M|, {}^b|df^{-1}_M| < \infty$
automatically imply the quasi isometry of $h_1, f^*_E h_2$ or $g_1,f^*_M g_2$, 
respectively. A simple consideration and a variant of the proof of transitivity
in 3.19 show that $d(E_1,E_2)=0$ is an equivalence relation $\sim$. Set
${\cal B}^{N,n}_{L,diff} (I,B_k) := {\cal B}^{N,n} (I,B_k)/\sim$
and for $\delta>0$
\[ V_\delta = \{ (E_1,E_2) \in ({\cal B}^{N,n}_{L,diff} (I,B_k))^2 \} \en | \en
   d^{p,r}_{L,diff} (E_1,E_2) < \delta . \]

\begin{prop}
${\mathfrak L} = \{ V_\delta \}_{\delta>0}$ is a basis for a metrizable uniform
structure ${\mathfrak U}^{p,r}_{L,diff}$.
\end{prop}

The proof is quite analogous to that of 5.16 in [6]. 
\qed

Denote ${\cal B}^{N,n,p}_{L,diff,r} (I,B_k)$ for the pair
$({\cal B}^{N,n} (I,B_k), {\mathfrak U}^{p,r}_{L,diff})$ and
${\cal B}^{N,n,p,r}_{L,diff} (I,B_k)$ for the completion.

The next task would be to prove the locally arcwise connectedness of
${\cal B}^{N,n,p,r}_{L,diff}$. If we restrict to 
$(E,h) \lra (M^n,g)$, i. e. we forget the metric connection $\nabla^k$,
then the corresponding space is locally arcwise connected according to
5.19 of [6]. Taking into account the metric connection $\nabla^h$,
the situation becomes much worse. Given $(g,h,\nabla^h), (g',h',\nabla^{h'})$
sufficiently neighbored, we have to prove that they could be connected
by a (sufficiently short) arc $\{ (g_t, h_t, \nabla^{h_t}) \}$.
Here $\nabla^{h_t}$ must be metric w. r. t. $h_t$. We were not able to construct
the arc $\{ \nabla^{h_t} \}_t$ for given $\{ h_t \}_t$. One could also try
to set $\nabla^t=(1-t)\nabla+t\nabla$ and to construct $h_t$ from $\nabla^t$
s. t. $\nabla^t$ is metric w. r. t. $h_t$. In local bases $e_1, \dots, e_n$,
$\Phi_1, \dots, \Phi_N$ this would lead to the system
\[ \nabla^t_{e_i} h_{t,\alpha \beta} = \Gamma^\gamma_{t,i\alpha} h_{t,\gamma \beta}
   + \Gamma^\gamma_{t, i \beta} h_{t,\alpha \gamma}, \quad i=1, \dots, n, 
   \alpha, \beta = 1, \dots, N , \]
where $h_{t,\alpha \beta} = h_t (\Phi_\alpha, \Phi_\beta)$, 
$\nabla^t_i \Phi_\alpha = \Gamma^\gamma_{t,i \alpha} \Phi_\gamma$.
This is a system of $n\frac{N(N+1)}{2}$ equations for the $\frac{N(N+1)}{2}$
components $h_{\alpha \beta}$, i. e. it is overdetermined. With other words,
we don't see a comparatively simple and natural proof for locally arcwise
connectedness. ${\cal B}^{N,n,p,r}_{L,diff} (I,B_k)$ is a complete metric space.
Hence locally and locally arcwise connectedness coincide. But to prove locally
connectedness amounts very soon to similar questions just discussed.

Consider for $E=((E,h,\nabla^h) \lra (M,g)) \in {\cal B}^{N,n} (I,B_k)$
\be
\{ E' \in {\cal B}^{N,n,p,r}_{L,diff} (I,B_k) \en | \en 
   d^{p,r}_{L,diff} (E,E') < \infty \} .
\end{equation}
The set is open and contains the arccomponent of $E$. If 
${\cal B}^{N,n,p,r}_{L,diff} (I,B_k)$ would be locally arcwise connected =
locally connected then we would have
\be
\arccomp (E) = \comp (E) = (3.27).
\end{equation}
If we endow the total spaces $E$ with the Kaluza--Klein metric
\[ g_E(X,Y) = h(X^V,Y^V) + g_M(\pi_*X,\pi_*Y), \quad X^V, Y^V 
   \mbox{ vertical components} , \]
then $(E,g_E)$ becomes a Riemannian manifold of bounded geometry, hence a proper
metric space. It follows from the definition that
\be
\{ E' \in {\cal B}^{N,n,p,r}_{L,diff} (I,B_k) \en | \en 
   d^{p,r}_{L,diff} (E,E') < \infty \} \subseteq \comp_L (E).
\end{equation}
(3.28), (3.29) and the foregoing considerations are for us motivation enough to
define the generalized component $\gencomp (E)$ by
\be
\gencomp(E) := \{ E' \in {\cal B}^{N,n,p,r}_{L,diff} (I,B_k) \en | \en
d^{p,r}_{L,diff}(E,E') < \infty \} .
\end{equation}
In particular $\gencomp (E)$ is a subset of a Lipschitz component and is endowed
with a well defined topology coming from ${\mathfrak U}^{p,r}_{L,diff}$.

The final step in this section consists in the additional admission of compact
topological perturbations. We consider pairs 
$E_i = ((E_i,h_i,\nabla^{h_i}) \lra (M^n_i, g_i)) \in {\cal B}^{N,n} (I,B_k)$, 
$i=1,2$, with the following property. There exist compact submanifolds
$K^n_i \subset M^n_i$ and $f=(f_E,f_M)$, 
$f|_{M_1 \setminus K_1} \in \tilde{\cal D}^{p,r+1} (E_1|_{M_1 \setminus K_1},
E_2|_{M_2 \setminus K_2}) \cap C^{\infty,r+1} (E_1,E_2)$.
For such pairs define
\bea
d^{p,r}_{L,diff,rel} (E_1,E_2) &=& \inf \{ \max \{ 0, \log {}^b|df_E| \}
+ \max \{ 0, \log {}^b|dh_E| \} 
+ \non \\
&& + \max \{ 0, \log {}^b|df_M| \} + \max \{ 0, \log {}^b|dh_M| \} 
+ \non \\
&& + \sup\limits_{e_1} d(h_Ef_E e_1,e_1) + \sup\limits_{e_2} d(f_Eh_E e_2,e_2) +
\non \\
&& + \sup\limits_{x_1} d(h_Mf_M x_1,x_1) + \sup\limits_{x_2} d(f_Mh_M x_2,x_2) +
\non \\
&& + |(f_M|_{M_1 \setminus K_1})^* g_2-g_1|_{M_1 \setminus K_1}|_{g_1,p,r} +
\non \\
&& + |(f_E|_{E}|_{M_1 \setminus K_1})^* h_2-h_1|_{E_1}|_{M_1 \setminus K_1}
|_{g_1,h_1,\nabla^{h_1},p,r}
+ \non \\
&& + |(f_E|_{E}|_{M_1 \setminus K_1})^* \nabla^{h_2}-\nabla^{h_1}|_{E_1}|_{M_1 \setminus K_1}
|_{g_1,h_1,\nabla^{h_1},p,r}
\non \\
&& | \en f=(f_E,f_M) \in C^{\infty,r+1}(E_1,E_2), \non \\
&& h=(h_E,h_M) \in C^{\infty,r}(E_1,E_2) \} \non \\
&& \mbox{bundle maps and for some } K_1 \subset M_1 \mbox{ holds} \non \\
&& f|_{E_1}|_{M_1 \setminus K_1} \in \tilde{\cal D}^{p,r} (E_1|_{M_1 \setminus K_1},
E_2|_{f_M(M_1 \setminus K_1)}) \mbox{ and} \non \\
&& h|_{E_2}|_{f({M_1 \setminus K_1})} = (f|_E|_{M_1 \setminus K_1})^{-1} \}
\eea
if $\{ \dots \} \neq \emptyset$ and $\inf \{ \dots \} < \infty$. In the other case set
$d^{p,r}_{L,diff,rel} (E_1,E_2) = \infty$. This definition seems to be quite lengthy
but it is quite natural. It measures outside a compact set the distinction of ''shape''
and the geometric objects in question. Set
\[ V_\delta = \{ (E_1,E_2) \in ({\cal B}^{N,n} (I,B_k))^2 \en | \en
   d^{p,r}_{L,diff,rel} (E_1,E_2) < \delta \}. \]

\begin{prop}
${\mathfrak L} = \{ V_\delta \}_{\delta>0}$ is a basis for a metrizable uniform
structure ${\mathfrak U}^{p,r}_{L,diff,rel}$ on ${\cal B}^{N,n}(I,B_k)/\sim$ where
$E_1 \sim E_2$ if $d^{p,r}_{L,diff,rel} (E_1,E_2) = 0$.
\qed
\end{prop}

Denote ${\cal B}^{N,n,p,r}_{L,diff,rel} (I,B_k)$ for the completed
${\cal B}^{N,n} (I,B_k)$ endowed with this uniform structure. We have again that
$d^{p,r}_{L,diff,rel} (E_1,E_2) < \infty$ implies $d_L(E_1,E_2) < \infty$ 
(here we consider $E_1,E_2$ as proper metric spaces). Hence
$E_2 \in \comp_L (E_1)$, i. e.
\be
\{ E_2 \in {\cal B}^{N,n,p,r}_{L,diff,rel} (I,B_k) \en | \en d^{p,r}_{L,diff,rel}
(E_1,E_2) < \infty \} \subseteq \comp_L (E_1).
\end{equation}
For this reason we denote again the left hand side $\{ \dots \}$ of (3.32) by
$\gencomp^{p,r}_{L,diff,rel} (E_1)$ keeping in mind that this is not an arc
component but a subset of a Lipschitz component endowed with the induceed from
${\mathfrak U}^{p,r}_{L,diff,rel}$ topology.

In our later applications we prove and thereafter use the trace class property of
$e^{-tD^2}-e^{-t{D'}^2}$. Here essentially enter estimates for $D-D'$, coming from
the explicit expression for $D-D'$. But in this expression only enter
$g,\nabla,\cdot,g',\nabla',\cdot'$. This is the reason that we consider in some of our
applications smaller generalized components, which are in fact arc components. 
Exactly spoken, we restrict in some of our applications to those uniform structures
and components where $h_1=f^*_Eh_2$, i. e. the fibre metric does not vary.
Nevertheless, the generalized components play the more important role as appropriate
equivalence classes in classification theory. We prove the trace class property
of $e^{-tD^2}-e^{-t{D'}^2}$ also for generalized components.

Set now
\bea
d^{p,r}_{L,diff,F} (E_1,E_2) &=& \inf \{ \max \{ 0, \log {}^b|df_M| \}
+ \max \{ 0, \log {}^b|df^{-1}_M| \} + \non \\
&& |g_1-f^*_M g_2|_{g_1,p,r} + |\nabla^{h_1}-f^*_M \nabla^{h_2}|_{g_1,h_1,\nabla^{h_1},p,r}
\non \\
&& | \en f=(f_E,f_M) \in \tilde{\cal D}^{p,r}(E_1,E_2), \non \\
&& f_E \mbox{ fibrewise an isometry } \}
\eea
if $\{ \dots \} \neq \emptyset$ and $\inf \{ \dots \} < \infty$. In the other case set
$d^{p,r}_{L,diff,F} (E_1,E_2) = \infty$. $d^{p,r}_{L,diff,F} (\cdot,\cdot) = 0$
is an equivalence relation $\sim$. Set 
${\cal B}^{N,n}_{L,diff,F} (I,B_k) = {\cal B}^{N,n} (I,B_k)/\sim$ and for
$\delta > 0$
\[ V_\delta = \{ (E_1,E_2) \in ({\cal B}^{N,n}_{L,diff,F} (I,B_k))^2 \en | \en
   d^{p,r}_{L,diff,F} (E_1,E_2) < \delta \}. \]

\begin{prop}
${\mathfrak L} = \{ V_\delta \}_{\delta>0}$ is a basis for a metrizable uniform
structure ${\mathfrak U}^{p,r}_{L,diff,F}$.
\qed
\end{prop}

Denote ${\cal B}^{N,n,p,r}_{L,diff,F} (I,B_k)$ for the corresponding completion.

\begin{prop}
{\bf a)} ${\cal B}^{N,n,p,r}_{L,diff,F} (I,B_k)$ is locally arcwise connected.

{\bf b)} In ${\cal B}^{N,n,p,r}_{L,diff,F} (I,B_k)$ coincide components with 
arccomponents.

{\bf c)} ${\cal B}^{N,n,p,r}_{L,diff,F} (I,B_k) = \sum\limits_{i \in I}
\comp^{p,r}_{L,diff,F} (E_i)$ as topological sum.

{\bf d)} For $E \in {\cal B}^{N,n}$
\[ \comp^{p,r}_{L,diff,F} (E) = \{ E' \in {\cal B}^{N,n,p,r}_{L,diff,F} (I,B_k) 
   \en | \en d^{p,r}_{L,diff,F} (E,E') < \infty \} . \]
\end{prop}

\begin{proof}
of a)
$g_t=(1-t)g_1+tf^*_Mg_2$, $\nabla_t=(1-t)\nabla^{h_1}+tf^*_E\nabla^{h_2}$
yield an arc between $E_1$ and $f^*E_2$. Here we use $h_1=f^*_Eh_2$ and that
$\nabla^{h_1}, \nabla^{h_2}$ are metric.

\qed
\end{proof}

Quite analogously we define - based on 
$h_1|_{E|_{M \setminus K}} = f^*_E (h_2|_{E'|_{M' \setminus K'}})$
the uniform space
${\cal B}^{N,n,p,r}_{L,diff,F,rel}(I,B_k)$
and its generalized components
\[ \gencomp^{p,r}_{L,diff,F,rel} (E) = \{ E' \en | \en d^{p,r}_{L,diff,F,rel}
   (E,E') < \infty \} . \]
Here $d^{p,r}_{L,diff,F,rel}(E,E')$ is defined as
$d^{p,r}_{L,diff,F,rel}(\cdot,\cdot)$ with the additional condition
$h|_E|_{M \setminus K} = f^*_E (h'|_{E'}|_{M' \setminus K'}$.
${\cal B}^{N,n,p,r}_{L,diff,F,rel}$ is not locally arcwise connected.
Now the classification of ${\cal B}^{N,n}(I,B_k)$ amounts to two tasks.

1. Classification (i. e. ''counting'') the (generalized) components 

$\gencomp(E)$ by invariants,

2. Classification of the elements inside a component by invariants, where
number valued invariants should be relative invariants.

Until now $g,h,\nabla^h$ could be fixed independently, keeping in mind that
$\nabla^h$ should be a metric connection with respect to $h$. The situation
rapidly changes if we restrict to Clifford bundles. The new ingredient is the
Clifford multiplication $\cdot$ which relates $g,h,\nabla^h$. This will be
studied in the next section.

\section{Uniform structures of Clifford bundles}

\setcounter{equation}{0}

As we know from the definition, a Clifford bundle
$(E,h,\nabla^h,\cdot) \lra (M^n,g)$
has as additional ingredient the module structure of 
$E_m$ over $\cl_m(g) = \cl(T_mM,g_m)$.
Change of $g$, $g \lra g'$, changes point by point the Clifford algebra,
$\cl_m(g) \lra \cl_m(g')$. Locally they are isomorphic since by radial
parallel transport of orthonormal bases in a normal neighborhood 
$U(m_0)$ always
\be
\cl(M,g)|_{U(m_0)} \cong U(m_0) \times \cl(\R^n) \cong 
U(m_0) \times \cl(M,g')|_{U(m_0)} .
\end{equation}
The same holds for bundles of Clifford modules if we fix the typical fibre, i. e.
\be
E|_U \cong E'|_U
\end{equation}
but globally (4.2) in general does not hold although as vector bundles $E$ and $E'$
can be isomorphic. The point is the module structure which includes $g$ (in 
$\cl_m(g)$) as operating algebra and $\cdot : T_m M \otimes E_m \lra E_m$. Therefore
we consider for a moment the following admitted deformations. Let
$(E,h,\nabla^h,\cdot) \lra (M^n,g)$ be a Clifford bundle of $rk N$. The vector
bundle structure $E \lra M$ (of $rk N$) shall remain fixed. We admit variation of
$g$, hence of $\cl(TM)$, variation of $\cdot \in \homm(TM \otimes E,E)$, hence 
variation of the structure of $E$ as bundle of Clifford modules, compatible
variation of $h, \nabla^h$. $\homm (TM \otimes E, E) \cong T^*M \otimes E^* \otimes E$
is for given $g, h$ a Riemannian vector bundle. Including $\nabla^h$, the notion of
Sobolev sections is well defined. For fixed $g, h, \nabla^h$ the space
$\Gamma (\mult,g,h,\nabla^h)$ of Clifford multiplications $\cdot$ is a well defined
subspace of $\Gamma ( \homm (TM \otimes E, E))$ described invariantly by the 
conditions
\bea
\sk{X \cdot \Phi}{\Psi} &=& - \sk{\Phi}{\Psi \cdot Y} , \\
\nabla_X (Y \cdot \Phi) &=& (\nabla_X Y) \cdot \Phi + Y \cdot \nabla_Y \Phi .
\eea
We describe this space locally as follows. Let $U(m_0) \subset M$ and
$e_1(m), \dots, e_n(m)$ be a field of orthonormal bases obtained from 
$e_1(m_0), \dots, e_n(m_0)$ by radial parallel translation, similarly
$\Phi_1, \dots, \Phi_N$ a field of orthonormal bases in $E|_{U(m_0)}$, 
obtained also by radial parallel translation. Then (4.3), (4.4) mean for 
the attachment $e_i \otimes \Phi_j \lra e_i \cdot \Phi_j$
\bea
\sk{e_i \cdot \Phi_j}{\Phi_k} &=& - \sk{\Phi_j}{e_i \cdot \Phi_k}, \non \\
&& i=1,\dots,n, \en j=1, \dots, N, \\
\nabla_{e_i} (e_j \cdot \Phi_k) &=& (\nabla_{e_i} e_j) \cdot \Phi_k
+e_i \cdot \nabla_{e_j} \Phi_k, \non \\
&& i,j=1,\dots,n, \en k=1,\dots,N.
\eea
If we write in the linear space $\homm (T_m \otimes E_m, E_m)$
$e_i \cdot \Phi_j = \sum\limits^N_{k=1} a^k_{ij} \Phi_k$ then (4.5) reduces to
$\frac{nN(N+1)}{2}$ independent linear equations between the
$a^k_{ij}$, $a^k_{ij}=-a^j_{ik}$. Hence the fibre $\mult_m (g,h,\nabla^h)$
is an $n \cdot N^2 - \frac{nN(N+1)}{2} = \frac{nN}{2}(N-1)$--dimensional
affine subspace at any point $m$.
This establishes a locally trivial fibre bundle $\mult(g,h,\nabla^h)$.
The charts for local trivialization arise from radial parallel translation
$P$ of $\cdot$ from $m$ to $m_0$ since
$\nabla_{\frac{\partial}{\partial r}} (P(e_i \cdot \Phi_j) 
- (Pe_i) \cdot (P \Phi_j) = 0$.
\clm $(g,h,\nabla^h) \subset \Gamma (\mult(g,h,\nabla_h))$
now are those sections of $\mult(g,h,\nabla^h)$
which additionally satisfy (4.4) or (4.6).
Consider the Riemannian vector bundle
$({\cal H} = \homm (TM \otimes E, E) = (T^* M \otimes E^* \otimes E, h_{\cal H})
\lra (M^n,g))$. Assume as always $(M^n,g)$ with $(I), (B_k)$, 
$(E,h,\nabla^h)$ with $(B_k)$, $k \ge \frac{n}{2}+2$ and
$(\pi: E \lra M) \in C^{\infty,k+1}(E,M)$.
Then with respect to the Kaluza--Klein metrix 
$g_{\cal H}(X,Y) = h_{\cal H} (X^V,Y^V) + g_M (\pi_* X, \pi_* Y)$,
$X^V$, $Y^V$ vertical components, the total space of ${\cal H}$ becomes
a Riemannian manifold satisfying $(I)$, $(B_k)$. The fibres ${\cal H}_m$
are totally geodesic submanifolds, moreover they are flat. The latter also
holds for the affine fibres $\mult_m (g,h,\nabla^h)$ of $\mult (g,h,\nabla^h)$.
If $\cdot$ and $\cdot'$ are sections of $\mult (g,h,\nabla^h)$ satisfying
(4.4), i. e. $\cdot, \cdot' \in \clm (g,h,\nabla^h)$ then
\be
(1-t) \cdot + t \cdot' \in \clm (g,h,\nabla^h),
\end{equation}
Let $k \ge r > \frac{n}{p}+2$, $\delta > 0$. Set
\[ V_\delta = \{ (\cdot,\cdot') \in \clm (g,h,\nabla^h)^2 \en | \en
   |\cdot - \cdot'|_{g,h,\nabla^h,p,r} < \delta \} . \]

\begin{lemma}
${\mathfrak L} = \{V_\delta\}_{\delta>0}$ is a basis for a metrizable uniform
structure ${\mathfrak U}^{p,r}(\clm (g,h,\nabla^h))$. 
\qed
\end{lemma}

Denote by $\clm^{p,r} (g,h,\nabla^h)$ the completion.

\begin{prop}
{\bf a)} $\clm^{p,r} (g,h,\nabla^h)$ is locally arcwise connected.

{\bf b)} In $\clm^{p,r} (g,h,\nabla^h)$ coincide components and arc components.

{\bf c)} $\clm^{p,r} (g,h,\nabla^h) = \sum\limits_{i \in I} \comp^{p,r}(\cdot_i)$.

{\bf d)} $\comp^{p,r}(\cdot) = \{ \cdot' \en | \en |\cdot - \cdot'|_{g,h,\nabla^h,p,r}
< \infty \}$.
\end{prop}

\begin{proof}
a) follows from (4.7), b) and c) follow from a), d) is a simple calculation.
\qed
\end{proof}

\begin{remark}
{\rm
In the language of the intrinsic Riemannian geometry of $\mult(g,h,\nabla^h)$ and 
of $\Gamma (\mult(g,h,\nabla^h))$ we can rewrite $|\cdot - \cdot'|_{g,h,\nabla^h,p,r} < \delta$
as $\cdot' = \exp X \circ \cdot$, 
$X \in \Gamma ( T ( \mult (g,h,\nabla^h)))$,
$|X|_{g,h,\nabla^h,p,r} < \delta$.
Here $X_m$ lies in the affine subspace $\mult_m$.
\qed
}
\end{remark}

Denote by $\cl {\cal B}^{N,n} (I,B_k)$ the set of (Clifford isometry classes) of all
Clifford bundles $(E,h,\nabla^h,\cdot) \lra (M^n,g)$ of (module) rank $N$ over 
$n$--manifolds, all with $(I)$ and $(B_k)$.

\begin{lemma}
Let $E_i=((E_i,h_i,\nabla^{h_i},\cdot_i) \lra (M^n_i,g_i)) \in \cl {\cal B}^{N,n}(I,B_k)$,
$i=1,2$ and $f=(f_E,f_M) \in \tilde{\cal D}^{p,r+1} (E_1,E_2) \cap C^{\infty,k+1} (E_1,E_2)$
be a vector bundle isomorphism between bundles of Clifford modules, 
$f_E(X \cdot_1 \Phi) = (f_M)_* X \cdot_2 f_E \Phi$. Then
$f^*E_2 := ((E_1, f^*_E h_2, f^*_E \nabla^{h_2}, f^*_E \cdot_2) \lra (M_1, f^*_M g_2)) \in
\cl {\cal B}^{N,n}(I, B_k)$.
\end{lemma}

\begin{proof}
The definitions of $f^*_E h_2$, $f^*_E \nabla^{h_2}$, $f^*_M g_2$ are clear.
$f^*_E \cdot_2$ is defined by $X (f^*_E \cdot_2) \Phi = f^{-1}_E (f_* X \cdot_2 f_E \Phi)$.
It is now an easy calculation that $f^*E_2 \in \cl {\cal B}^{N,n}(I,B_k)$.
\qed
\end{proof}

Let $k \ge r > \frac{n}{p}+2$ and define for $E_1,E_2 \in \cl {\cal B}^{N,n}(I,B_k)$
\bea
d^{p,r}_{L,diff} (E_1,E_2) &=& \inf \{ \max \{ 0, \log {}^b|df_E| \}
+ \max \{ 0, \log {}^b|df^{-1}_E| \} + \non \\
&& \max \{ 0, \log {}^b|df_M| \} + \max \{ 0, \log {}^b|df^{-1}_M| \} + \non \\
&& |g_1-f^*_M g_2|_{g_1,p,r} + |h_1-f^*_E h_2|_{g_1,h_1,\nabla^{h_1},p,r} + \non \\
&& |\nabla^{h_1}-f^*_E \nabla^{h_2}|_{g_1,h_1,\nabla^{h_1},p,r} + 
|\cdot_1-f^*_E \cdot_2|_{g_1,h_1,\nabla^{h_1},p,r} 
\non \\
&& | \en f=(f_E,f_M) \in \tilde{\cal D}^{p,r}(E_1,E_2) \mbox{ is a } (k+1) - 
\mbox{bounded} \non \\
&& \mbox{isomorphism of Clifford bundles} \} \non
\eea
if $\{ \dots \} \neq \emptyset$ and $\inf \{ \dots \} < \infty$. In the other case set
$d^{p,r}_{L,diff} (E_1,E_2) = \infty$. $d^{p,r}_{L,diff}$
is numerically not symmetric but nevertheless it defines a uniform
structure which is by definition symmetric. Set for $\delta > 0$
\[ V_\delta = \{ (E_1,E_2) \in \cl{\cal B}^{N,n} (I,B_k))^2 \} \en | \en
   d^{p,r}_{L,diff} (E_1,E_2) < \delta \} . \]

\begin{prop}
${\mathfrak L} = \{ V_\delta \}_{\delta>0}$ is a basis for a metrizable uniform
structure ${\mathfrak U}^{p,r}_{L,diff,F} (\cl {\cal B}^{N,n}(I,B_k))$.
\qed
\end{prop}

The proof is quite analogous to that of 5.16 in [6] using additionally 4.1. 

\qed

Denote $\cl {\cal B}^{N,n,p}_{L,diff,r}(I,B_k)$ for the pair 
$(\cl {\cal B}^{N,n}(I,B_k), {\cal U}^{p,r})$ and
$\cl {\cal B}^{N,n,p,r}_{L,diff}(I,B_k)$ for the completion.
By the same motivation as above and at the end of section 3 we 
introduce again the generalized component
$\gencomp (E) = \gencomp^{p,r}_{L,diff} ((E,h,\nabla^h) \lra (M,g))
\subset \cl {\cal B}^{N,n,p,r}_{L,diff}(I,B_k)$ by
\be
\gencomp^{p,r}_{L,diff} (E) = \{ E' \in \cl {\cal B}^{N,n,p,r}_{L,diff}(I,B_k) \en | \en
d^{p,r}_{L,diff} (E,E') < \infty \}.
\end{equation}
$\gencomp (E)$ contains $\arccomp(E)$ and is endowed with a Sobolev topology
induced from ${\cal U}^{p,r}_{L,diff}$.

The absolutely last step in our uniform structures approach is the additional
admission of compact topological pertulations. We proceed as in (3.31), i. e.
we assume additionally 
$E_i =((E_i,h_i,\nabla^{h_i},\cdot_i) \lra (M^n_i, g_i)) \in \cl {\cal B}^{N,n}(I,B_k)$, 
add in (3.31) still
$|(f_E|_{E|_{M_1 \setminus K_1}})^* \cdot_2 - \cdot_1 |_{E|_{M_1 \setminus K_1}}
|_{g_1,h_1,\nabla^{h_1},p,r}$
and assume 
$f=(f_E,f_M)|_{M_1 \setminus K_1}$, 
$h=(h_E,h_M)|_{M_2 \setminus K_2 = f_M (M_1 \setminus K_1)}$
vector bundle isomorphisms (not necessary Clifford isometric).
Then we get 
$d^{p,r}_{L,diff,rel}(E_1,E_2)$,
define $V_\delta$, ${\mathfrak L} = \{ v_\delta \}_{\delta > 0}$,
obtain the metrizable uniform structure
${\mathfrak U}^{p,r}_{L,diff,rel} (\cl {\cal B}^{N,n}(I,B_k))$
and finally the completion
$\cl {\cal B}^{N,n,p,r}_{L,diff,rel}$.
We set again
$\gencomp (E) = \gencomp^{p,r}_{L,diff,rel} (E) = \{ E' \in \cl 
{\cal B}^{N,n,p,r}_{L,diff,rel}(I,B_k)) \en | \en d^{p,r}_{L,diff,rel} (E,E') < \infty \}$
which contains the arc component and inherits a Sobolev topology from
${\mathfrak U}^{p,r}_{L,diff,rel}$.

As in section 3, we obtain by requiring additionally 
$h_1 = f^*_E h_2$ or
$h_1|_{E_1|_{M_1 \setminus K_1}} = f^*_E (h_2|_{E_2|_{M_2 \setminus K_2}})$
local distances
$d^{p,r}_{L,diff,F} (\cdot, \cdot)$ or
$d^{p,r}_{L,diff,F,rel} (\cdot, \cdot)$ 
and corresponding uniform spaces
$\cl {\cal B}^{N,n,p,r}_{L,diff,F} (I,B_k)$ or
$\cl {\cal B}^{N,n,p,r}_{L,diff,F,rel} (I,B_k)$
respectively. We obtain generalized components
\be
\gencomp^{p,r}_{L,diff,F} (E)
\end{equation}
and
\be
\gencomp^{p,r}_{L,diff,F,rel} (E)
\end{equation}
as before. One of our main technical results 
will be that $E$ and $E'$ in the same generalized
component (4.10) or (4.11) respectively implies that after transforming
$e^{-t{D'}^2}$ into the Hilbert space $L_2((M,E),g,h)$,
$e^{-tD^2} - e^{-t{D'}^2}$ and $e^{-tD^2} D - e^{-t{D'}^2} D'$ 
are of trace class and their trace norm is uniformly bounded on compact
$t$--Intervalls $[a_0,a_1]$, $a_0>0$. For our later applications the components
(4.9), (4.10) are most important, excluded one case, the case
$D^2 = \Delta (g)$, ${D'}^2 = \nabla (g')$.
In this case variation of $g$ automatically induces variation of 
the fibre metric and we have to consider (4.8) and 
$\gencomp^{p,r}_{L,diff,rel}(E)$.

\section{General heat kernel estimates}

\setcounter{equation}{0}

We collect some standard facts concerning the heat kernel of 
$e^{-tD^2}$. The best references for this are [3], [4].

We consider the self-adjoint closure of $D$ in $L_2(E)=H^0(E)$,
$D=\int\limits^{+ \infty}_{- \infty} \lambda E_\lambda$.

\begin{lemma}
$\{e^{itD}\}_{t \in {\bf R}}$ defines a unitary group on the spaces
$H^r(E)$, for $0 \le h \le r$ holds
\be
|D^h e^{itD} \Psi|_{L_2} \cdot  |e^{itD} D^h \Psi|_{L_2} =
|D^h \Psi|_{L_2} .
\end{equation}
\qed
\end{lemma}

We can extend this action to $H^{-r}(E)$ by means of duality.

\begin{lemma} 
$e^{-tD^2}$ maps $L_2(E) \equiv H^0(E) \rightarrow H^r(E)$
for any $r>0$ and
\be
|e^{-tD^2}|_{L_2 \rightarrow H^r} \le C \cdot t^{- \frac{r}{2}}, \en
t \in ]0, \infty[, \en C=C(r) .
\end{equation}
\end{lemma}

\begin{proof} 
Insert into $e^{-tD^2} = \int e^{-t \lambda^2} \, dE_\lambda $
the equation
\[ e^{-t \lambda^2} = \frac{1}{\sqrt{4 \pi t}} 
   \int\limits^{+ \infty}_{- \infty} e^{i \lambda s} 
   e^{- \frac{s^2}{4t}} \,\, ds \]
and use
\[ \sup |\lambda^r e^{-t \lambda^2}| \le C \cdot t^{- \frac{r}{2}} . \]
\qed
\end{proof}

\begin{coro}
Let $r, s \in {\bf Z}$ be arbitrary. Then 
$e^{-tD^2}: H^r(E) \rightarrow H^s(E)$ continuously.
\end{coro}

\begin{proof} 
This follows from 5.2., duality and the semi group property
of $\{e^{-tD^2}\}_{t \ge 0}$. 
\qed
\end{proof}

$e^{-tD^2}$ has a Schwartz kernel 
$W \in \Gamma ({\bf R} \times M \times M, E \boxtimes E)$, 
\[ W(t,m,p) = \langle \delta(m), e^{-tD^2} \delta(p) \rangle ,\]
where $\delta (m) \in H^{-r}(E) \otimes E_m$ is the map
$\Psi \in H^r(E) \rightarrow \langle \delta(m), \Psi \rangle = \Psi(m)$,
$r> \frac{n}{2}$. The main result of this section is the fact that for
$t>0, W(t,m,p)$ is a smooth integral kernel in $L_2$ with good decay
properties if we assume bounded geometry.

Denote by $C(m)$ the best local Sobolev constant of the map 
$\Psi \rightarrow \Psi(m), r > \frac{n}{2}$, and by $\sigma(D^2)$
the spectrum.

\begin{lemma} 

{\bf a)} $W(t,m,p)$ is for $t>0$ smooth in all Variables.

{\bf b)} For any $T>0$ and sufficiently small $\varepsilon > 0$ there 
exists $C > 0$ such that
\be
|W(t,m,p)| \le e^{-(t-\varepsilon) \inf \sigma(D^2)} \cdot
C \cdot C(m) \cdot C(p) \en \mbox{for all} \en t \in ]T, \infty[ .
\end{equation}

{\bf c)} Similar estimates hold for 
$(D^i_m D^j_p W)(t,m,p)$. 
\end{lemma}

\begin{proof} 

a) First one shows $W$ is continuous, which follows from
$\langle \delta(m), \cdot \rangle$ continuous in $m$ and 
$e^{-tD^2} \delta(p)$ continuous in $t$ and $p$. Then one applies
elliptic regularity.

b) Write 

$|\langle \delta(m), e^{-tD^2} \delta(p) \rangle| =
|\langle (1+D^2)^{-\frac{r}{2}} \delta(m), 
(1+D^2)^r e^{-tD^2} (1+D^2)^{\frac{r}{2}} 
\delta(p) \rangle|$

c) Follows similar as b. 
\qed
\end{proof}

\begin{lemma} 
For any $\varepsilon>0, T>0, \delta>0$ there exists $C>0$ such that for 
$r>0, m \in M, T>t>0$ holds
\be
\int\limits_{M \setminus B_r(m)} |W(t,m,p)|^2 \,\, dp \le
C \cdot C(m) \cdot e^{-\frac{(r-\varepsilon)^2}{(4+\delta)t}} .
\end{equation}
A similar estimate holds for $D^i_m D^j_p W(t,m,p)$.  
\end{lemma}

We refer to [3] for the proof.
\qed

\begin{lemma} 
For any $\varepsilon>0, T>0, \delta>0$ there exists $C>0$ such that for 
all $m, p \in M$ with $dist(m,p) > 2 \varepsilon, T>t>0$ holds   
\be
|W(t,m,p)|^2  \le
C \cdot C(m) \cdot C(p) \cdot 
e^{-\frac{(dist(m,p)-\varepsilon)^2}{(4+\delta)t}} .
\end{equation}
A similar estimate holds for $D^i_m D^j_p W(t,m,p)$.  
\end{lemma}

We refer to [3] for the proof.
\qed

\begin{prop} 
Assume $(M^n,g)$ with (I) and $(B_K)$, $(E,\nabla)$ with $(B_K)$, 
$k \ge r > \frac{n}{2}+1$. Then all estimates in 5.4, 5.5, 5.6 hold
with uniform constants.
\end{prop}

\begin{proof} 
>From the assumptions $H^r(E) \cong W^r(E)$ and
$\sup_m C(m) = C =$ global Sobolev constant for $W^r(E)$,
according to 2.13.
\qed
\end{proof}

Let $U \subset M$ be precompact, open, $(M^+,g^+)$ closed with 
$U \subset M^+$ isometrically and $E^+ \rightarrow M^+$ a
Clifford bundle with $E^+|_U \cong E|_U$ isometrically. Denote by
$W^+(t,m,p)$ the heat kernel of $e^{-t{D^+}^2}$.

\begin{lemma}   
Assume $\varepsilon>0, T>0, \delta>0$. Then there exists $C>0$ such
that for all $T>t>0, m,p \in U$ with $B_{2\varepsilon}(m),
B_{2\varepsilon}(p) \subset U$ holds
\be 
|W(t,m,p)-W^+(t,m,p)| \le 
C \cdot e^{- \frac{\varepsilon^2}{(4+\delta)t}}
\end{equation}
\end{lemma}

We refer to [3] for the simple proof. \qed

\begin{coro}   
$tr W(t,m,m)$ has for $t \rightarrow 0^+$ the same asymptotic expansion
as for $tr W^+(t,m,m)$. 
\qed
\end{coro}

\section{Trace class property under variation of the Clifford connection}

\setcounter{equation}{0}

Our procedure in the next sections is as follows. We admit step by step larger
perturbations of a given Clifford structure $((E,h,\nabla^h) \lra (M^n,g))$
where the perturbations $E'=((E',h',\nabla^{h'}) \lra ({M'}^n,g'))$ are
elements of $\gencomp (E)$. Fixing $\gencomp(E)$ has the consequence that we
start with $(E,h,\nabla^h) \lra (M^n,g)$, fix the vector bundle $E \lra M$ and 
admit step by step perturbation of $\nabla^h$, then simultaneously of
$g,h,\nabla^h,\cdot$ and finally compact topological perturbations. The
reason for this is that if we have
$f=(f_E,f_M): ((E,h,\nabla^h) \lra (M^n,g)) \lra ((E',h',\nabla^{h'}) \lra
(M',g'))$ 
which is a vector bundle isomorphism (not necessarily an isometry) then 
we can write $f: E \lra E'$ as $f: E \lra f^*E' \overset{i} {\lra} E'$
where $i$ is a Clifford isometry, i. e. the interesting perturbation takes
place from $E$ to $f^*E'$. In this section we start with the simplest case,
perturbation of the metric Clifford connection $\nabla = \nabla^h$,
$\nabla \lra \nabla'$, $\nabla - \nabla' \in \Omega^{1,p,r} ({\mathfrak G}^{\cl}_E)$.
This case has model character for the other cases and for this reason we
present carefully the sometimes very complicated estimates.

We come now to the first main result of this paper.

\begin{theorem}
Assume $(E,\nabla) \rightarrow (M^n,g)$, $(M^n,g)$ with (I) and $(B_k)$  
$(E,\nabla)$ with $(B_k)$, $k \ge r > n+2, n \ge 2$, 
$\nabla' \in comp(\nabla) \cap C_E(B_k) \subset C^{1,r}_E(B_k)$,
$D=D(g, \nabla), D'=D'(g, \nabla')$ generalized Dirac operators.
Then
\[ e^{-tD^2} - e^{-t{D'}^2}\]
is for $t>0$ a trace class operator and its trace norm ist uniformly
bounded on compact $t$--intervalls $[a_0, a_1], a_0>0$.
\end{theorem}

\begin{remark}
{\rm 
The condition 
$\nabla' \in comp(\nabla) \cap C_E(B_k) \subset C^{1,r}_E(B_k)$, i. e. 
$\nabla' \in comp(\nabla)$ and additionally $\nabla'$ smooth and satisfying $(B_K)$
can be weakened to $\nabla' \in comp(\nabla) \subset C^{1,r}_E(B_k)$.
The main reason for this is that we can write
$\nabla'=\nabla'_0+(\nabla'-\nabla'_0)$, $\nabla'_0 \in  C_E(B_k)$, 
$|\nabla'-\nabla'_0|_{1,r,\nabla} < \varepsilon$. Then one can reestablish the 
whole Sobolev theory etc. extensively using the module structure theorem.
}
\end{remark}

We refer to the forthcomming paper [14]. 
\qed

The proof of theorem 6.1 will occupy the remaining part of this section.
We always assume the assumptions of 6.1. According to 3.13,
\[ {\cal D}_D = {\cal D}_{D'}, \qu {\cal D}_{D^2} = {\cal D}_{{D'}^2} . \]

\begin{lemma} 
Assume $t>0$. Then 
\be
e^{-tD^2} -  e^{-t{D'}^2} = \int\limits^t_0 e^{-sD^2} 
({D'}^2-D^2) e^{-(t-s){D'}^2} \,\, ds .         
\end{equation}                                              
\end{lemma}

\begin{proof} 
(6.1) means at heat kernel level
\be
W(t,m,p) - W'(t,m,p) = - \int\limits^t_0 \int\limits_M (W(s,m,q),
(D^2-{D'}^2) W'(t-s,q,p))_q \,\, dq \,\, ds ,
\end{equation}
where $( , )_q$ means the fibrewise scalar product at $q$ and 
$dq=dvol_q(g)$. Hence for (6.1) we have to prove (6.2). (6.2) is an
immediate consequence of Duhamel's principle. Only for completeness,
we present the proof of (6.2), which is the last of the following
7 facts and implications.

\me
{\bf 1.} For $t>0$ is $W(t,m,p) \in L_2(M,E,dp) \cap {\cal D}_D^2$

\me
{\bf 2.} If $\Phi, \Psi \in {\cal D}_D^2$ then
$\int (D^2 \Phi, \Psi) - (\Phi, D^2 \Psi) \,\, dvol = 0$ 
(Greens formula). 
         
\me
{\bf 3.} $
   ((D^2+ \frac{\partial}{\partial \tau}) \Phi(\tau, g) 
   \Psi(t-\tau, q))_q - (\Phi(\tau,g), (D^2+ \frac{\partial}{\partial t})
   \Psi(t-\tau, q))_q = $

  $ = (D^2 (\Phi(\tau,q), \Psi(t-\tau,q))_q - (\Phi(\tau,q),
   D^2 \Psi(t-\tau,q))_q + \frac{\partial}{\partial \tau}
   (\Phi(\tau,g), \nonumber $

  $ \Psi(t-\tau,q))_q$ .

\me
{\bf 4.}
$ \int\limits^\beta_\alpha \int\limits_M 
   ((D^2+\frac{\partial}{\partial \tau}) \Phi(\tau,q), 
   \Psi(t-\tau, q))_q - (\Phi(\tau,q),  
    (D^2+\frac{\partial}{\partial t}) \Psi(t-\tau,q))_q \,\, dq \,\, d\tau = 
   $

$  = \int\limits_M [(\Phi(\beta,q), \Psi(t-\beta,q)_q -
   (\Phi(\alpha,q), \Psi(t-\alpha,q))_q ] \,\, dq$.

\me 
{\bf 5.}
$ \Phi(t,q) = W(t,m,q), \Psi(t,q) = W'(t,q,p) $ yields 

$  - \int\limits^\beta_\alpha \int\limits_M (W(\tau,m,q), 
   (D^2+\frac{\partial}{\partial t}) W'(t-\tau,q,p) \,\, dq \,\, d\tau = $

$   = \int\limits_M [(W(\beta,m,q), W'(t-\beta,q,p))_q -
   (W(\alpha,m,q), W'(t-\alpha,q,p))_q] \,\, dq$  . 

\me
{\bf 6.}
Performing $\alpha \rightarrow 0^+, \beta \rightarrow t^-$ in 5. yields

$ - \int\limits^t_0 \int\limits_M (W(s,m,q),  
   (D^2+\frac{\partial}{\partial t}) W'(t-s,q,p))_q \,\, dq \,\, ds =
    W(t,m,p) - W'(t,m,p)$. 

\me
{\bf 7.}
Finally, using $D^2 + \frac{\partial}{\partial t} = 
D^2 - {D'}^2 + {D'}^2 + \frac{\partial}{\partial t} $ and
$({D'}^2+\frac{\partial}{\partial t}) W' = 0$ we obtain

$ W(t,m,p) - W'(t,m,p) = - \int\limits^t_0 \int\limits_M (W(s,m,q),
(D^2-{D'}^2) W'(t-s,q,p))_q \,\, dq \,\, ds $

which is (6.2).

If we write $D^2-{D'}^2 = D(D-D')+(D-D')D'$ then
\bea
e^{-tD^2}-e^{-t{D'}^2} 
&=& - \int\limits^t_0 e^{-sD^2} (D^2-{D'}^2) e^{-(t-s){D'}^2)} \,\, ds \non \\
&=& - \int\limits^t_0 e^{-sD^2} D(D-D') e^{-(t-s){D'}^2)} \,\, ds \non \\
&& - \int\limits^t_0 e^{-sD^2} (D-D')D' e^{-(t-s){D'}^2)} \,\, ds \non \\
&=& \int\limits^t_0 e^{-sD^2} D \eta e^{-(t-s){D'}^2)} \,\, ds \non \\
&& + \int\limits^t_0 e^{-sD^2} \eta D' e^{-(t-s){D'}^2)} \,\, ds , \non
\eea
where $\eta=\eta^{op}$ in the sense of section 3, 
$\eta^{op} (\Psi)|_x = \sum\limits^n_{i=1} e_i \eta_{e_i} (\Psi)$
and $|\eta^{op}|_{op,x} \le C \cdot |\eta|_x$, $C$ independent of $x$.
We split
$\int\limits^t_0 = \int\limits^{\frac{t}{2}}_0 + \int\limits^t_{\frac{t}{2}}$,
\def\theequation{$I_1$}   
\be
e^{-tD^2}-e^{-t{D'}^2} 
= \int\limits^{\frac{t}{2}}_0 e^{-sD^2} D \eta e^{-(t-s){D'}^2} \,\, ds 
\end{equation}
\def\theequation{$I_2$}   
\be
+ \int\limits^{\frac{t}{2}}_0 e^{-sD^2} \eta D' e^{-(t-s){D'}^2} \,\, ds 
\end{equation}
\def\theequation{$I_3$}   
\be 
+ \int\limits^t_{\frac{t}{2}} e^{-sD^2} D \eta e^{-(t-s){D'}^2} \,\, ds 
\end{equation}
\def\theequation{$I_4$}   
\be
+ \int\limits^t_{\frac{t}{2}} e^{-sD^2} \eta D' e^{-(t-s){D'}^2} \,\, ds.
\end{equation}
\def\theequation{\thesection.\arabic{equation}}  
\setcounter{equation}{2}

We want to show that each integral ${\rm (I_1) - (I_4)}$ is a product of
Hilbert--Schmidt operators and to estimate their Hilbert--Schmidt norm.
Consider the integrand of ${\rm (I_4)}$,
\be
(e^{-sD^2} \eta)(D' e^{-(t-s){D'}^2}).
\end{equation}
According to 5.2,
\bea
|e^{-(t-s){D'}^2}|_{L_2 \rightarrow H^1} & \le & C \cdot (t-s)^{-\frac{1}{2}} \\
|D' e^{-(t-s){D'}^2}|_{L_2 \rightarrow L_2} & \le & |D'|_{H^1 \rightarrow L_2}
\cdot |e^{-(t-s) {D'}^2}|_{L_2 \rightarrow H^1} \non \\
& \le & C \cdot (t-s)^{-\frac{1}{2}} .
\eea
Write
\be
(e^{-sD^2} \eta)(D' e^{-(t-s){D'}^2}) = (e^{-\frac{s}{2}D^2} f)(f^{-1}
e^{-\frac{s}{2}D^2} \eta)(D' e^{-(t-s){D'}^2}).
\end{equation}

Here $f$ shall be a scalar function which acts by multiplikation. (6.5) estimates the  
right hand factor in (6.6). The main point is the right choice of $f$.
$e^{-\frac{s}{2}D^2} f$ has the integral kernel
\be
W(\frac{s}{2},m,p) f(p)
\end{equation}
and $f^{-1} e^{-\frac{s}{2}D^2} \eta$ has the kernel
\be
f^{-1}(m) W(\frac{s}{2},m,p) \eta(p).
\end{equation}

We have to make a choice such that (6.7), (6.8) are square integrable
over $M \times m$ and that their $L_2$--norm is on compact $t$--intervals uniformly bounded.

We decompose the $L_2$--norm of (6.7) as
\begin{eqnarray}
  &{}& \int\limits_M \int\limits_M | W(\frac{s}{2},m,p) |^2 |f(m)|^2 \,\, dm \,\, dp = 
  \\
  &{}& \int\limits_M \int\limits_{dist(m,p) \ge c} | W(\frac{s}{2},m,p) |^2 
  |f(m)|^2 \,\, dp \,\, dm = \\
  &{}& \int\limits_M \int\limits_{dist(m,p) < c} | W(\frac{s}{2},m,p) |^2 
  |f(m)|^2 \,\, dp \,\, dm 
\end{eqnarray}
We obtain from 5.4 for $s \in ]\frac{t}{2},t[$
\[ (6.11) \le \int\limits_M C_1 |f(m)|^2 vol B_c(m) \,\, dm \le
   C_2 \int\limits_M |f(m)|^2 \,\, dm \]
and from 5.5.
\[ \int\limits_M \int\limits_{dist(m,p) \ge c} | W(\frac{s}{2},m,p) |^2  
   |f(m)|^2 \,\, dp \,\, dm  \le 
   \int\limits_M C_3 e^{- \frac{-(c-\varepsilon)^2}{4+\delta} \frac{s}{2}} 
   |f(m)|^2 \,\, dm \,\, \le \]
\be   
   \le C_3 \cdot e^{- \frac{-(c-\varepsilon)^2}{4+\delta} \frac{s}{2}} 
   \int\limits_M |f(m)|^2 \,\, dm, \quad c > \varepsilon .
\end{equation}
Hence the estimate of $\int\limits_M \int\limits_M
| W(\frac{s}{2},m,p) |^2 |f(m)|^2 dp dm    $   
for $s \in [\frac{t}{2},t]$ is done if 
\[ \int\limits_M |f(m)|^2 \,\, dm < \infty  \]
and then $|e^{-\frac{s}{2}D^2}f|_2 \le C_4 \cdot |f|_{L_2}$, where, according to
5.4, $C_4=C_4(t)$ contains a factor $e^{-at}$, $a>0$, if $\inf \sigma (D^2) > 0$.

For (6.8) we have to estimate 
\be
\int\limits_M \int\limits_M |f(m)|^{-2} 
|(W(\frac{s}{2},m,p), \eta^{op}(p) \cdot)_p|^2 \,\, dp \,\, dm
\end{equation}          
We recall a simple fact about Hilbert spaces. Let $X$ be a Hilbert space, 
$x \in X, x \ne 0$. Then $|x|= \sup\limits_{|y|=1} | \langle x,y \rangle |$,
\be
|x|^2 = \Big( \sup_{|y|=1} | \langle x,y \rangle | \Big)^2  .
\end{equation}
This follows from $|\langle x,y \rangle| \le |x| \cdot |y|$
and equality for $y = \frac{x}{|x|}$. We apply this to $E \rightarrow M$,
$X = L_2 (M,E,dp)$, $x=x(m)=W(t,m,p),  \eta^{op}(p) \cdot)_p
= W(t,m,p) \circ \eta^{op}(p)$ and have to estimate
\be
\sup_{{\scriptsize
\begin{array}{c}
\Phi \in C^\infty_c(E) \\
|\Phi|_{L_2} = 1
\end{array}}
}      
N(\Phi) = 
\sup_{{\scriptsize
\begin{array}{c}
\Phi \in C^\infty_c(E) \\
|\Phi|_{L_2} = 1
\end{array}}
} 
|\langle \delta(m), e^{-tD^2} \eta ^{op} \Phi \rangle|_{L_2}
\end{equation}
According to 5.4, 5.5 and 5.7 
\be
W(t,m,\cdot) \in H^{\frac{r}{2}}(E), \quad 
|W(t,m,\cdot)|_{H^{\frac{r}{2}}} \le C_5(t) .                                     
\end{equation}
Hence we have can restrict in (6.15) to
\be
\sup_{{\scriptsize
\begin{array}{c}
\Phi \in C^\infty_c(E) \\
|\Phi|_{L_2} = 1 \\
|\Phi|_{H^{\frac{r}{2}}} \le C_5
\end{array}}
}      
N(\Phi) 
\end{equation}
In the sequel we estimate (6.17). For doing this, we recall some simple
facts concerning the wave equation
\be
\frac{\partial \Phi_s}{\partial s} = i D \Phi_s,
\quad \Phi_0 = \Phi,
\quad \Phi \,\,\, C^1 \en 
\mbox{with compact support.}
\end{equation}
It is well known that (6.18) has a unique solution $\Phi_s$ which ist given
by
\be
\Phi_s = e^{i s D} \Phi
\end{equation}
and
\be
\mbox{supp} \en \Phi_s \subset U_{|s|} \en( \mbox{supp} \en \Phi)
\end{equation}                           
$U_{|s|} = |s|$ -- neighborhood. Moreover,
\be
|\Phi_s|_{L_2} =  |\Phi|_{L_2}, \quad
|\Phi_s|_{H^{\frac{r}{2}}} = |\Phi|_{H^{\frac{r}{2}}} .       
\end{equation}
We fix a uniformly locally finite cover
${\cal U} = \{ U_\nu \}_\nu = \{B_d(x_\nu)\}_\nu$
by normal charts of radius $d < r_{inj}(M,g)$ and associated 
decomposition of unity $\{ \varphi_\nu \}_\nu$ satisfying
\be
|\nabla^i \varphi_\nu| \le C \en \mbox{for all} \en \nu, \en
0 \le i \le k+2
\end{equation}

Write
\bea
   N(\Phi) & = & | \langle \delta(m) , e^{-tD^2} \eta^{op} \Phi \rangle | 
   \nonumber \\
   & = & \frac{1}{\sqrt{4 \pi t}} \en \Big| \langle \delta(m),
   \int\limits^{+ \infty}_{- \infty} e^{\frac{-s^2}{4t}} e^{isD}
   (\eta^{op} \Phi) \,\, ds \rangle \Big|_{L_2(dp)} \en = \nonumber \\
   & = & \frac{1}{\sqrt{4 \pi t}} \en \Big| 
   \int\limits^{+ \infty}_{- \infty} e^{\frac{-s^2}{4t}} (e^{isD}
   \eta^{op} \Phi) (m) \,\, ds \Big|_{L_2(dp)} .
\eea
We decompose
\be
 \eta^{op}(\Phi) = \sum\limits_\nu \varphi_\nu \eta^{op} \Phi .
\end{equation}  
(6.24) is a locally finite sum, (6.18) linear. Hence
\be
 (\eta^{op}(\Phi))_s = \sum\limits_\nu (\varphi_\nu \eta^{op} \Phi)_s .
\end{equation}  
Denote as above
\[ | \en |_{p,i} \equiv | \en |_{W^{p,i}} , \]  
in particular
\be 
  | \en |_{2,i} \equiv | \en |_{W^{2,i}} \sim | \en |_{H^i} , \qu i \le k .  
\end{equation}
Then we obtain from (6.21), (6.22), (2.1)
\[ |(\varphi_\nu \eta^{op} \Phi)_s |_{H^{\frac{r}{2}}}    
   = |\varphi_\nu \eta^{op} \Phi |_{H^{\frac{r}{2}}}     
   \le C_6 |\varphi_\nu \eta^{op} \Phi|_{2, \frac{r}{2}} \le \]
\be                                               
  \le C_7 | \eta^{op} \Phi|_{2, \frac{r}{2}, U_\nu}     
  \le C_8 | \eta |_{2, \frac{r}{2}, U_\nu}                                                           
  \le C_9 | \eta |_{1, r-1, U_\nu}                
\end{equation}
since $r-1-\frac{n}{i} \ge \frac{r}{2} - \frac{n}{2}, 
r-1 \ge \frac{r}{2} , 2 \ge i$ for $r > n+2$ and 
$|\Phi|_{H^{\frac{r}{2}}} \le C_5$. This yields together with (2.3)
the estimate
\bea   
  |(\eta^{op} \Phi)_s (m)| & \le & C_{10} \cdot
   \sum\limits_{{\scriptsize
   \begin{array}{c}
      \nu \\
      m \in U_s(U_\nu)
    \end{array}}}      
   |(\varphi_\nu \eta^{op} \Phi)_s|_{2,\frac{r}{2}} \en \le  \nonumber \\
   & \le & C_{11} \cdot
   \sum\limits_{{\scriptsize
   \begin{array}{c}
      \nu \\          
      m \in U_s(U_\nu) 
    \end{array}}}
   | \eta |_{1, r-1, U_\nu}  
   \le C_{12} \cdot |\eta|_{1, r-1, B_{2d+|s|}(m)} \en = \nonumber \\
  & = & C_{12} \cdot vol( B_{2d+|s|}(m)) \cdot
  \left( \frac{1}{vol B_{2d+|s|}(m)} \cdot
  |\eta|_{1, r-1, B_{2d+|s|}(m)} \right) . \nonumber \\
       { }   
\eea 
There exist constants $A$ and $B$, independent of $m$ s. t.
\[ vol( B_{2d+|s|}(m)) \le A \cdot e^{B_{|s|}} . \]
Write 
\be
  e^{- \frac{s^2}{4t}} \cdot  vol( B_{2d+|s|}(m)) \le     
  C_{13} \cdot e^{-\frac{9}{10}\frac{s^2}{4t}} ,  \quad
  C_{13} = A \cdot e^{10B^2t} ,
\end{equation}                                                  
thus obtaining
\[  N(\Phi) \le C_{14} \int\limits^\infty_0 
   e^{-\frac{9}{10}\frac{s^2}{4t}}    
   \left( \frac{1}{vol B_{2d+|s|}(m)} \cdot   
   |\eta|_{1, r-1, B_{2d+|s|}(m)} \right) \,\, ds , \]
\[ C_{14} = C_{12} \cdot C_{13} = C_{12} \cdot A \cdot e^{10B^2t} . \] 
Now we apply Lemma 2.7 with $R=3d+s$ and infer
\bea
   & {} & \int\limits_M  \frac{1}{vol B_{2d+|s|}(m)} \cdot
   |\eta|_{1, r-1, B_{2d+|s|}(m)} \,\, dm \en \le \nonumber \\
   & \le & |\eta|_{1,r-1} + C (3d+s) \cdot (2d+s) 
   |\nabla \eta|_{1,r-1} \en \le \nonumber \\
   & \le & |\eta|_{1,r-1} + C (3d+s) \cdot (2d+s) |\eta|_{1,r-1} .
\eea
$C (3d+s)$ depends on $3d+s$ at most linearly exponentielly, i. e.
\[ C (3d+s) \cdot (2d+s) \le A_1 e^{B_1s} . \]
This implies
\bea
&& \int\limits^\infty_0 
   e^{-\frac{9}{10}\frac{s^2}{4t}}    
   \int\limits_M        
   \frac{1}{vol B_{2d+|s|}(m)} \cdot   
   |\eta|_{1, r-1, B_{2d+|s|}(m)} \,\, dm \,\, ds \\
& \le & = \int\limits^\infty_0 
   e^{-\frac{9}{10}\frac{s^2}{4t}} 
   (|\eta|_{1,r-1} + C (3d+s) \cdot (2d+s)
   |\eta|_{1,r-1}) \,\, ds \non \\
& \le & \int\limits^\infty_0
   e^{-\frac{8}{10}\frac{s^2}{4t}} \,\, ds ( |\eta|_{1,r-1} + A_1
   e^{10B^2_1t} |\eta|_{1,r} \non \\\
& = & \sqrt{t} \cdot \frac{1}{2} \sqrt{5 \pi}
   ( |\eta|_{1,r-1} + A_1 e^{10B^2_1t} |\eta|_{1,r}) \en < \en \infty . \non
\eea
The function
$\R_+ \times M \rightarrow \R$, 
\[ (s,m) \rightarrow e^{-\frac{9}{10}\frac{s^2}{4t}}
   \left( \frac{1}{vol B_{2d+|s|}(m)} \cdot   
   |\eta|_{1, r-1, B_{2d+|s|}(m)} \right) \]
is measurable, nonnegative, the integrals (6.30), (6.31) exist, hence
according to the principle of Tonelli, this function is $1$--summable,
the Fubini theorem is applicable and 
\[ \tilde \eta := C_{10} \cdot  \int\limits^\infty_0 
   e^{-\frac{9}{10}\frac{s^2}{4t}}    
   \left( \frac{1}{vol B_{2d+|s|}(m)} \cdot   
   |\eta|_{1, r-1, B_{2d+|s|}(m)} \right) \,\, ds \]
is (for $\eta \not\equiv 0)$ everywhere $\not= 0$ and $i$--summable.
We proved
\be
  \int |(W(t,m,p), \eta^{op} \cdot )_p|^2 \le \tilde \eta (m)^2 .
\end{equation}
Now we set
\be
  f(m) = (\tilde \eta (m))^{\frac{1}{2}}
\end{equation}
and infer $f(m) \not= 0$ everywhere, $f \in L_2$ and 
\bea
|f^{-1} e^{-\frac{s}{2} D^2} \circ \eta|^2_{L_2}
&=& \int\limits_M \int\limits_M f(m)^{-2} |((W(\frac{s}{2},m,p),\eta^{op})_p|^2
 \,\, dp \,\, dm \non \\
& \le & \int\limits_M \frac{1}{\tilde{\eta}(m)} \tilde{\eta} (m)^2 \,\, dm
= \int\limits_M \tilde{\eta} (m) \,\, dm \non \\
& \le & C_{12} \cdot A \cdot e^{10B^2s} \sqrt{s} \cdot \frac{1}{2} \sqrt{5 \pi}
(|\eta|_{1,r-1} + A_1 e^{10B^2_1s} |\eta|_{1,r}) \non \\
& \le & C_{15} \sqrt{s} e^{10B^2s} |\eta|_{1,r}, \no
\eea
i. e. 
\be
|f^{-1} e^{-\frac{s}{2}D^2} \circ \eta|_2 \le C^{\frac{1}{2}}_{15} \cdot s^{\frac{1}{4}}
\cdot e^{5B^2s} \cdot |\eta|^{\frac{1}{2}}_{1,r}.
\end{equation}
Here according to the term $A_1 e^{10B^2_1s}$, $C_{15}$ still depends on $s$.

We obtain
\bea
|e^{-\frac{s}{2}D^2} \circ f|_{L_2} \cdot |f^{-1} \circ e^{-\frac{s}{2}D^2} \circ \eta| 
& \le & C_4 |f|_{L_2} \cdot C^{\frac{1}{2}}_{15} \cdot s^{\frac{1}{4}} \cdot
e^{5B^2s} \cdot |\eta|^{\frac{1}{2}}_{1,r} \non \\
& \le & C_4 \cdot C_{15} \sqrt{s} e^{10B^2s} |\eta|_{1,r} \non \\
& = & C_{16} \cdot \sqrt{s} \cdot e^{10 B^2 s} |\eta|_{1,r} .
\eea
This yields $e^{-sD^2} \circ \eta$ is of trace class,
\bea
|e^{-sD^2} \eta|_1 & \le & |e^{-\frac{s}{2}D^2} \circ f|_2 \cdot
|f^{-1} e^{-\frac{s}{2}D^2} \eta|_2 \non \\
& \le & C_{16} \sqrt{s} e^{10B^2s} |\eta|_{1,r},
\eea
$e^{-sD^2} \circ \eta \circ D' \circ e^{-(t-s){D'}^2}$ is of trace class,
\bea
|e^{-sD^2} \circ \eta \circ D' \circ e^{-(t-s){D'}^2}|_1
& \le & |e^{-sD^2} \eta|_1 \cdot |D' e^{-(t-s){D'}^2}|_{op} \non \\
& \le & C_{16} \sqrt{s} e^{10B^2s} |\eta|_{1,r} \cdot C' \cdot \frac{1}{\sqrt{t-s}},
\eea
\bea
\left| \int\limits^t_{\frac{t}{2}} (e^{-sD^2} \circ \eta \circ D' \circ e^{-(t-s){D'}^2}
\,\, ds \right|_1
& \le &  \int\limits^t_{\frac{t}{2}} |e^{-sD^2} \eta \circ D' e^{-(t-s){D'}^2}|_1
\,\, ds \non \\
& \le & C_{16} \cdot C' \cdot e^{10B^2t} |\eta|_{1,r} \cdot \int\limits^t_{\frac{t}{2}} 
\left( \frac{s}{t-s} \right)^{\frac{1}{2}} \,\, ds , \non \\
&&
\eea
\[ \int\limits^t_{\frac{t}{2}} \left( \frac{s}{t-s} \right)^{\frac{1}{2}} \,\, ds 
   = [ \sqrt{s(t-s)} + \frac{t}{2} \arcsin \frac{2s-t}{t} ]^t_{\frac{t}{2}}
   = - \frac{t}{2} + \frac{t}{2} \frac{\pi}{2} = \frac{t}{2} ( \frac{\pi}{2} - 1) , \]
\bea
\left| \int\limits^t_{\frac{t}{2}} (e^{-sD^2} \circ \eta \circ D' \circ e^{-(t-s){D'}^2}
\,\, ds \right|_1
& \le & C_{16} \cdot C' \cdot e^{10B^2t} \cdot (\frac{\pi}{2}-1) \cdot \frac{t}{2}
|\eta|_{1,r} \non \\
& = & C_{17} e^{10B^2t} \cdot t \cdot |\eta|_{1,r}.
\eea
Here $C_{17} = C_{17}(t)$ and $C_{17}(t)$ can grow exponentially in $t$ if the 
volume grows exponentially. (6.40) expresses the fact that ${\rm (I_4)}$ is of
trace class and its trace norm is uniformly bounded on any $t$--intervall
$[a_0,a_1]$, $a_0>0$. The treatment of ${\rm (I_1) - (I_3)}$ is quite parallel to
that of ${\rm (I_4)}$. Write the integrand
of ${\rm (I_3)}$, ${\rm (I_2)}$ or 
${\rm (I_1)}$ as 
\be
(D e^{-\frac{s}{2}D^2}) [(e^{-\frac{s}{4}D^2}f) (f^{-1} e^{-\frac{s}{4}D^2} \eta)]
e^{-(t-s){D'}^2}
\end{equation}
or
\be
(e^{-sD^2}) [(\eta e^{-\frac{(t-s)}{4}{D'}^2}f^{-1}) (f e^{-\frac{(t-s)}{4}{D'}^2})]
D' e^{-\frac{t-s}{2}{D'}^2}
\end{equation}
or
\be
(e^{-sD^2} D) [( \eta e^{-\frac{(t-s)}{2}{D'}^2}f^{-1}) (f e^{-\frac{(t-s)}{2}{D'}^2} )],
\end{equation}
respectively. Then in the considered intervals the expression $[ \dots ]$ are of trace
class which can literally be proved as for ${\rm (I_4)}$. The main point in ${\rm (I_4)}$
was the estimate of $f^{-1} e^{-\pi D^2} \eta$. In (6.42), (6.43) we have to estimate
expressions $\eta e^{-\tau {D'}^2} f^{-1}$. Here we use the fact that $\eta=\eta^{op}$
is symmetric with respect to the fibre metric $h$: the endomorphism $\eta_{e_i}(\cdot)$
is skew symmetric as the Clifford multiplication $e_i \cdot$ which yields together that
$\eta^{op}$ is symmetric. Then the $L_2$--estimate of $(\eta^{op} \cdot W'(\tau,m,p),\cdot)$
is the same as that of $W'(\tau,m,p), \eta^{op} (p) \cdot)$ and we can perform the same
procedure as that starting with (6.12). The only distinction are other constants.
Here essentially enters the equivalence of the $D$-- and $D'$--Sobolev spaces i. e.
the symmetry of our uniform structure. The factors outside $[ \dots ]$ produce
$\frac{1}{\sqrt{s}}$ on $[ \frac{t}{2}, t ]$, $\frac{1}{\sqrt{t-s}}$ and $\frac{1}{\sqrt{s}}$
on $[0,\frac{t}{2}]$ (up to constants). Hence ${\rm (I_1) - (I_3)}$ are of trace class
with uniformly bounded trace norm on any $t$--intervall $[a_0,a_1]$, $a_0>0$.
This finishes the proof of theorem 6.1.
\qed
\end{proof}

For our later applications we need still the trace class property of
\be
e^{-t D^2} D - e^{-t {D'}^2} D' .
\end{equation}
Consider the decomposition
\bea
e^{-t D^2} D - e^{-t {D'}^2} D' &=& e^{-\frac{t}{2}D^2} D 
(e^{-\frac{t}{2} D^2} - e^{-\frac{t}{2} {D'}^2}) \\
&+& (e^{-\frac{t}{2} D^2} D - e^{-\frac{t}{2} {D'}^2} D') e^{-\frac{t}{2}{D'}^2} .
\eea
According to 6.1, $e^{-\frac{t}{2} D^2} - e^{-\frac{t}{2} {D'}^2}$ is for $t>0$
of trace class. Moreover, $e^{-\frac{t}{2} D^2} D = D e^{-\frac{t}{2} D^2}$
is for $t>0$ bounded, its operator norm is $\le \frac{C}{\sqrt{t}}$. Hence their 
product is for $t>0$ of trace class and has bounded trace norm for $t \in [a_0,a_1]$,
$a_0>0$. (6.45) is done. We can write (6.46) as
\bea
(e^{-\frac{t}{2} D^2} D - e^{-\frac{t}{2} {D'}^2} D') e^{-\frac{t}{2}{D'}^2} 
&=&
[ e^{-\frac{t}{2}D^2} (D-D') + ( e^{-\frac{t}{2}D^2} - e^{-\frac{t}{2}{D'}^2})D' ]
\cdot e^{-\frac{t}{2}{D'}^2} 
\non \\
&=& 
[-e^{\frac{t}{2}D^2} \eta] e^{-\frac{t}{2}{D'}^2}
+ [ \int\limits^{\frac{t}{2}}_0 e^{-sD^2} D \eta e^{-(\frac{t}{2}-s){D'}^s} \,\, ds \non \\
&&
+ \int\limits^{\frac{t}{2}}_0 e^{-sD^2} \eta D' e^{-(\frac{t}{2}-s){D'}^2} \,\, ds ]
(D' e^{-\frac{t}{2}{D'}^2}).
\eea
Now
\be
[e^{-\frac{t}{2}D^2} \eta] \cdot e^{-\frac{t}{2}{D'}^2} = [ (e^{-\frac{t}{4}D^2} f) 
(f^{-1} e^{-\frac{t}{4}D^2} \eta ) ] e^{-\frac{t}{2}{D'}^2}.
\end{equation}
(6.48) is of trace class and its trace norm is uniformly bounded on any
$[a_0,a_1]$, $a_0>0$, according the proof of 6.1. If we decompose
$\int\limits^{\frac{t}{2}}_0 = \int\limits^{\frac{t}{4}}_0 + 
\int\limits^{\frac{t}{2}}_{\frac{t}{4}}$
then we obtain back from the integrals in (6.47) the integrals $(I_1) - (I_4)$, 
replacing $t \ra \frac{t}{2}$. These are done. 
$D' e^{-\frac{t}{2}{D'}^2}$ generates $C/\sqrt{t}$ in the estimate of
the trace norm. Hence we proved

\begin{theorem}
Assume $(E,\nabla) \lra (M^n,g)$ with $(I),(B_k)$, $(E,\nabla)$ with $(B_k)$,
$k \ge r > n+2$, $n \ge 2$, 
$\nabla' \in \comp(\nabla) \cap {\cal C}_E(B_k) \subset {\cal C}^{1,r}_E(B_k)$,
$D=D(g,\nabla)$, $D'=D(g,\nabla')$ generalized Dirac operators. Then
\[ e^{-tD^2} - e^{-t{D'}^2} \]
and
\[ D e^{-tD^2} - D' e^{-t{D'}^2} \] 
are trace class operators for $t > 0$ and their trace norm is uniformly bounded
on compact $t$--intervalls $[a_0,a_1]$, $a_0>0$.
\qed
\end{theorem}

\section{Trace class property for variation of the Clifford structure}

\setcounter{equation}{0}

Our intention is to admit much more general perturbations than those of
$\nabla = \nabla^h$ only. Nevertheless, the discussion of more general perturbations
is modelled by the case of $\nabla$--perturbation. In this section, we admit 
perturbations of $g, \nabla^h, \cdot$, fixing $h$, the topology and vector bundle
structure of $E \lra M$. The main result of this section shall be formulated
as follows.

\begin{theorem}
Let $E=(E,h,\nabla=\nabla^h,\cdot) \lra (M^n,g)$ be a Clifford bundle with
$(I)$, $(B_k(M,g))$, $(B_k(E,\nabla))$, $k \ge r+1 > n+3$, 
$E'=(E,h,\nabla'={\nabla'}^h,\cdot') \lra (M^n,g') \in \gencomp^{1,r+1}_{L,diff,F}
(E) \cap \cl {\cal B}^{N,n} (I,B_k)$, 
$D=D(g,h,\nabla=\nabla^h,\cdot)$, $D'=D(g',h,\nabla'={\nabla'}^h,\cdot')$
the associated genera\-li\-zed Dirac operators. Then for $t>0$
\be
e^{-tD^2} - e^{-t{D'}^2_{L_2}}
\end{equation}
is of trace class and the trace norm is uniformly bounded on compact $t$--intervalls
$[a_0,a_1]$, $a_0>0$.
\end{theorem}

Here ${D'}^2_{L_2}$ is the unitary transformation of ${D'}^2$ to 
$L_2 = L_2((M,E),g,h)$. 7.1. needs some explanations. $D$ acts in 
$L_2 = L_2((M,E),g,h)$, $D'$ in $L'_2 = L_2((M,E),g',h)$. $L_2$ and $L'_2$
are quasi isometric Hilbert spaces. As vector spaces they coincide, their
scalar products can be quite different but must be mutually bounded at the 
diagonal after multiplication by constants. $D$ is self adjoint on ${\cal D}_D$
in $L_2$, $D'$ is self adjoint on ${\cal D}_{D'}$ in $L'_2$ but not necessarily
in $L_2$. Hence $e^{-t{D'}^2}$ and $e^{-tD^2} - e^{-t{D'}^2}$ are not defined
in $L_2$. One has to graft $D^2$ or ${D'}^2$. Write 
$dvol_q(g) \equiv d q(g) = \alpha (q) \cdot d q(g') \equiv dvol_q(g')$.
Then
\bea
&& 0 < c_1 \le \alpha (q) \le c_2, \alpha, \alpha^{-1} \mbox{ are }
(g,\nabla^g)- \mbox{ and } (g',\nabla^{g'})- \mbox{bounded} \non \\
&& \mbox{up to order } 3, |\alpha-1|_{g,1,r+1}, |\alpha-1|_{g',1,r+1} < \infty,
\eea
since $g' \in \comp^{1,r+1}(g)$. Define 
$U : L_2 \lra L'_2$, $U \Phi = \alpha^{\frac{1}{2}} \Phi$. 
Then $U$ is a unitary equivalence between $L_2$ and $L'_2$, 
$U^* = U^{-1}$. $D'_{L_2} := U^* D' U$ acts in $L_2$, is self adjoint on 
$U^{-1} ({\cal D}_{D'})$, since $U$ is a unitary equivalence. The same
holds for ${D'}^2_{L_2} = U^* {D'}^2 U = (U^* D' U)^2$. 
It follows from the definition of the spectral measure, the spectral
integral and the spectral representations
${D'}^2 = \int \lambda^2 \,\, d E'_\lambda$, 
$e^{-t{D'}^2} = \int e^{-t\lambda^2} \,\, d E'_\lambda$
that
${D'}^2_{L_2} = U^* {D'}^2 U = U^* \int \lambda^2 \,\, d E'_\lambda U =
\int \lambda^2 \,\, d(U^* E'_\lambda U)$ 
and
\be 
e^{-t{D'}^2_{L_2}} = \int e^{-t\lambda^2} \,\, d (U^* E'_\lambda U) 
= U^* (\int e^{-t\lambda^2} \,\,
dE'_\lambda) U = U^* e^{{-tD'}^2} U.
\end{equation}
In (7.1) $e^{-t{D'}^2_{L_2}}$ means 
$e^{-t{D'}^2_{L_2}} = e^{-t(U^* D' U)^2} = U^* e^{-t{D'}^2} U$.
We obtain from $g' \in \comp^{1,r+1}(g)$, 
${\nabla'}^h \in \comp^{1,r+1}(\nabla^hg)$, 
$\cdot' \in \comp^{1,r+1}(\cdot)$, 
$D-\alpha^{-\frac{1}{2}} D' \alpha^{\frac{1}{2}} = 
D - D' - \frac{grad' \alpha \cdot'}{2 \alpha}$
and (7.2) the following lemma

\begin{lemma}
$W^{1,i}(E,g,h,\nabla^h) = W^{1,i}(E,g',h,{\nabla'}^h)$
as equivalent Banach spaces, $0 \le i \le r+1$. 
\qed
\end{lemma}

\begin{coro}
$W^{2,i}(E,g,h,\nabla^h) = W^{2,i}(E,g',h,{\nabla'}^h)$
as equivalent Hilbert spaces, $0 \le j \le \frac{r+1}{2}$. 
\qed
\end{coro}

\begin{coro}
$H^j(E,D) \cong H^jK(E,D')$, $0 \le j \le \frac{r+1}{2}$. 
\qed
\end{coro}

7.2 has a parallel version for the endomorphism bundle $\mbox{End} E$.

\begin{lemma}
$\Omega^{1,1,i} (\mbox{End} E, g, h, \nabla^h) \cong \Omega^{1,1,i} (\mbox{End} E, g', h, {\nabla'}^h)$ 
$0 \le i \le r+1$.

\qed
\end{lemma}

\begin{lemma}
$\Omega^{1,2,j} (\mbox{End} E, g, h, \nabla^h) \cong \Omega^{1,2,j} (\mbox{End} E, g', h, {\nabla'}^h)$ 
$0 \le j \le \frac{r+1}{2}$.

\qed
\end{lemma}

$e^{-t{D'}^2_{L_2}} : L_2 \lra L_2$ has evidently the heat kernel
\[ W'_{L_2} (t,m,p) = \alpha^{-\frac{1}{2}} (m) W'(t,m,p) \alpha^{\frac{1}{2}} (p) \]
$W' \equiv W_{L'_2}$. Our next task is to obtain an explicit expression for
$e^{-tD^2} - e^{-t{D'}^2_{L_2}}$. For this we apply again Duhamel's principle. The
steps 1. -- 4. in the proof of 6.3. remain. Then we set 
$\Phi(t,q)=W(t,m,q)$, $\Psi(t,q)=W'_{L_2}(t,m,q)$ 
and obtain
\bea
&&  
- \int\limits^\beta_\alpha \int\limits_M h_q (W(\tau,m,q),(D^2+\frac{\partial}{\partial t})
W'_{L_2} (t-\tau,q,p)) \,\, dq(g) \,\, d\tau = 
\non \\
&& 
= \int\limits_M [h_q(W(\beta,m,q),W'_{L_2}(t-\beta,q,p) - h_q(W(\alpha,m,q),
W'_{L_2}(t-\alpha,q,p)] \,\, dq(g).
\non
\eea
Performing $\alpha \lra 0^+, \beta \lra t$ and using $dq(g) = \alpha(q) dq(g')$
yields
\bea
&&  
- \int\limits^t_0 \int\limits_M h_q (W(s,m,q),(D^2+\frac{\partial}{\partial t})
W' (t-s,q,p)) \,\, dq(g) \,\, ds = 
\non \\
&& 
= - \int\limits^t_0 \int\limits_M [h_q(W(s,m,q), (D^2-{D'}^2_{L_2}) W'_{L_2}(t-s,q,p)
\,\, dq(g) \,\, ds
\non \\
&&
= W(t,m,p) \alpha(p) - W'_{L_2}(t,m,p). 
\eea
(7.4) expresses the operator equation
\bea
e^{-tD^2} \alpha - e^{-t{D'}^2_{L_2}} &=& - \int\limits^t_0 e^{-sD^2} 
(D^2-{D'}^2_{L_2}) e^{-(t-s){D'}^2_{L_2}} \,\, ds. \non \\
e^{-tD^2} \alpha - e^{-t{D'}^2_{L_2}} &=& e^{-tD^2} (\alpha-1) + 
e^{-tD^2} - e^{-t{D'}^2_{L_2}}, \mbox{ hence} \non \\
e^{-tD^2} - e^{-t{D'}^2_{L_2}} &=& - e^{-tD^2} (\alpha-1) - 
\int\limits^t_0 e^{-sD^2} (D^2-{D'}^2_{L_2}) e^{-(t-s){D'}^2_{L_2}} \,\, ds. \non \\
&& 
\eea
As we mentioned in (7.2),
$(\alpha-1) = \frac{dq(g)}{dq(g')} - 1 = \frac{\sqrt{\det g}}{\sqrt{\det g'}} - 1 
\in \Omega^{0,1,r+1}$ since $g \in \comp^{1,r+1} (g)$. We write
$e^{-tD^2} (\alpha-1) = (e^{-\frac{t}{2}D^2}f) (f^{-1} e^{-\frac{t}{2}D^2} (\alpha-1))$,
determine $f$ as in section 6 from $\eta_\alpha = \alpha -1$ and obtain
$e^{-tD^2} (\alpha-1)$ is of trace class with trace norm uniformly bounded on any
$t$--interval $[a_0,a_1]$, $a_0>0$. Decompose 
$D^2 - {D'}^2_{L_2} = D (D - D'_{L_2}) + (D - D'_{L_2}) D'_{L_2}$.
We need explicit analytic expressions for this.
$D (D - D'_{L_2}) = D (D - \alpha^{-\frac{1}{2}} D' \alpha^{\frac{1}{2}})
= D (D - D') - D \frac{\grad' \alpha \cdot'}{2 \alpha}$, 
$(D - D'_{L_2}) D'_{L_2} = ((D - D') - \frac{\grad' \alpha \cdot'}{2 \alpha} ) 
\alpha^{-\frac{1}{2}} D' \alpha^{\frac{1}{2}} $.
If we set again $D-D'=-\eta$ then we have to consider as in section 6 with
$\frac{\grad' \alpha}{2 \alpha} = \frac{\grad' \alpha \cdot'}{2 \alpha}$
where $\grad' \equiv \grad_{g'}$
\bea
&& \int\limits^{\frac{t}{2}}_0 e^{-sD^2} D (\eta - \frac{\grad' \alpha}{2 \alpha}) 
e^{-(t-s){D'}^2_{L_2}} \,\, ds + \non \\
&& \int\limits^{\frac{t}{2}}_0 e^{-sD^2} (\eta - \frac{\grad' \alpha}{2 \alpha}) 
D'_{L_2} e^{-(t-s){D'}^2_{L_2}} \,\, ds + \non \\
&& \int\limits^t_{\frac{t}{2}} e^{-sD^2} D (\eta - \frac{\grad' \alpha}{2 \alpha}) 
e^{-(t-s){D'}^2_{L_2}} \,\, ds + \non \\
&& \int\limits^t_{\frac{t}{2}} e^{-sD^2} (\eta - \frac{\grad' \alpha}{2 \alpha}) 
D'_{L_2} e^{-(t-s){D'}^2_{L_2}} \,\, ds . \non 
\eea
It follows immediately from $g' \in \comp^{1,r+1}(g)$ that the vector field
$\frac{\grad' \alpha}{\alpha} \in \Omega^{0,1,r} (TM)$.
If we write
$\eta^{op}_0 = - \frac{\grad' \alpha \cdot'}{\alpha}$ 
then $\eta^{op}_0$ is a zero order operator, 
$|\eta_0|_r < \infty$
and we literally repeat the procedure for $(I_1) - (I_4)$ in section 6, 
inserting $\eta_0 = - \frac{\grad' \alpha \cdot'}{\alpha}$
for $\eta$ there. Hence there remains to discuss the integrals
\bea
&& \int\limits^t_0 e^{-sD^2} D \eta e^{-(t-s){D'}^2_{L_2}} \,\, ds + \non \\
&& \int\limits^t_0 e^{-sD^2} \eta D'_{L_2} e^{-(t-s){D'}^2_{L_2}} \,\, ds.
\eea
The first main step is to insert explicit expressions for $D-D'$. 
Let $m_0 \in M$, $U=U(m_0)$ a manifold and bundle coordinate neighborhood with
coordinates $x^1, \dots, x^n$ and local bundle basis
$\Phi_1, \dots, \Phi_n : U \lra E|_U$. Setting
$\nabla_{\frac{\partial}{\partial x_i}} \Phi_\alpha \equiv \nabla_i \Phi_\alpha =
\Gamma^\beta_{i\alpha} \Phi_\beta$, 
$\nabla \Phi_\alpha = d x^i \otimes \Gamma^\beta_{i\alpha} \Phi_\beta$,
we can write 
$D \Phi_\alpha = \Gamma^\beta_{i\alpha} g^{ik} 
\frac{\partial}{\partial x^k} \cdot \Phi_\beta$,
$D' \Phi_\alpha = {\Gamma'}^\beta_{i\alpha} {g'}^{ik} 
\frac{\partial}{\partial x^k} \cdot' \Phi_\beta$,
or for a local section $\Phi$
\be
D \Phi = g^{ik} \frac{\partial}{\partial x^k} \cdot \nabla_i \Phi, \quad
D' \Phi = {g'}^{ik} \frac{\partial}{\partial x^k} \cdot' \nabla'_i \Phi.
\end{equation}
This yields
\bea
- (D-D') \Phi &=& g^{ik} \frac{\partial}{\partial x^k} \cdot \nabla_i \Phi
- {g'}^{ik} \frac{\partial}{\partial x^k} \cdot' \nabla'_i \Phi \non \\
&=& [(g^{ik}-{g'}^{ik}) \frac{\partial}{\partial x^k} \cdot \nabla_i 
+ {g'}^{ik} \frac{\partial}{\partial x^k} \cdot (\nabla_i - \nabla'_i) \non \\
&& + {g'}^{ik} \frac{\partial}{\partial x^k} (\cdot - \cdot') \nabla'_i] \Phi, 
\eea
i. e. we can write
\be
- (D-D') \Phi = (\eta^{op}_1 + \eta^{op}_2 + \eta^{op}_3) \Phi,
\end{equation}
where locally
\bea
\eta^{op}_1 \Phi &=& (g^{ik} - {g'}^{ik}) \frac{\partial}{\partial x^k} \cdot \nabla_i \Phi , \\
\eta^{op}_2 \Phi &=& {g'}^{ik} \frac{\partial}{\partial x^k} \cdot (\nabla_i - \nabla'_i) \Phi , \\
\eta^{op}_3 \Phi &=& {g'}^{ik} \frac{\partial}{\partial x^k} (\cdot - \cdot') \nabla'_i \Phi .
\eea
Here $({g'}^{ik}) = (g'_{jl})^{-1}$.         
We simply write $\eta_\nu$ instead $\eta^{op}_\nu$, hence
\bea
   \mbox{(7.6)} &=& 
   \int\limits^t_0 e^{-sD^2} D (\eta_1+\eta_2+\eta_3)                       
   e^{-(t-s){D'}^2_{L_2}} \,\, ds \en + \\
   &+& \int\limits^t_0 e^{-sD^2} (\eta_1+\eta_2+\eta_3)                       
   D'_{L_2} e^{-(t-s){D'}^2_{L_2}} \,\, ds.  
\eea                           
We have to estimate
\be
  \int\limits^t_0 e^{-sD^2} D \eta_\nu e^{-(t-s){D'}^2_{L_2}} \,\, ds 
\end{equation}                                                     
and
\be
  \int\limits^t_0 e^{-sD^2} \eta_\nu D'_{L_2} e^{-(t-s){D'}^2_{L_2}} \,\, ds .
\end{equation}
Decompose $\int\limits^t_0 = \int\limits^{\frac{t}{2}}_0 + 
\int\limits^t_{\frac{t}{2}}$ which yields
\def\theequation{$I_{\nu,1}$}   
\be
   \int\limits^{\frac{t}{2}}_0 e^{-sD^2} D \eta_\nu e^{-(t-s){D'}^2_{L_2}} \,\, ds ,
\end{equation}
\def\theequation{$I_{\nu,2}$}  
\be
   \int\limits^{\frac{t}{2}}_0 e^{-sD^2} \eta_\nu D'_{L_2} e^{-(t-s){D'}^2_{L_2}} \,\, ds ,
\end{equation}
\def\theequation{$I_{\nu,3}$}  
\be
   \int\limits^t_{\frac{t}{2}} e^{-sD^2} D \eta_\nu  e^{-(t-s){D'}^2_{L_2}} \,\, ds ,
\end{equation}
\def\theequation{$I_{\nu,4}$}  
\be
   \int\limits^t_{\frac{t}{2}} e^{-sD^2} \eta_\nu D'_{L_2} e^{-(t-s){D'}^2_{L_2}} \,\, ds .
\end{equation}                         
\def\theequation{\thesection.\arabic{equation}}  
\setcounter{equation}{16}

$(I_{\nu,1}) - (I_{\nu,4})$ look as $(I_1) - (I_4)$ in section 6.
But in distinction to section 6, not all $\eta_\nu = \eta^{op}_\nu$ are operators
of order zero. Only $\eta_2$ is a zero order operator, 
generated by an $\mbox{End} E$
valued 1--form $\eta_2$. $\eta_1$ and $\eta_3$ are first order operators. We
start with $\nu=2$, $\eta_2 \cdot |\eta_2|_{1,r} < \infty$ is a consequence of
$E' \in \comp^{1,r+1}_{L,diff} (E)$ and we are from an analytical point of view
exactly in the situation of section 6. $(I_{2,1}) - (I_{2,4})$ can be estimated
quite parallel to $(I_1) - (I_4)$ in section 6 and we are done. There remains to
estimate $(I_{\nu,j})$, $\nu \neq 2$, $j=1,\dots,4$. We start with $\nu=1$, $j=3$
and write
\be
e^{-sD^2} D \eta_1 e^{-(t-s){D'}^2} = (De^{-\frac{s}{2}D^2}) (e^{-\frac{s}{4}D^2} 
\cdot f) (f^{-1} e^{-\frac{s}{4}D^2} \eta_1) (e^{-(t-s) {D'}^2}).
\end{equation}
$De^{-\frac{s}{2}D^2}$ and $e^{-(t-s) {D'}^2}$ are bounded in $[\frac{t}{2}, t]$
and we perform their estimate as in section 6. $e^{-\frac{s}{4}D^2} \cdot f$ is
Hilbert--Schmidt if $f \in L_2$. There remains to show that for appropriate $f$
\[ f^{-1} e^{-\frac{s}{4}D^2} \eta_1 \]
is Hilbert--Schmidt. Recall $r+1>n+3$, $n \ge 2$, which implies 
$\frac{r}{2} > \frac{n}{2}+1$, $r-1-n \ge \frac{r}{2} - \frac{n}{2}$, 
$r-1 \ge \frac{r}{2}$, $2 \ge i$, i. e. 2.4, 2.5 are available. If we write
in the sequel pointwise or Sobolev norms we should always write 
$|\Psi|_{g',h,m'}$, $|\Psi|_{H^\nu(E,D')}$, $|\Psi|_{g',h,\nabla',2,\frac{r}{2}}$, 
$|g-g'|_{g',m}$, $|g-g'|_{g',1,r}$ etc. or the same with respect to 
$g,h,\nabla,D$, depending on the situation. But we often omit the reference to
$g',h,\nabla',D,m,g,h \dots$ in the notation. The justification for doing 
this in the Sobolev case is the symmetry of our uniform structure.

Now
\bea
(\eta_1 \Phi) (m) &=& ((g^{ik} - {g'}^{ik}) \frac{\partial}{\partial x^k} \cdot 
\nabla_i \Phi)|_m , \\
|\eta_1 \Phi|_m &=& |\eta_1 \Phi|_{g,h,m} \non \\
& \le & C_1 \cdot |g-g'|_{g,m} \cdot 
\left( \sum\limits^n_{k=1} \left| \frac{\partial}{\partial x^k} \right|^2_{g,m} 
\right)^{\frac{1}{2}} \cdot
\left( \sum\limits^n_{i=1} \left| \nabla_i \Phi \right|^2_{h,m} 
\right)^{\frac{1}{2}} .
\non
\eea
To estimate $\sum\limits^n_{k=1} \left| \frac{\partial}{\partial x^k} \right|^2_{g,m}$ 
more concretely we assume that $x^1, \dots, x^n$ are normal coordinates with respect 
to $g$, i. e. we assume a (uniformly locally finite) cover of $M$ by normal charts
of fixed radius $\le r_{inj}(M,g)$. Then
$\left| \frac{\partial}{\partial x^k} \right|^2_{g,m} = g \left( 
\frac{\partial}{\partial x^k}, \frac{\partial}{\partial x^k} \right)
= g_{kk}(m)$, 
and there is a constant
$C_2 = C_2 (R, r_{inj}(M,g))$
s. t.
$\left( \sum\limits^n_{i=1} \left| \nabla_i \Phi \right|^2_{h,m} 
\right)^{\frac{1}{2}} \le C_2$.
Using finally
$|\nabla_X \Phi| \le |X| \cdot |\nabla \Phi|$,
we obtain
\be
|\eta_1 \Phi|_m \le C \cdot |g-g'|_g \cdot |\nabla \Phi|_{h,m}.
\end{equation}
(7.19) extends by the Leibniz rule to higher derivatives 
$|\nabla^k \eta_1 \Phi|_m$,
where the polynomials on the right hand side are integrable by the module structure
theorem (this is just the content of this theorem). (7.18), (7.19) also hold (with
other constants) if we perform some of the replacements $g \lra g'$, 
$\nabla \lra \nabla'$: We remark that the expressions
$D(g,h,\nabla^h,\cdot$, $D(g',h,\nabla^h,\cdot)$ 
are invariantly defined, hence
\be
[D(g,h,\nabla^h,\cdot) - D(g,h,\nabla^h,\cdot)] (\Phi|_U) = ((g^{ik} - {g'}^{ik})
\partial_k) \cdot \nabla_i (\Phi|_U) .
\end{equation}
We have to estimate the kernel of
\be
h_p (W(t,m,p), \eta^{op}_1 \cdot)
\end{equation}
in $L_2((M,E),g,h)$ and to show that this represents the product of two
Hilbert--Schmidt operators in $L_2 = L_2((M,E),g,h)$. We cannot immediately apply
the procedure starting with (6.27), (6.28) since $\eta^{op}_1$ is not of zero
order but we would be done if we could write (7.21) as
\be
(\eta^{op}_{1,1} (p) W(t,m,p), \eta^{op}_{1,0} \cdot),
\end{equation}
$\eta^{op}_{1,1}$ of first order, $\eta^{op}_{1,0}$ of zeroth order. Then we would
replace in (6.23), $\dots$, 
$W$ by $\eta^{op}_{1,1} (p) W(t,m,p)$, 
apply $k \ge r+1 > n+3$, 5.4, 5.5, 5.7 and obtain
\be
\eta^{op}_{1,1}W(t,m,\cdot) \in H^{\frac{r}{2}} (E), \quad
|W(t,m,\cdot)|_{H^{\frac{r}{2}}} \le C(t)
\end{equation}
and would then literally proceed as in (6.27) -- (6.40).

Let $\Phi \in C^\infty_c (U)$. Then
\bea
&& \int (W(t,m,p), \eta^{op}_1 (p) \Phi (p) )_p \,\, dvol_p(g) = \non \\
&& \int (((g^{ik} - {g'}^{ik}) \partial_k ) \cdot \nabla_i W, \Phi)_p \,\, dvol_p(g) - \non \\
&& - \int ( W, (\nabla_i (g^{ik} - {g'}^{ik}) \partial_k ) \cdot \Phi ) \,\, dvol_p(g) = \non \\
&& - \int ( \nabla_i W, (g^{ik} - {g'}^{ik}) \partial_k \cdot \Phi)_p \,\, dvol_p(g) - \non \\
&& - \int ( W, ( \nabla_i ((g^{ik} - {g'}^{ik}) \partial_k )) \cdot \Phi)_p \,\, dvol_p(g). \non
\eea
This can easily be globalized by introducing a u. l. f. cover by normal charts
$\{ U_\alpha \}_\alpha$ of fixed radius, an associated decomposition of unity
$\{ \phi_\alpha \}_\alpha$ as follows:
\bea
&& \int ( W, \eta^{op}_1 ( \sum \phi_\alpha \Phi )) = \sum\limits_\alpha \int
(W, \eta^{op}_1 (\phi_\alpha \Phi)) \non \\
&& = \sum\limits_\alpha \int ( \nabla_{\alpha,i} W, ((g^{ik}_\alpha - {g'}^{ik}_\alpha)
\partial_{\alpha,k} \cdot \phi_\alpha \Phi) - \non \\ 
&& - \sum\limits_\alpha \int ( W, (\nabla_{\alpha,i}
((g^{ik}_\alpha - {g'}^{ik}_\alpha) \partial_k )) \cdot \phi_\alpha \Phi) \non 
\eea
\bea
&& = - \int ( \sum\limits_\alpha \nabla_{\alpha,i} W, \phi_\alpha ((g^{ik}_\alpha - 
{g'}^{ik}_\alpha) \partial_{\alpha,k} ) \cdot \Phi) - \\
&& - \int ( W, \sum\limits_\alpha \phi_\alpha ( \nabla_{\alpha,i} ((g^{ik}_\alpha -
{g'}^{ik}_\alpha ) \partial_k )) \cdot \Phi )
\eea
Using (7.24), (7.25), we write
\bea
N(\Phi) &=& | \sk{\delta(m)}{e^{tD^2} \eta^{op}_1 \Phi} |_{L_2(M,E,dp)} \non \\
&=& | (W(t,m,p),\eta^{op}_1 \Phi)_p |_{L_2(M,E,dp)} \non \\
&=& | (\eta^{op}_{1,1}(p) W(t,m,p), \eta^{op}_{1,0} \Phi)_p + \non \\
&+& (W( t,m,p,\eta^{op}_{1,0,0} \Phi ))_p |_{L_2(M,E,dp)} .
\eea

Now we use 
$|\nabla_X \raisebox{0.3ex}{$\chi$}| \le |X| \cdot |\nabla \raisebox{0.3ex}{$\chi$}|$,
that the cover is u.l.f. and 
$|\nabla W| \le C_1 \cdot (|DW|+W)$
(since we have bounded geometry) and obtain
\bea
N(\Phi) & \le & C \cdot ( | ( D W(t,m,p), \eta^{op}_{1,0} \Phi |_{L_2(M,E,dp)}
+ | W(t,m,p,\eta^{op}_{1,0,0} \Phi |_{L_2(dp)} \non \\
& \equiv & C \cdot (N_1(\Phi) + N_2(\Phi)) .
\eea
Hence we have to estimate
\be
\sup\limits_{{\scriptsize \begin{array}{c} \Phi \in C^\infty_c(E) \\ |\Phi|_{L_2}=1
\end{array}}} N_1 (\Phi) =
\sup\limits_{{\scriptsize \begin{array}{c} \Phi \in C^\infty_c(E) \\ |\Phi|_{L_2}=1
\end{array}}}
| \sk{\delta(m)}{(De^{-tD^2}) \eta^{op}_{1,0} \Phi} |_{L_2{dp}}
\end{equation}
and
\be
\sup\limits_{{\scriptsize \begin{array}{c} \Phi \in C^\infty_c(E) \\ |\Phi|_{L_2}=1
\end{array}}} N_2 (\Phi) =
\sup\limits_{{\scriptsize \begin{array}{c} \Phi \in C^\infty_c(E) \\ |\Phi|_{L_2}=1
\end{array}}}
| \sk{\delta(m)}{(e^{-tD^2}) \eta^{op}_{1,0,0} \Phi} |_{L_2{dp}}.
\end{equation}
According to $k>r+1>n+3$, 5.4, 5.5, 5.7,
\bea
&& D(W(t,m,\cdot), W(t,m,\cdot) \in H^{\frac{r}{2}}(E), \non \\
&& |(D(W(t,m,\cdot)|_{H^{\frac{r}{2}}}, |W(t,m,\cdot)|_{H^{\frac{r}{2}}} \le C_1(t)
\eea 
and we can restrict in (7.28), (7.29) to 
\be
\sup\limits_{{\scriptsize \begin{array}{c} \Phi \in C^\infty_c(E) \\ |\Phi|_{L_2}=1 \\ 
| \Phi |_{H^{\frac{r}{2}}} \le C_1(t) \end{array}}} N_i (\Phi).
\end{equation}
$\eta^{op}_{1,0}$, $\eta^{op}_{1,0,0}$ are of order zero and we estimate them
by
\bea
C \cdot |g-g'|_{g,2,\frac{r}{2}} & \le & C' |g-g'|_{g,1,r-1} \\
\mbox{and} \qquad 
D \cdot |\nabla (g-g')|_{g,2,\frac{r}{2}} & \le & D' |\nabla (g-g')|_{g,1,r-1} \non \\
& \le & D'' |g-g'|_{g,1,r}
\eea
respectively. As we have seen already, into the estimate (7.33) enters 
$|\nabla \eta|_{1,r-1}$, i. e. in our case
$|\nabla^2 (g-g')|_{r-1} \sim |g-g'|_{r+1}$. For this reason we assumed
$E' \in \comp^{1,r+1}_{L,diff,F} (E)$, 
not as in section 6,
$E \in \comp^{1,r} (E')$.
In the expression for $N_1(\Phi)$ corresponding to (6.23) there is now a
slight deviation,
\be
N_1(\Phi) = \frac{1}{\sqrt{4 \pi t}} \frac{1}{2t} \left| \int\limits^{+ \infty}_{- \infty} 
s \cdot e^{-\frac{s^2}{4t}} e^{isD} \eta^{op}_{1,0} \Phi(m) \,\, ds \right| .
\end{equation}
We estimate in (7.34) 
$s \cdot e^{-\frac{1}{18} \frac{s^2}{4t}}$ by a constant, write instead of (6.29)
\[ e^{-\frac{17}{18} \frac{s^2}{4t}} \cdot vol (B_{2d+s}(m)) \le C \cdot
e^{-\frac{9}{10} \frac{s^2}{4t}} \]
and proceed now for $N_1(\Phi)$, $N_2(\Phi)$ literally as in (6.29) -- (6.40).
Hence (7.17) is of trace class, its trace norm in uniformly bounded on any
$t$--intervall $[a_0,a_1]$, $a_0>0$. $(I_{1,3})$ is done. $(I_{1,4})$ is
absolutely parallel to $(I_{1,3})$, even better, since the left hand factor
$D$ is missing. 
$|D'_{L_2} e^{-(t-s) D^{'2}_{L_2}}|_{op}$ 
now produces the factor
$\frac{1}{\sqrt{t-s}}$ which is integrable over $[\frac{t}{2}, t]$. Write the
integrand of $(I_{1,1})$ as 
\be
(D e^{-sD^2}) (\eta_1 e^{-\frac{(t-s)}{2} D^{'2}_{L_2}} f^{-1}) 
(f e^{-\frac{(t-s)}{2} D^{'2}_{L_2}}).
\end{equation}
We proceed with (7.35) as before. Here $\eta_1$ already stands at the right place,
we must not perform partial integration. Into the estimate enters again the first
derivative of $W'$. $D e^{-sD^2}$ generates the factor $\frac{1}{\sqrt{s}}$ which
is intgrable on $[0, \frac{t}{2}]$.
We write $(I_{1,2})$ as 
\be
\int\limits^{\frac{t}{2}}_0 e^{-sD^2} 
[(\eta_1 e^{-\frac{(t-s)}{4} D^{'2}_{L_2}} f^{-1})
(f e^{-\frac{(t-s)}{4} D^{'2}_{L_2}})] 
e^{-\frac{(t-s)}{2} D^{'2}_{L_2}} D^{'2}_{L_2} \,\, ds
\end{equation}
and proceed as before.

Consider finally the case $\nu=3$, locally
\[ \eta^{op}_3 \Phi = {g'}^{ik} \frac{\partial}{\partial x^k} (\cdot - \cdot')
\nabla'_i \Phi . \]
The first step in this procedure is quite similar as in the case $\nu=1$ to shift
the derivation to the left of $W$ and to shift all zero order terms to the right.

Let $X$ be a tangent vector field and $\Phi$ a section.

\begin{lemma}
$X (\cdot - \cdot') \nabla'_i \Phi = \nabla'_i (X (\cdot - \cdot') \Phi)$ +
zero order terms.
\end{lemma}

\begin{proof}
$X (\cdot - \cdot') \nabla'_i \Phi = [ X (\cdot - \cdot') \nabla'_i \Phi -
\nabla'_i (X (\cdot - \cdot') \Phi) ] + \nabla'_i (X (\cdot - \cdot') \Phi)$.
We are done if $[ \dots ]$ on the right hand side contains no derivatives of $\Phi$.
But an easy calculation yields
\bea
&& [X (\cdot - \cdot') \nabla'_i \Phi - \nabla'_i (X (\cdot - \cdot') \Phi ) ] = \non \\
&& = X \cdot (\nabla'_i - \nabla_i) \Phi - (\nabla'_i - \nabla_i) (X \cdot \Phi) \non \\
&& + (\nabla'_i - \nabla_i) X \cdot' \Phi + (\nabla_i X) (\cdot' - \cdot) \Phi.
\eea
\qed
\end{proof}

Hence for $\Phi, \Psi \in C^{\infty}_c(U)$
\bea
&& \int h (\Psi, {g'}^{ik} \frac{\partial}{\partial x^k} (\cdot - \cdot') 
\nabla'_i \Phi)_p dp(g) \,\, = \non \\
&& = \int h (\Psi, \nabla'_i ( {g'}^{ik} \frac{\partial}{\partial x^k} (\cdot - \cdot') 
\Phi)_p \,\, dp(g) + \\
&& + \int h (\Psi, {g'}^{ik} \frac{\partial}{\partial x^k} (\nabla'_i - \nabla_i) \Phi 
- (\nabla'_i - \nabla_i) {g'}^{ik} \frac{\partial}{\partial x^k} \cdot \Phi)_p + \non \\
&& + (\nabla'_i - \nabla_i) X \cdot' \Phi + \left( \nabla_i {g'}^{ik} 
\frac{\partial}{\partial x^k} \right) (\cdot' - \cdot) \Phi)_p \,\, dp(g).
\eea
(7.38) equals to
\be
\int h ( {\nabla'_i}^* \Psi, {g'}^{ik} \frac{\partial}{\partial x^k}
(\cdot - \cdot') \Phi)_p \,\, dp(g).
\end{equation}
If $\Phi$ is Sobolev and $\Psi = W$ then we obtain again by a u.l.f. cover by normal charts
$\{ U_\alpha \}_\alpha$ and an associated decomposition of unity $\{ \phi_\alpha \}_\alpha$
\bea
&& \int h (W,\eta^{op}_3 \Phi)_p \,\, dp(g) = \non \\
&& = \int h (W, \sum\limits_\alpha {g'}^{ik}_\alpha \frac{\partial}{\partial x^k}
(\cdot - \cdot') \nabla'_{\alpha, i} (\phi_\alpha \Phi))_p \,\, dp(g) = \non \\
&& = \int h ( {\nabla'_{\alpha,i}}^* W, \sum\limits_\alpha \phi_\alpha {g'}^{ik}_\alpha
\frac{\partial}{\partial x^k_\alpha} (\cdot - \cdot') \Phi)_p \,\, dp(g) + \\
&& + \int h (W, \eta^{op}_{3,0} \Phi)_p \,\, dp(g), 
\eea
where $\eta^{op}_{3,0} \Phi$ is the right component in $h(\cdot, \cdot)$ under the
integral (7.41), multiplied with $\phi_\alpha$ and summed up over $\alpha$.

Now we proceed literally as before. Start with 
\bea
(I_{3,3}) &=& \int\limits^t_{\frac{t}{2}} e^{-sD^2} D \eta^{op}_3 
e^{-(t-s) D^{'2}_{L_2}} \,\, ds =  \non \\
&=& \int\limits^t_{\frac{t}{2}} (D e^{-\frac{s}{2}D^2}) [(e^{-\frac{s}{4}D^2} f)
(f^{-1} e^{-\frac{s}{4}D^2} \eta^{op}_3)] e^{-(t-s) D^{'2}_{L_2}} \,\, ds.
\eea
We want that for suitable $f \in L_2$, 
$f^{-1} e^{\frac{s}{4}D^2} \eta^{op}_3$ 
is Hilbert--Schmidt. For this we have to estimate
$h(W(t,m,p),\eta^{op}_3 \cdot)_p$ 
and to show it defines an integral operator with finite
$L_2((M,E),dp)$--norm. We estimate
\bea
N(\Phi) &=& | \sk{\delta(m)}{e^{-tD^2} \eta^{op}_3 \Phi} |_{L_2((M,E),dp)} = \\
&=& | h(W(t,m,p), \eta^{op}_3 \Phi)_p |_{L_2((M,E),dp)}.
\eea
Using (7.41) and (7.42), we write
\bea
N(\Phi) &=& | h (W(t,m,p), \eta^{op}_3 \Phi)_p |_{L_2(dp)} = \non \\
&=& | h (\eta^{op}_{3,1} W(t,m,p), \eta^{op}_{3,0} \Phi)_p \non \\
&+& h (W(t,m,p), \eta^{op}_{3,0,0} \Phi)_p |_{L_2(dp)}.
\eea
Now we use
$|{\nabla'_X}^* \raisebox{0.3ex}{$\chi$}| \le 
C_1 |\nabla'_X \raisebox{0.3ex}{$\chi$}| \le 
C_2 |X| \cdot |\nabla' \raisebox{0.3ex}{$\chi$}| \le
C_3 |X|(|\nabla \raisebox{0.3ex}{$\chi$}| + |\raisebox{0.3ex}{$\chi$}|)$,
that the cover is u.l.f. and
$|\nabla W| \le C_4 (|D W| + |W|)$
and obtain
\bea
N(\Phi) & \le & C ( | h D W(t,m,p), \eta^{op}_{3,0} \Phi)_p |_{L_2(dp)}
+ | h (W(t,m,p), \eta^{op}_{3,0,0} \Phi)_p |_{L_2(dp)} = \non \\
&=& C ( N_1(\Phi) + N_2(\Phi)). \non
\eea
Here we again essentially use the bounded geometry and refer to [3] and [21].

Hence we have to estimate
\be
\sup\limits_{{\scriptsize \begin{array}{c} \Phi \in C^\infty_c(E) \\ |\Phi|_{L_2}=1 
\end{array}}}
N_1(\Phi) = 
\sup\limits_{{\scriptsize \begin{array}{c} \Phi \in C^\infty_c(E) \\ |\Phi|_{L_2}=1 
\end{array}}}
| \sk{\delta(m)}{(De^{-tD^2}) \eta^{op}_{3,0} \Phi} |_{L_2(dp)}
\end{equation}
and
\be
\sup\limits_{{\scriptsize \begin{array}{c} \Phi \in C^\infty_c(E) \\ |\Phi|_{L_2}=1 
\end{array}}}
N_2(\Phi) = 
\sup\limits_{{\scriptsize \begin{array}{c} \Phi \in C^\infty_c(E) \\ |\Phi|_{L_2}=1 
\end{array}}}
| \sk{\delta(m)}{e^{-tD^2} \eta^{op}_{3,0,0} \Phi} |_{L_2(dp)}.
\end{equation}
According to $k>r+1>n+3$, 5.4, 5.5, 5.7
\bea
&& D W(t,m,\cdot), W(t,m,\cdot) \in H^{\frac{r}{2} (E)} \non \\
&& |D W(t,m,\cdot)|_{H^{\frac{r}{2}}}, |W(t,m,\cdot)|_{H^{\frac{r}{2}}} \le C_1(t)
\eea
and we can restrict on (7.48), (7.49) to
\be
\sup\limits_{{\scriptsize \begin{array}{c} \Phi \in C^\infty_c(E) \\ |\Phi|_{L_2}=1 \\
|\Phi|_{H^{\frac{r}{2}}} \le C_1(t) \end{array}}}
\end{equation}
$\eta^{op}_{3,0}$, $\eta^{op}_{3,0,0}$ are of order zero and can be estimated by
\be
C_0 |\cdot - \cdot'|_{2,\frac{r}{2}} \le C_1 |\cdot - \cdot'|_{1,r-1}
\end{equation}
and
\be
D_0 \cdot ( |\nabla - \nabla'|_{2,\frac{r}{2}} + |\cdot - \cdot'|_{2,\frac{r}{2}} \le
D_1 \cdot ( |\nabla - \nabla'|_{2,r-1} + |\cdot - \cdot'|_{1,r-1}
\end{equation}
respectively.

Now we proceed literally as for $(I_{1,3})$, replacing (7.34) by
\be
N_1(\Phi) = \frac{1}{\sqrt{4 \pi t}} \frac{1}{2t} \left| \int\limits^{+\infty}_{-\infty}
s e^{-\frac{s^2}{4t}} e^{isD} \eta^{op}_{3,0} \Phi(m) \,\, ds \right| .
\end{equation}
$(I_{3,3})$ is done, $(I_{3,4})$, $(I_{3,1})$, $(I_{3,2})$ are absolutely parallel to
the case $\nu=1$.

This finishes the proof of 7.1.
\qed

We need in later sections the theorem analogous to 6.3 for the case of additional
variation of $g,\nabla^h,\cdot$.

\begin{theorem}
Suppose the hypothesises of 7.1. Then
\[ D e^{-tD^2} - D'_{L_2} e^{-tD^{'2}_{L_2}} \]
is of trace class and the trace norm is uniformly bounded on compact $t$--intervalls
$[a_0,a_1]$, $a_0>0$.
\end{theorem}

\begin{proof}
The proof is a simple combination of the proofs of 6.3 and 7.1.
\qed
\end{proof}

\begin{example}
The simplest standard example is 
$E = (\Lambda^* T^* M \otimes \C, g_{\Lambda^*}, \nabla^{g_{\Lambda^*}}) \lra (M^n,g)$
with Clifford multiplication
\[ X \otimes \omega \in T_m M \otimes \Lambda^* T^* M \otimes \C \lra X \cdot \omega
= \omega_X \wedge \omega - i_X \omega, \]
where $\omega_X := g(, X)$. In this case $E$ as a vector bundle remains fixed but the
Clifford module structure varies smoothly with $g,g' \in \comp(g)$. It is well known
that in this case $D=d+d^*$, $D^2=(d+d^*)^2$ Laplace operator $\Delta$.
\end{example}

\begin{theorem}
Assume $(M^n,g)$ with $(I)$, $(B_k)$, $k \ge r+1 > n+3$, 
$g' \in {\cal M}(I,B_k)$,
$g' \in \comp^{1,r+1}(g) \subset {\cal M}(I,B_k)$. Denote by
$\Delta'_{L_2(g)} = U^* i^* \Delta (g') i U$
the transformation of $\Delta' = \Delta (g')$ from
$L_2(g',g') \equiv L_2((M, \Lambda^* T^* M \otimes \C), g', g'_{\Lambda^*})$
to $L_2(g) \equiv L_2(g,g) = L_2 ((M, \Lambda^* T^* M \otimes \C), g, g_{\Lambda^*})$, 
where 
$i: L_2(g,g') = L_2 ((M, \Lambda^* T^* M \otimes \C), g',g_{\Lambda^*}) \lra L_2(g',g')$ 
and 
$U: L_2(g,g) \lra L_2(g',g)$,
$U = \alpha^{\frac{1}{2}}$,
$dq(g)=\alpha(q) \cdot dq(g')$, 
are the canonical maps, 
$i^*$, $U^*$
their adjoints. Then for $t>0$
\[ e^{-t \Delta(g)} - e^{-t \Delta'_{L_2}(g)} \]
is of trace class and the trace norm is uniformly bounded on compact $t$--intervalls
$[a_0,a_1]$, $a_0>0$.
\end{theorem}

Unfortunately the proof would not follow from 7.1 since 
$E' = E'(g') \notin \comp^{1,r+1}_{L,diff,F} (E)$ but
$E' \in \comp^{1,r+1}_{L,diff} (E)$, the fibre metric 
$g_{\Lambda^*}$ varies simultaneously with $g$. Now there are two possibilities.
1. A complete direct proof for this special case which will be even much easier
than the proof of 7.1 since $h',\nabla',\cdot'$ have very explicit expressions
depending only on $g'$. 2. A still more general version of 7.1, admitting even
variation of the fibre metric, 
$g \lra g', h \lra h', \nabla^h \lra \nabla^{h'}, \cdot \lra \cdot'$.
We decide to establish the general version. Before the formulation of the theorem
we must give some explanations. Consider the Hilbert spaces
$L_2(g,h) = L_2((M,E),g,h)$, 
$L_2(g',h) = L_2((M,E),g',h)$, 
$L_2(g',h') = L_2((M,E),g',h') \equiv L'_2$
and the maps
\bea
&& i_{(g',h),(g',h')} : L_2(g',h) \lra L_2(g',h'), \en 
i_{(g',h),(g',h')} \Phi = \Phi \non \\
&& U_{(g,h),(g',h)} : L_2(g,h) \lra L_2(g',h), \en 
U_{(g,h),(g',h)} \Phi = \alpha^{\frac{1}{2}} \Phi \non
\eea
where $dp(g) = \alpha(p) dp(g')$. Then we set
\bea
D'_{L_2(g,h)} &=& D'_{L_2} := U^*_{(g,h),(g',h)} i^*_{(g',h),(g',h')}
D' i_{(g',h),(g',h')} U_{(g,h),(g',h)} \equiv \non \\
& \equiv & U^* i^* D' i U.
\eea
Here $i^*$ is even locally defined (since $g'$ is fixed) and
$i^*_p = dual^{-1}_h \circ i' \circ dual_{h'}$, 
where $dual_h(\Phi(p)) = h_p (\cdot, \Phi(p))$.
In a local basis field $\Phi_1, \dots, \Phi_N$, 
$\Phi(p) = \xi^i(p) \Phi_i(p)$, 
\be
i^*_p \Phi(p) = h^{kl} h'_{ik} \xi^i \Phi_l(p).
\end{equation}
It follows from (7.55) that for $h' \in \comp^{1,r+1}(h)$
$i^*$, ${i^*}^{-1}$ are bounded up to order $k$, 
\bea
&& i^*-1, {i^*}^{-1}-1 \in \Omega^{0,1,r+1} (Hom((E,h',\nabla^{h'}) \lra \non \\
&& \lra (M,g'), (E,h,\nabla^h) \lra (M,g')))
\eea
and
\bea
&& i^*-1, {i^*}^{-1}-1 \in \Omega^{0,2,\frac{r+1}{2}} (Hom((E,h',\nabla^{h'}) \lra \non \\
&& \lra (M,g'), (E,h,\nabla^h) \lra (M,g'))).
\eea
$D' \equiv \overline{D'}$ is self adjoint on 
$D_{\overline{D'}} = \overline{C^\infty_c(E)}^{|\,\,|_{D'}}$,
where
$|\Phi|^2_{D'} = |\Phi|^2_{L'_2} + |D' \Phi|^2_{L'_2}$.
$i: L_2(g',h) \lra L_2(g',h') \equiv L'_2$ 
and
$i^*: L_2(g',h') \lra L_2(g',h)$ 
are for 
$h' \in \comp^{1,r+1}(h)$
quasi isometries with bounded derivatives, they map
$C^\infty_c (E)$
1--1 onto 
$C^\infty_c (E)$
and
$i^* D' i$
is self adjoint on
$\overline{C^\infty_c(E)}^{|\,\,|_{i^*D'i}} = D_{i^*D'i} \subset L_2((M,E),g',h)
\equiv L_2(g',h)$.
We obtain as a consequence that
$e^{-t(i^*D'i)^2}$ 
is defined and selfadjoint in 
$L_2((M,E),g',h) = L_2(g',h)$,
maps for $t>0$ and $i,j \in \Z$
$H^i (E,i^*D'i)$ 
continuously into
$H^j (E,i^*D'i)$ 
and has the heat kernel
$W'_{g',h}(t,m,p) = \sk{\delta(m)}{e^{-t(i^*D'i)^2} \delta(p)}$,
$W'(t,m,p)$
satisfies the same general estimates as 
$W(t,m,p)$
in section 5. By exactly the same arguments we obtain that
$e^{-t U^* (i^*D'i)^2U} = e^{-t (U^* i^*D'i U)^2} = U^* e^{-t (i^*D'i)^2} U$
is defined in
$L_2 = L_2((M,E),g,h)$,
self adjoint and has the heat kernel
$W'_{L_2} (t,m,p) = W'_{g,h} (t,m,p) = \alpha^{-\frac{1}{2}}(m) W'_{g',h} (t,m,p) 
\alpha (p)^{\frac{1}{2}}$.
Here we assume
$g' \in \comp^{1,r+1}(g)$.
Now we are able to formulate our main theorem.

\begin{theorem}
Let $E=((E,h,\nabla = \nabla^h,\cdot) \lra (M^n,g))$
be a Clifford bundle with 
$(I)$, $(B_k(M,g))$, $(B_k(E,\nabla))$, $k \ge r+1 >n+3$, 
$E' = ((E,h,\nabla'=\nabla^{h'},\cdot') \lra (M^n,g)) \in 
\gencomp^{1,r+1}_{L,diff} (E) \cap \cl {\cal B}^{N,n} (I,B_k)$, 
$D=D(g,h,\nabla=\nabla^h,\cdot)$,
$D'=D(g',h,\nabla'=\nabla^{h'},\cdot')$
the associated generalized Dirac operators, 
$dp(g) = \alpha (p) dp(g')$,
$U = \alpha^{\frac{1}{2}}$.
Then for $t>0$
\be
e^{-tD^2} - U^* e^{-t(i^* D' i)^2} U
\end{equation}
is of trace class and the trace norm is uniformly bounded on compact 
$t$--intervalls $[a_0,a_1]$, $a_0>0$.
\end{theorem}

\begin{proof}
We are done if we could prove the assertions for
\be
e^{-t(U D' U^*)^2} - e^{-t(i^* D' i)^2} = U e^{-tD^2} U^* - e^{-t(i^* D' i)^2}
\end{equation}
since $U^* (7.59) U = (7.58)$. To get a better explicit expression for (7.59), 
we apply again Duhamel's principle. This holds since Greens formula for 
$U D^2 U^*$ holds, 
\[ \int h_q (U D^2 U^* \Phi, \Psi) - h (\Phi, U D^2 U^* \Psi) \,\, dq(g')  = 0. \]
We obtain
\bea
&& - \int\limits^t_0 \int\limits_M h_q ( \alpha^{\frac{1}{2}} (m) W(s,m,q)
\alpha^{-\frac{1}{2}} (q), \left( U D^2 U^* + \frac{\partial}{\partial t} \right) 
 W'_{g',h} (t-s,q,p) \,\, dq(g') \,\, ds = \non \\
&& = - \int\limits^t_0 \int\limits_M h_q ( \alpha^{\frac{1}{2}} (m) W(s,m,q)
\alpha^{-\frac{1}{2}} (q), ( U D^2 U^* - (i^* D' i)^2 ) 
 W'_{g',h} (t-s,q,p) \,\, dq(g') \,\, ds = \non \\
&& = \alpha^{\frac{1}{2}} (m) W(s,m,q) \alpha^{-\frac{1}{2}} (q)
- W'_{g',h'} (t,m,p) = \non \\
&& =W_{g',h} (t,m,p) -  W'_{g',h} (t,m,p).
\eea
(7.60) expresses the operator equation
\bea
&& e^{-t(U D U^*)^2} - e^{-t(i^* D' i)^2} = \non \\
&=& - \int\limits^t_0 e^{-s(U^* D U)^2} ((U D U^*)^2 - (i^* D' i)^2) 
e^{-(t-s)(i^* D' i)^2} \,\, ds = \non \\
&=& - \int\limits^t_0 e^{-s(U D U^*)^2} U D U^* (U D U^* - i^* D' i) 
e^{-(t-s)(i^* D' i)^2} \,\, ds - \\
&-& \int\limits^t_0 e^{-s(U D U^*)^2} (U D U^* - i^* D' i) (i^* D' i) 
e^{-(t-s)(i^* D' i)^2} \,\, ds .
\eea
We write (7.62) as 
\bea
&& - \int\limits^t_0 \alpha^{\frac{1}{2}} e^{-sD^2} D \alpha^{-\frac{1}{2}}
(\alpha^{\frac{1}{2}} D \alpha^{-\frac{1}{2}} - i^* D' i ) e^{-(t-s)(i^*D'i)^2} \,\, ds = \non \\
&& - \int\limits^t_0 \alpha^{\frac{1}{2}} e^{-sD^2} D \alpha^{-\frac{1}{2}}
(D - i^* D' i - \frac{\grad \alpha \cdot}{2 \alpha}) e^{-(t-s)(i^*D'i)^2} \,\, ds = \non \\
&&  - \int\limits^t_0 \alpha^{\frac{1}{2}} e^{-sD^2} D \alpha^{-\frac{1}{2}}
i^* ((i^{*-1}-1) D + (D-D') - i^{*-1} \frac{\grad \alpha \cdot}{2\alpha}) 
e^{-(t-s)(i^*D'i)^2} \,\, ds = \non \\
&& \int\limits^t_0 \alpha^{\frac{1}{2}} e^{-sD^2} D (\eta_0+\eta_1+\eta_2+\eta_3+\eta_4)
e^{-(t-s)(i^*D'i)^2} \,\, ds , \non 
\eea
$\eta_0 = \frac{\grad \alpha \cdot}{2\alpha^{\frac{3}{2}}}$, 
$\eta_i = - \alpha^{-\frac{1}{2}} i^* \eta_i (7)$, $i=1,2,3$, 
$\eta_1(7)=(7.10)$, $\eta_2(7)=(7.11)$, $\eta_3(7)=(7.12)$, 
$\eta_4=\alpha^{-\frac{1}{2}} i^{*-1} (i^*-1)D$.
Here $\eta_0$ and $\eta_2$ are of zeroth order. $\eta_1$ and $\eta_3$
can be discussed as in (7.18)--(7.54). $\eta_4$ can be discussed analogous to 
$\eta_1$, $\eta_3$ in section 7, i.e. $\eta_4$ will be shifted via partial integration
to the left (up to zero order terms) and 
$\alpha^{-\frac{1}{2}} i^* (i^*-1)$
thereafter again to the right. In the estimates one has to replace $W$ by $D W$ and
nothing essentially changes as we exhibited in (7.34). We perform in (7.62) the same
decomposition and have to estimate 20 integrals,
\def\theequation{$I_{\nu,1}$}   
\be
 \int\limits^{\frac{t}{2}}_0 \alpha^{\frac{1}{2}} e^{-sD^2} D \eta_\nu 
 e^{-(t-s)(i^*D'i)^2} \,\, ds ,
\end{equation}
\def\theequation{$I_{\nu,2}$}  
\be
 \int\limits^{\frac{t}{2}}_0 \alpha^{\frac{1}{2}} e^{-sD^2} \eta_\nu 
 (i^* D' i) e^{-(t-s)(i^*D'i)^2} \,\, ds ,
\end{equation}
\def\theequation{$I_{\nu,3}$}  
\be
 \int\limits^t_{\frac{t}{2}} \alpha^{\frac{1}{2}} e^{-sD^2} D \eta_\nu  
 e^{-(t-s)(i^*D'i)^2} \,\, ds ,
\end{equation}
\def\theequation{$I_{\nu,4}$}  
\be
 \int\limits^t_{\frac{t}{2}} \alpha^{\frac{1}{2}} e^{-sD^2} \eta_\nu 
 (i^* D' i) e^{-(t-s)(i^*D'i)^2} \,\, ds, 
\end{equation}                         
\def\theequation{\thesection.\arabic{equation}}  
\setcounter{equation}{62}

$\nu=0,\dots,4$ and to show that these are products of Hilbert--Schmidt operators
and have uniformly bounded trace norm on compact $t$--intervals. This has been
completely modelled in the proof of 7.1.
\qed
\end{proof}

Combining the proof of 7.9 and the decomposition (6.45)--(6.48), we obtain

\begin{theorem}
Assume the hypothesises of 7.10. Then for $t>0$
\[ e^{tD^2} D - U^* e^{-t(i^*D'i)^2} (i^*D'i) U \]
is of trace class and its trace norm is uniformly bounded on compact
$t$--intervalls $[a_0,a_1]$, $a_0>0$.
\qed
\end{theorem}

The operators $i^* {D'}^2 i$ and $(i^* D' i)^2$ are different in general. We 
should still compare $e^{-t i^* {D'}^2 i}$ and $e^{-t(i^*D'i)^2}$.

\begin{theorem}
Assume the hypothesises of 7.10. Then for $t>0$
\[ e^{-t(i^*{D'}^2i)} - e^{-t(i^*D'i)^2} \]
is of trace class and the trace norm is uniformly bounded on compact
$t$--intervalls $[a_0,a_1]$, $a_0>0$.
\end{theorem}

\begin{proof}
We obtain again immediately from Duhamel's principle
\bea
&& e^{-t i^* {D'}^2 i} - e^{-t(i^*D'i)^2} = \non \\
&& = - \int\limits^t_0 e^{-s (i^* {D'}^2 i)} (i^* {D'}^2 i - (i^* D' i)^2)
e^{-(t-s)(i^*D'i)^2} \,\, ds = \non \\
&& = - \int\limits^t_0 e^{-s (i^* {D'}^2 i)} i^* D' (1 - ii^*) D' i 
e^{-(t-s)(i^*D'i)^2} \,\, ds = \non \\
&& = - \int\limits^t_0 e^{-s (i^* {D'}^2 i)} (i^* D' i) i^{-1} (1 - ii^*) i^{*-1} 
(i^* D' i) e^{-(t-s)(i^*D'i)^2} \,\, ds .
\eea
In $[\frac{t}{2}, t]$ we shift $i^*D'i$ again to the left of the kernel
$W'_{e^{-s (i^* D^2 i)}}$ via partial integration and estimate
\bea
&& (i^*D'i e^{-\frac{s}{2} (i^* {D'}^2 i)}) 
[(e^{-\frac{s}{4} (i^* {D'}^2 i)}) f) (f^{-1} e^{-\frac{s}{4} (i^* {D'}^2 i)} 
i^{-1}  (1-ii^*) i^{*-1})] \non \\
&& ((i^*D'i) e^{-(t-s) (i^*D'i)^2} ) \non
\eea
as before. In $[0,\frac{t}{2}]$ we write the integrand of (7.63) as
\bea
&& (e^{-s (i^* {D'}^2 i)} i^*D'i) [((i^*i)^{-1} e^{-\frac{t-s}{4} (i^*D'i)^2} f^{-1})
(f e^{-\frac{t-s}{4} (i^*D'i)^2} )] \non \\
&& (e^{-\frac{t-s}{2} (i^*D'i)^2} (i^*D'i)) \non 
\eea
and proceed as in the corresponding cases.
\qed
\end{proof}

\begin{theorem}
Assume the hypothesises of 7.10. Then for $t>0$
\[ e^{-tD^2} - e^{-t(U^*i^*D^2iU)} \equiv e^{-tD^2} - U^* e^{-t(i^*D^2i)} U \]
is of trace class and the trace norm is uniformly bounded on any
$t$--intervall $[a_0,a_1]$, $a_0>0$.
\end{theorem}

\begin{proof}
This immediately follows form 7.10 and 7.12.
\qed
\end{proof}

Theorem 7.9 is now a special case of the more general theorem 7.12. 
\qed

Finally the last class of admitted perturbations are compact topological
perturbations which will be the content of the next section.

\section{Trace class property in the class of additional compact topological
perturbations}

\setcounter{equation}{0}

Let $E=((E,h,\nabla^h) \lra (M^n,g)) \in \cl {\cal B}^{N,n} (I,B_k)$ be a
Clifford bundle, $k \ge r+1>n+3$, 
$E'=((E,h',\nabla^{h'}) \lra ({M'}^n,{g'})) \in \comp^{1,r+1}_{L,diff,rel} (E) \cap
\cl {\cal B}^{N,n} (I,B_k)$. Then there exist $K \subset M$, $K' \subset M'$
and a vector bundle isomorphism (not necessarily an isometry)
$f=(f_E,f_M) \in \tilde{\cal D}^{1,r+2} (E|_{M \setminus K}, E'|_{M' \setminus K'}$
s. t.
\bea
&& g|_{M \setminus K} \mbox{ and } f^*_M g'|_{M \setminus K} 
\mbox{ are quasi isometric,} \\
&& h|_{E|_{M \setminus K}} \mbox{ and } f^*_E h'|_{E|_{M \setminus K}}
\mbox{ are quasi isometric,} \\
&& | g|_{M \setminus K} - f^*_M g'|_{M \setminus K} |_{g,1,r+1} < \infty ,\\
&& | h|_{E_{M \setminus K}} - f^*_E h'|_{_E|{M \setminus K}} |_{g,h,\nabla^h,1,r+1} < \infty ,\\
&& | \nabla^h|_{E|_{M \setminus K}} - f^*_E \nabla^{h'}|_{E|_{M \setminus K}} 
|_{g,h,\nabla^h,1,r+1} < \infty ,\\
&& | \cdot|_{M \setminus K} - f^*_E \cdot'|_{M \setminus K} |_{g,h,\nabla^h,1,r+1} < \infty. 
\eea
(8.1) -- (8.6) also hold if we replace $f$ by $f^{-1}$, $M \setminus K$
by $M' \setminus K'$ and $g,h,\nabla^h,\cdot$ by $g',h',\nabla^{h'},\cdot'$.
If we consider the complete pull back
$f^*_E(E'|_{M' \setminus K'})$,
i.e. the pull back together with all Clifford datas, then we have on $M \setminus K$
two Clifford bundles, $E|_{M \setminus K}$, $f^*_E(E'|_{M' \setminus K'}$ which are
as vector bundles isomorphic and we denote $f^*_E(E'|_{M' \setminus K'}$
again by $E'$ on $M \setminus K$, i.e.
$g'_{new} = (f_M|_{M \setminus K})^* g'_{old}$ etc..
(8.1) -- (8.6) and the symmetry of our uniform structure
${\mathfrak U}^{1,r+1}_{L,diff,rel}$ imply
\bea
&& W^{1,i} (E|_{M \setminus K}) \cong W^{1,i} (E'|_{M \setminus K}), 
\quad 0 \le i \le r+1 ,\non \\
&& W^{2,j} (E|_{M \setminus K}) \cong W^{2,j} (E'|_{M \setminus K}), 
\quad 0 \le j \le \frac{r+1}{2}, \non \\
&& H^j (E|_{M \setminus K},D) \cong H^j (E'|_{M \setminus K},D'), 
\quad 0 \le j \le \frac{r+1}{2},  
\eea
\bea
&& \Omega^{1,1,i} (End (E|_{M \setminus K})) \cong \Omega^{1,1,i} (End (E'|_{M \setminus K})), 
\quad 0 \le i \le r+1, \non \\
&& \Omega^{1,2,j} (End (E|_{M \setminus K})) \cong \Omega^{1,2,j} (End (E'|_{M \setminus K})), 
\quad 0 \le j \le \frac{r+1}{2} . \non 
\eea
Here the Sobolev spaces are defined by restriction of corresponding Sobolev sections.

We now fix our set up for compact topological pertulations. Under more special
assumptions this has already been done in [2], [5].
Set ${\cal H} = L_2((K,E|_K),g,h) \oplus L_2((K',E'|_{K'}),g',h') \oplus
L_2((M \setminus K, E),g,h) $
and consider the following maps
\bea
&& i_{L_2,K'} : L_2((K',E'|_{K'}),g',h') \lra {\cal H}, \non \\
&& i_{L_2,K'} (\Phi) = \Phi, \non \\
&& i^{-1} : L_2((M' \setminus K',E'|_{M' \setminus K'}),g',h') \lra 
L_2((M' \setminus K',E'|_{M' \setminus K'}),g',h), \non \\
&& i^{-1} \Phi = \Phi, \non \\
&& U^* : L_2((M' \setminus K',E'|_{M' \setminus K'}),g',h) \lra 
L_2((M' \setminus K',E'|_{M' \setminus K'}),g,h), \non \\
&& U^* \Phi = \alpha^{-\frac{1}{2}},  \non
\eea
where $dq(g)=\alpha(q) dq(g')$. 
We identify $M \setminus K$ and $M' \setminus K'$ as manifolds and
$E'|_{M' \setminus K'}$ and $E|_{M \setminus K}$ as vector 
bundles. Then we have natural embeddings
\bea
&& i_{L_2,M} : L_2((M,E),g,h) \lra {\cal H} , \non \\
&& i_{L_2,K'} \oplus U^*i^{-1} : L_2((M',E'),g',h') \lra {\cal H} , \non \\
&& (i_{L_2,K'} \oplus U^*i^{-1}) \Phi = i_{L_2,K'} \raisebox{0.3ex}{$\chi$}_{K'} \Phi +
U^*i^{-1} \raisebox{0.3ex}{$\chi$}_{M' \setminus K'} \Phi. \non
\eea
The images of these two embeddings are closed subspaces of ${\cal H}$. Denote
by $P$ and $P'$ the projection onto these closed subspaces. $D$ is defined on
${\cal D}_D \subset \im P$. We extend it onto $(\im P)^\perp$ as zero
operator. The definition of (the shifted) $D'$ is a little more complicated.
For the sake of simplicity of notation we write
$U^*i^{-1} \equiv i_{L_2,K'} \oplus U^*i^{-1} = id \oplus U^*i^{-1}$, 
keeping in mind that $i_{L_2,K'}$ fixes $\raisebox{0.3ex}{$\chi$}_K, \Phi$ and the scalar product.
Moreover we set also 
$i U \raisebox{0.3ex}{$\chi$}_{K'} \Phi = U^* i^* \raisebox{0.3ex}{$\chi$}_{K'} \Phi = \raisebox{0.3ex}{$\chi$}_{K'} \Phi$.
Let $\Phi \in {\cal D}_{D'}$, 
$\raisebox{0.3ex}{$\chi$}_{K'} \Phi + U^* i^{-1} \raisebox{0.3ex}{$\chi$}_{M' \setminus K'} \Phi$ 
its image in ${\cal H}$. Then 
$(U^*i^*D'iU) (\raisebox{0.3ex}{$\chi$}_{K'} \Phi + U^* i^{-1} \raisebox{0.3ex}{$\chi$}_{M' \setminus K'}, \Phi) =
U^*i^*D'\Phi = \raisebox{0.3ex}{$\chi$}_{K'} D' \Phi + U^*i^* \raisebox{0.3ex}{$\chi$}_{M' \setminus K'} D' \Phi$.
Now we set as 
${\cal D}_{U^*i^*D'iU} \subset H$ 
\be
{\cal D}_{U^*i^*D'iU} = \{ \raisebox{0.3ex}{$\chi$}_{K'} \Phi + U^* i^{-1} \raisebox{0.3ex}{$\chi$}_{M' \setminus K'} \Phi
| \Phi \in {\cal D}_{D'} \} \oplus (\im P')^\perp.
\end{equation}
It follows very easy from the selfadjointness of $D'$ on ${\cal D}_{D'}$ and (8.7)
that $U^*i^*D'iU$ is self adjoint on ${\cal D}_{U^*i^*D'iU}$, if we additionally
set $U^*i^*D'iU=0$ on $(\im P')^\perp$.

\begin{remark}
If $g$ and $h$ do not vary then we can spare the whole $i-U$--procedure, 
$i=U=id$. Nevertheless this case still includes interesting pertulations.
Namely pertulations of $\nabla,\cdot$ and compact topological pertulations.

\qed
\end{remark}

We set for the sake of simplicity
$\tilde{D'} = U^*i^*D'iU$.
The first main result of this section is the following

\begin{theorem}
Let $E=((E,h,\nabla^h) \lra (M^n,g)) \in \cl {\cal B}^{N,n}(I,B_k)$,
$k \ge r+1>n+3$, 
$E' \in \gencomp^{1,r+1}_{L,diff,rel}(E) \cap \cl {\cal B}^{N,n}(I,B_k)$.
Then for $t>0$
\be
e^{-tD^2}P - e^{-t{\tilde{D'}}^2}P'           
\end{equation}
and
\be
e^{-tD^2}D - e^{-t{\tilde{D'}}^2} \tilde{D'}
\end{equation}
are of trace class and their trace norms are uniformly bounded on any
$t$--intervall $[a_0,a_1]$, $a_0>0$.
\end{theorem}

For the proof we make the following construction. Let
$V \subset M \setminus K$ be open,
$\overline{M \setminus K} \setminus V$ compact,
$\dist (V, \overline{M \setminus K} \setminus (M \setminus K)) \ge 1$
and denote by $B \in L({\cal H})$ the multiplication operator
$B=\raisebox{0.3ex}{$\chi$}_\nu$. The proof of 8.2 consists of two steps. First we prove
8.2 for the restriction of (8.9), (8.10) to $V$, i.e. for $B(8.9)B$,
thereafter for $(1-B)(8.9)B$, $B(8.9)(1-B)$ and the same for (8.10).

\begin{theorem}
Assume the hypothesis of 8.2. Then
\bea
&& B(e^{-tD^2}P - e^{-t{\tilde{D'}}^2}P')B, \\
&& B(e^{-tD^2}D - e^{-t{\tilde{D'}}^2}\tilde{D'})B, \\
&& B(e^{-tD^2}P - e^{-t{\tilde{D'}}^2}P')(1-B), \\
&& (1-B)(e^{-tD^2}P - e^{-t{\tilde{D'}}^2}P')B, \\
&& B(e^{-tD^2}D - e^{-t{\tilde{D'}}^2}\tilde{D'})(1-B), \\
&& (1-B)(e^{-tD^2}D - e^{-t{\tilde{D'}}^2}\tilde{D'})B, \\
&& (1-B)(e^{-tD^2}P - e^{-t{\tilde{D'}}^2}P')(1-B), \\
&& (1-B)(e^{-tD^2}D - e^{-t{\tilde{D'}}^2}\tilde{D'})(1-B)
\eea
are of trace class and their trace norms are uniformly bounded on any
$t$--intervall $[a_0,a_1]$, $a_0>0$.
\end{theorem}

8.2 immediately follows from 8.3.
We start with the assertion for (8.11). Introduce functions
$\phi, \psi, \gamma \in C^{\infty}(M,[0,1])$
with the following properties.

1. $\supp \phi \subset M \setminus K$, 
$(1-\phi) \in C^\infty_c (M \setminus K)$, 
$\phi|_V=1$.

2. $\psi$ with the same properties as $\phi$ and additionally
$\psi=1$ on $\supp \phi$, i.e. 
$\supp(1-\psi) \cap \supp \phi=0$.

3. $\gamma \in C^\infty_c(M)$, $\gamma=1$ on 
$\supp(1-\phi)$, $\gamma|_V=0$.

Define now as in [5] an approximate heat kernel $E(t,m,p)$ on $M$ by
\[ E(t,m,p) := \gamma(m) W(t,m,p) (1-\phi(p)) + \psi(m) \tilde{W'}(t,m,p) \phi(p). \]
Applying Duhamel's principle yields
\bea
&& -\int\limits^\beta_\alpha \int\limits_M h_q (W(s,m,q), 
\left( \frac{\partial}{\partial t} + D^2 \right) E(t-s,q,p)) 
\raisebox{0.3ex}{$\chi$}_\nu(p) \,\, dq(g) \,\, ds \non \\
&& = \int\limits_M [h_q (W(\beta,m,q), E(t-\beta,q,p)) - h_q (W(\alpha,m,q), \non \\
&& E(t-\alpha,q,p))] \raisebox{0.3ex}{$\chi$}_\nu(p) \,\, dq(g).
\eea
Performing $\alpha \lra 0^+$, $\beta \lra t^-$ in (8.19), we obtain
\bea
&& -\int\limits^\beta_\alpha \int\limits_M h_q (W(s,m,q), 
\left( D^2 + \frac{\partial}{\partial t} \right) E(t-s,q,p)) 
\raisebox{0.3ex}{$\chi$}_\nu(p) \,\, dq(g) \,\, ds = \non \\
&& \lim\limits_{\beta \ra t^-} \int\limits_M [h_q (W(\beta,m,q), E(t-\beta,q,p)) 
\raisebox{0.3ex}{$\chi$}_\nu(p) \,\, dq(g) - \non \\ 
&& E(t,m,p) \raisebox{0.3ex}{$\chi$}_\nu(p).
\eea
Now we use
\be
\raisebox{0.3ex}{$\chi$}_V (p) (1-\phi(p)) = 0
\end{equation}
and obtain
\bea
&& \lim\limits_{\beta \ra t^-} \int\limits_M [h_q (W(\beta,m,q), E(t-\beta,q,p)) 
\raisebox{0.3ex}{$\chi$}_\nu(p) \,\, dq(g) = \non \\
&& = \lim\limits_{\beta \ra t^-} \int\limits_M [h_q (W(\beta,m,q), \psi(q) \tilde{W'}(t-\beta,q,p)) 
\phi(p) \raisebox{0.3ex}{$\chi$}_\nu(p) \,\, dq(g) = \non \\
&& = W(t,m,p) \non
\eea
since $\tilde{W'}(\tau,q,p)$ is the heat kernel of $e^{-\tau \tilde{D'}^2}$.
This yields
\bea
&& - \int\limits^t_0 \int\limits_M h_q (W(s,m,q), \left( D^2+ \frac{\partial}{\partial t} \right)
E(t-s,q,p)) \raisebox{0.3ex}{$\chi$}_V(p) \,\, dq(g) \,\, ds = \non \\
&& - \int\limits^t_0 \int\limits_M h_q (W(s,m,q), (D^2\psi(q)-\psi(q){\tilde{D'}}^2)
\tilde{W'}(t-s,q,p)) \raisebox{0.3ex}{$\chi$}_V(p) \,\, dq(g) \,\, ds = \non \\
&& = [W(t,m,p) - \tilde{W'}(t,m,p)] \cdot \raisebox{0.3ex}{$\chi$}_V(p) .
\eea
(8.22) expresses the operator equation
\be
(e^{-tD^2}P - e^{-t{\tilde{D'}}^2}P')B = 
- \int\limits^t_0 e^{-sD^2} (D^2 \psi - \psi \tilde{D'}^2) e^{-(t-s){\tilde{D'}}^2} B \,\, ds
\end{equation}
in ${\cal H}$ at kernel level.

We rewrite (8.23) as in the foregoing cases.
\bea
(8.23) &=& - \int\limits^t_0 e^{-sD^2} (D(D-\tilde{D'}) \psi + 
(D-\tilde{D'})\tilde{D'} \psi
+ \tilde{D'}^2 \psi - \psi \tilde{D'}^2) e^{-(t-s){\tilde{D'}}^2} \,\, ds = \non \\
& = &  - \left[ \int\limits^{\frac{t}{2}}_0 e^{-sD^2} (D(D-\tilde{D'}) \psi
e^{-(t-s){\tilde{D'}}^2} \,\, ds + \right. \\
&& + \int\limits^{\frac{t}{2}}_0 e^{-sD^2} (D-\tilde{D'})\tilde{D'} \psi
e^{-(t-s){\tilde{D'}}^2} \,\, ds + \\
&& + \int\limits^{\frac{t}{2}}_0 e^{-sD^2} (\tilde{D'}^2 \psi - \psi \tilde{D'}^2) \,\, ds \\
&& + \int\limits^t_{\frac{t}{2}} e^{-sD^2} D(D-\tilde{D'}) \psi
e^{-(t-s){\tilde{D'}}^2} \,\, ds + 
\eea
\bea
&& + \int\limits^t_{\frac{t}{2}} e^{-sD^2} (D-\tilde{D'}) \tilde{D'} \psi
e^{-(t-s){\tilde{D'}}^2} \,\, ds + \\
&& \left. + \int\limits^t_{\frac{t}{2}} e^{-sD^2} (\tilde{D'}^2 \psi - \psi \tilde{D'}^2)
e^{-(t-s){\tilde{D'}}^2} \,\, ds \right] . 
\eea 
Write the integrand of (8.27) as 
\[ (e^{-\frac{s}{2}D^2}D) [ ( e^{-\frac{s}{4}D^2} f)(f^{-1} e^{-\frac{s}{4}D^2} 
(D - \tilde{D'}) \psi ) ] e^{-(t-s)} \tilde{D'}^2 , \]
$|e^{-\frac{s}{2}D^2}D|_{op} \le \frac{C}{\sqrt{s}}$,
$|e^{-(t-s){\tilde{D'}}^2}|_{op} \le C'$
and $[\dots]$ is the product of two Hilbert--Schmidt operators if $f$ can be choosen
$\in L_2$ and such that 
$f^{-1} e^{-\frac{s}{4}D^2} (D-\tilde{D'}) \psi$
is Hilbert--Schmidt. We know from sections 6 and 7, sufficient for this is that
$(D-\tilde{D'}) \psi$ has Sobolev coefficients of order $r+1$ (and $p=1$).
\bea
(D-\tilde{D'}) \psi &=& 
(D-\alpha^{-\frac{1}{2}} i^* D' i \alpha^{\frac{1}{2}}) \psi \non \\
&=& \left( D - i^* \frac{\grad' \alpha}{2 \alpha} \cdot' -i^* D' \right) \psi \non \\
&=& i^* \left( (i^{*-1}-1)D + (D-D') - \frac{\grad' \alpha \cdot'}{2\alpha} \right) \psi \non \\
&=& i^* \Bigg[ i^{*-1} ( \grad \psi \cdot + \psi D) + \grad \Psi ( \cdot - \cdot') \non \\
&+& (\grad \psi - \grad' \psi) \cdot' 
+ \psi (D-D') - \frac{\grad' \alpha \cdot'}{2 \alpha} \psi \Bigg] . \non 
\eea
$i^*$ is bounded up to order $k$, $i^{*-1}-1$ is $(r+1)$--Sobolev, 
$\grad \psi$, $\grad' \psi$ have compact support, $0 \le \psi \le 1$, 
$\frac{\grad \alpha \cdot'}{2 \alpha}$ is $(r+1)$--Sobolev and 
$\psi (D-D')$ is completely discussed in (7.9) -- (7.53). Hence (8.27)
is completely done.

Write the integrand of (8.28) as
\[ [ (e^{\frac{s}{2}D^2} f)(f^{-1} e^{\frac{s}{2}D^2} (D-D')) ] 
(D' e^{-(t-s) \tilde{D'}^2}). \] 
$[\dots]$ is the product of two Hilbert--Schmidt operators with bounded trace norm on 
$t$--intervalls $[a_0,a_1]$, $a_0>0$. An easy calculation yields
\[ \tilde{D'} \psi = \alpha^{-\frac{1}{2}} i^* D' i \alpha^{\frac{1}{2}} \psi =
i^* \grad \psi \cdot' + \psi \tilde{D'},  \]
hence
\[ |\tilde{D'} \psi e^{-(t-s){\tilde{D'}}^2}|_{op} = |(i^* \grad' \psi \cdot' + \psi
\tilde{D'}) e^{-(t-s){\tilde{D'}}^2}|_{op} \le  C + \frac{C'}{\sqrt{t-s}} , \]
(8.28) is done.

Rewrite finally the integrands of (8.24), (8.25) as
\bea
&& (e^{-sD^2}D) [ ((D-\tilde{D'}) \psi e^{-\frac{t-s}{2}{\tilde{D'}}^2}f^{-1})
(f e^{-\frac{t-s}{2}{\tilde{D'}}^2}) ] = \non \\
&& = e^{-sD^2}D i^* [ (((i^{*-1}-1) (\grad \psi \cdot + \psi D) + \non \\
&& + \grad \psi
(\cdot - \cdot') + (\grad \psi - \grad' \psi) \cdot' + \non \\
&& + \psi (D-D') - \frac{\grad \alpha}{2 \alpha} \cdot' \psi)
e^{-\frac{t-s}{2}{\tilde{D'}}^2} f^{-1})(f e^{-\frac{t-s}{2}{\tilde{D'}}^2} )] \non 
\eea
and 
\bea
&& e^{-sD^2} i^* [ ((D-\tilde{D'}) \psi e^{-\frac{t-s}{4}{\tilde{D'}}^2}f^{-1})
(f e^{-\frac{t-s}{4}{\tilde{D'}}^2}) ] (\tilde{D'} e^{-\frac{t-s}{2}{\tilde{D'}}^2}) = \non \\
&& = e^{-sD^2} i^* [ (((i^{*-1}-1) (\grad \psi \cdot + \psi D) + \grad \psi
(\cdot - \cdot') + \non \\
&& (\grad \psi - \grad' \psi) \cdot' + \psi (D-D') - 
\frac{\grad \alpha}{2 \alpha} \cdot' \psi) \non \\
&& e^{-\frac{t-s}{4}{\tilde{D'}}^2} f^{-1}) (f e^{-\frac{t-s}{4}{\tilde{D'}}^2} )]
(\tilde{D'} e^{-\frac{t-s}{2}{\tilde{D'}}^2}) ,  \non 
\eea
respectively, and (8.24), (8.25) are done. The remaining integrals are 
(8.26) and (8.29). We have to find an appropriate expression for 
$\tilde{D'}^2 \psi - \psi \tilde{D'}^2$.
\bea
\tilde{D'}^2 &=& ( \alpha^{-\frac{1}{2}} i^* D' i \alpha^{\frac{1}{2}} )
( \alpha^{-\frac{1}{2}} i^* D' i \alpha^{\frac{1}{2}} ) \\
&=& \alpha^{-\frac{1}{2}} i^* D' i^* \left( \frac{\grad' \alpha}{2\alpha^{\frac{1}{2}}}
\cdot' + \alpha^{\frac{1}{2}} D' \right)  \non \\
&=&
i^* \left( D' \alpha^{-\frac{1}{2}} + \frac{\grad' \alpha}{2\alpha^{\frac{3}{2}}}
\cdot' \right) i^* \left( \frac{\grad' \alpha}{2\alpha^{\frac{1}{2}}}
\cdot' + \alpha^{\frac{1}{2}} D' \right)  \non \\
&=& i^* D' i^* D' + i^* D' i^* \frac{\grad' \alpha}{2\alpha} \cdot' + 
i^* \frac{\grad' \alpha}{2\alpha} \cdot' i^* \frac{\grad' \alpha}{2\alpha} \cdot' +
 \non \\
&& \quad + i^* \frac{\grad' \alpha}{2\alpha} \cdot' i^* D' .
\eea
Hence
\bea
&& \tilde{D'}^2 \psi - \psi \tilde{D'}^2 = i^* D' i^* D' \psi - \psi i^* D' i^* D' + \non \\
&& + i^* D' i^* \frac{\grad' \alpha}{2 \alpha} \cdot' \psi - \psi i^*  D' i^* 
\frac{\grad' \alpha}{2 \alpha} \cdot' + \\
&& + i^* \frac{\grad' \alpha}{2 \alpha} \cdot' i^* \frac{\grad' \alpha}{2 \alpha} \cdot'
\psi - \psi i^* \frac{\grad' \alpha}{2 \alpha} \cdot' i^* \frac{\grad' \alpha}{2 \alpha}
\cdot' + \non \\
&& + i^* \frac{\grad' \alpha}{2 \alpha} \cdot' i^* D' \psi - \psi i^*
\frac{\grad' \alpha}{2 \alpha} \cdot' i^* D' = \non \\
&& = i^* D' i^* \grad' \psi \cdot' + i^* \grad' \psi \cdot' i^* D' + \psi i^* D'
i^* D' - \psi i^* D' i^* D' + \non \\
&& + i^* ( \grad' \psi \cdot' + \psi D' ) i^* \frac{\grad' \alpha}{2 \alpha} \cdot'
- \psi i^* D' i^* \frac{\grad' \alpha}{2 \alpha} \cdot' + \\
&& + i^* \frac{\grad' \alpha}{2 \alpha} \cdot' i^* ( \grad' \psi \cdot' + \psi D' )
- \psi i^* \frac{\grad' \alpha}{2 \alpha} \cdot' i^* D' = \non \\
&& = i^* D' i^* \grad' \psi \cdot' + i^* \grad' \psi \cdot' i^* D' + \\
&& + i^* \grad' \psi \cdot' i^* \frac{\grad' \alpha}{2 \alpha} \cdot' + \\
&& + i^* \frac{\grad' \alpha}{2 \alpha} \cdot' i^* \grad' \psi \cdot'.
\eea
The terms in (8.34) are first order operators but $\grad' \psi$ has compact support
and we are done. The terms in (8.35), (8.36) are zero order operators and we are
also done since $\grad' \psi$ has compact support.

Hence 
$(e^{-tD^2} P - e^{-t \tilde{D'}^2} P')B$, 
$B(e^{-tD^2} P - e^{-t \tilde{D'}^2} P')B$
are of trace class and the trace norm in uniformly bounded on any compact $t$--interval
$[a_0,a_1]$, $a_0>0$. The assertions for (8.11) are done.

Next we study the operator
\be
(e^{\frac{t}{2}D^2} P - e^{-\frac{t}{2}{\tilde{D'}}^2} P') (1-B).
\end{equation}
Denote by $M_\epsilon$ the multiplication operator with 
$\exp ( - \epsilon \dist (m,K)^2)$.
Refering to [2], we state that for $\epsilon$ small enough
\be
M_\epsilon e^{-tD^2} B, \quad M_\epsilon e^{-t{D'}^2} B
\end{equation}
and
\be
M^{-1}_\epsilon e^{-tD^2} \raisebox{0.3ex}{$\chi$}_G, \quad M^{-1}_\epsilon e^{-t{D'}^2} \raisebox{0.3ex}{$\chi$}_G
\end{equation}
are Hilbert--Schmidt for every compact $G \subset M$ or $G' \subset M'$.
Write
\bea
&& (e^{\frac{t}{2}D^2} P - e^{-\frac{t}{2}{\tilde{D'}}^2} P') (1-B) = \non \\
&& = [e^{-\frac{t}{2}D^2} P M_\epsilon] \cdot [M^{-1}_\epsilon (e^{-\frac{t}{2}D^2}P 
- e^{-\frac{t}{2}{\tilde{D'}}^2} P') (1-B) ] + \\
&& + [ e^{-\frac{t}{2}D^2} P - e^{-\frac{t}{2}{\tilde{D'}}^2} P') M_\epsilon ] \cdot 
[M^{-1}_\epsilon e^{-\frac{t}{2}{\tilde{D'}}^2} P' (1-B)] .
\eea
According to (6.7) -- (6.9) and (8.38), each of the factors $[\cdots]$ in (8.40), (8.41)
is Hilbert--Schmidt and we obtain that (8.34) is of trace class and has uniformly bounded
trace norm in any $t$--interval $[a_0,a_1]$, $a_0>0$. The same holds for
\bea
&& B (e^{\frac{t}{2}D^2} P - e^{-\frac{t}{2}{\tilde{D'}}^2} P') (1-B) \\
&& (1-B) (e^{\frac{t}{2}D^2} P - e^{-\frac{t}{2}{\tilde{D'}}^2} P') B \\
&& (1-B) (e^{\frac{t}{2}D^2} P - e^{-\frac{t}{2}{\tilde{D'}}^2} P') (1-B)
\eea
by multiplication of (8.37) from the left by $B$ etc., i.e. the assertions for
(8.13), (8.14), (8.17) are done. Write now
\bea
&& (e^{-\frac{t}{2}D^2} D - e^{-\frac{t}{2}{\tilde{D'}}^2} D') B = \non \\
&& (e^{-\frac{t}{2}D^2} D ) (e^{\frac{t}{2}D^2} P - e^{-\frac{t}{2}{\tilde{D'}}^2} P') B + \\
&& (e^{-\frac{t}{2}D^2} D - e^{-\frac{t}{2}{\tilde{D'}}^2} D') (e^{-\frac{t}{2}{D'}^2} P) B.
\eea
(8.45) is done already by (8.23) and 
$|e^{-\frac{t}{2}D^2}D|_{op} \le \frac{C}{\sqrt{t}}$. 
Decompose (8.46) as the sum of 
\be
e^{-\frac{t}{2}D^2}P (D-\tilde{D'}) \cdot (e^{-\frac{t}{2}{\tilde{D'}}^2} \tilde{D'})
= [e^{-\frac{t}{2}D^2}P (-\eta)] \cdot (e^{-\frac{t}{2}{\tilde{D'}}^2} \tilde{D'}) B
\end{equation}
and
\be
(e^{-\frac{t}{2}D^2} P - e^{-\frac{t}{2}{\tilde{D'}}^2} P') (e^{-\frac{t}{2}{\tilde{D'}}^2}
\tilde{D'}) B
\end{equation}
$[\dots]$ in (8.47) is done. Rewrite 
$e^{-\frac{t}{2}D^2}P - e^{-\frac{t}{2}{\tilde{D'}}^2}P'$
as
\bea
&& (e^{-\frac{t}{2}D^2} P - e^{-\frac{t}{2}{\tilde{D'}}^2} P) B + \\
&& + (e^{\frac{t}{2}D^2} P - e^{-\frac{t}{2}{\tilde{D'}}^2} P') (1-B).
\eea
(8.49), (8.50) are done already, hence (8.48) and hence
$(e^{\frac{t}{2}D^2} P - e^{-\frac{t}{2}{\tilde{D'}}^2} P') B$,
(8.12), (8.15), (8.16), (8.18). 
This finishes the proof of 8.3.
\qed

The proof of theorem 8.2 now follows from 8.3 by adding up the four terms containing
$e^{\frac{t}{2}D^2} P - e^{-\frac{t}{2}{\tilde{D'}}^2} P'$ 
or
$e^{\frac{t}{2}D^2} D - e^{-\frac{t}{2}{\tilde{D'}}^2} \tilde{D'}$, respectively.
\qed

\begin{remark}
We could perform the proof of 8.2, 8.3 also along the lines of (7.59) -- (7.62),
performing first a unitrary transformation, proving the trace class porperty and
performing the back transformation, as we indicate in (7.59). This procedure is
completely equivalent to the proof of 8.2, 8.3 presented above.
\qed
\end{remark}

The operators $U^*i^*{D'}^2iU$ and $(U^*i^*{D'}^2iU)^2$ are distinct in general
and we have still to compare 
$e^{-t(U^*i^*{D'}^2iU)} P'$ and $e^{-t(U^*i^*{D'}^2iU)^2} P'$.
According to our remark above, it is sufficient to prove the trace class property of
\be
e^{-t(i^*{D'}^2i)} P - e^{-t(i^*D'i)^2} P'
\end{equation}
in 
\[ {\cal H}' = L_2 ((K,E),g,h) \oplus L_2((K',E'),g',h') \oplus 
L_2((M \setminus K,E),g',h) . \]
Here we have an embedding
\bea
&& i_{L_2,K'} \oplus i^{-1} : L_2((M',E'),g',h') \lra {\cal H}' \non \\
&& (i_{L_2,K'} \oplus i^{-1}) \Phi = i_{L_2,K'} \raisebox{0.3ex}{$\chi$}_{K'} \Phi + i^{-1} 
\raisebox{0.3ex}{$\chi$}_{M' \setminus K'} \Phi, 
\eea
where
\bea
&& i^{-1} : L_2((M' \setminus K',E'|_{M' \setminus K'}),g',h') \lra
L_2((M' \setminus K',E'|_{M' \setminus K'}),g'h), \non \\
&& i^{-1} \Phi = \Phi, \non
\eea
and
\[ i^* D' i (\raisebox{0.3ex}{$\chi$}_K, \Phi 
+ i^{-1} \raisebox{0.3ex}{$\chi$}_{M' \setminus K'} \Phi) :=
i^* D' \Phi = \raisebox{0.3ex}{$\chi$}_{K'} D' \Phi 
+ i^* \raisebox{0.3ex}{$\chi$}_{M' \setminus K'} D' \Phi , \]
$i^* {D'}^2 i$ similar, all with the canonical domains of 
definition analogous to (8.8). $P'$ is here the projection onto 
$\im (i_{L_2,K'} \oplus i^{-1})$. 
We define
$i^* {D'}^2 i$, $(i^*D'i)^2$
to be zero on ${\im P'}^\perp$.

\begin{remark}
Quite similar we could embed
$L_2((M,E),g,h)$ into ${\cal H}'$, 
define $P$, $UDU^*$ and the assertion 8.2 would be equivalent to
the assertion for
\be
e^{-t(UDU^*)^2} P - e^{-t(i^*D'i)^2} P'.
\end{equation}
Applying the (extended) $U^*$ from the right, $U$ from the left, 
yields just the expression (8.9).
\qed
\end{remark}

\begin{theorem}
Assume the hypothesises of 8.2. Then
\be
e^{-t(i^*{D'}^2i)} P' - e^{-t(i^*D'i)^2} P'
\end{equation}
is of trace class and its trace norm is uniformly bounded on compact
$t$--intervals $[a_0,a_1]$, $a_0>0$.
\end{theorem}

\begin{proof}
We prove this by establishing the assertion for the four cases arising
from multiplication by $B$ ,$1-B$. Start with (8.54). B. Duhamel's principle
again yields
\bea
&& (e^{-t(i^*{D'}^2i)} P' - e^{-t(i^*D'i)^2} P') B = \non \\
&& = -\int\limits^t_0 e^{-s(i^*{D'}^2i)} ((i^*{D'}^2i) - (i^* D' i)^2)
e^{-(t-s)(i^*D'i)^2} \,\, ds.
\eea
An easy calculation yields
\bea
&& (i^* {D'}^2 i) \psi - \psi (i^* D' i)^2 = 
i^* {D'}^2 \psi - \psi i^* D' i^* D' = \non \\
&& = i^* D' \grad' \psi \cdot' + i^* \grad' \psi \cdot' D + \psi i^* {D'}^2 - \non\\
&& - (\psi i^* {D'}^2 + \psi i^* D' (i^*-1) D') = \non \\
&& = i^* D' \grad' \psi \cdot' + i^* \grad' \psi \cdot' D' - \\
&& - \psi i^* D' (i^*-1) D' . 
\eea
The first order operators in (8.56) contain the compact support factor 
$\grad' \psi$ and we are done. Here $i^*D'$ (coming from the first term or
from 
$\grad' \psi \cdot' D' = \grad' \psi \cdot' i^{*-1} i^* D$) 
will be connected with 
$e^{-s(i^* {D'}^2 i)}$
or
$e^{-(t-s)(i^*D'i)^2}$,
depending on the interval 
$[\frac{t}{2},t]$ or $[0,\frac{t}{2}]$.
The ($D'$)'s of the second order operator (8.57) can be distributed analogous
to the proof of 7.12. The remaining main point is 
$0 \le \psi \le 1$ and $i^*-1$
Sobolev of order $r+1$, i. e. 
$i^*-1 \in \Omega^{0,1,r+1} (Hom((E'|_{M' \setminus K'},g',h'),
(E|_{M \setminus K},g',h)))$.

The assertion for (8.54) $\cdot B$ is done. Quite analogously 
(and parallel to the proofs of
(8.11), (8.13), (8.14), (8.17)) one discusses the other 3 cases.
\qed
\end{proof}

We obtain as a corollary from 8.2 and 8.6

\begin{theorem}
Assume the hypothesis of 8.2. Then for $t>0$
\[ e^{-tD^2} P - e^{-t(U^*i^*{D'}^2iU)} P' \]
is of trace class in ${\cal H}'$ and the trace norm is uniformly bounded on 
compact $t$--intervals $[a_0,a_1]$, $a_0>0$.
\qed
\end{theorem}

Applying 8.7 to the case 
$E=\Lambda^*T^*M \otimes \C$, $D^2 = \Delta$, 
we obtain

\begin{theorem}
Let $(\Lambda^*T^*M \otimes \C, g_{\Lambda^*}) \in \cl {\cal B}^{2^n} (I,B_k)$, 
$k \ge r+2 > n+3$, 
$(\Lambda^*T^*M' \otimes \C, g'_{\Lambda^*}) \in \gencomp^{1,r+1}_{L,diff,rel}
(\Lambda^*T^*M \otimes \C, g_{\Lambda^*}) \cap \cl {\cal B}^{2^n,n} (I,B_k)$, 
$\Delta = \Delta(g)$, $\Delta' = \Delta(g')$ the graded Laplace operator. Then
for $t>0$
\[ e^{-t \Delta} P - e^{-t(U^*i^*\Delta'iU} P' \]
is of trace class in ${\cal H}'$ and the trace norm is uniformly bounded on 
compact $t$--intervals $[a_0,a_1]$, $a_0>0$.
\qed
\end{theorem}

Roughly or more concretely spoken, as one prefers, this means the following.
Given an open manifold $(M^n,g)$ satisfying $(I)$ and $(B_k)$, $k \ge r+2 > n+3$.
Cut out a compact submanifold $K$ and glue the compact submanifold $K'$ along
$\partial M = \partial K$, getting thus $M'$, endow $M'$ with a metric $g'$
satisfying $(I)$ and $(B_k)$ and
\[ |g-g'|_{M \setminus K,g,1,r+1} < \infty. \]
Then for $t>0$
\[ e^{-t\Delta} P - e^{-tU^*i^*\Delta'iU} P' \]
has the asserted properties.

\section{Relative index theory}

\setcounter{equation}{0}

We proved in 8.4 that after fixing 
$E \in \cl {\cal B}^{N,n}(I,B_k)$, 
$k \ge r+1>n+3$, 
we can attach to any 
$E' \in \gencomp^{1,r+1}_{L,diff,rel}(E)$
two number valued invariants, namely
\be
E' \lra \tr (e^{-tD^2} P - e^{-t(U^*i^*D'iU)^2} P')
\end{equation}
and
\be
E' \lra \tr (e^{-tD^2} P - e^{-tU^*i^*{D'}^2iU} P') .
\end{equation}
This is a contribution to the classification inside a component but still
unsatisfactory insofar as it

1. could depend on $t$.

2. will depend on the $K \subset M$, $K' \subset M'$ in question,

3. is not yet clear the meaning of this invariant.

We are in a much nore comfortable situation if we additionally assume that 
the Clifford bundles under consideration are endowed with an involution
$\tau : E \lra E$, s.t.
\bea
&& \tau^2 = 1, \quad \tau^* = \tau \\
&& [\tau,X]_+ = 0 \en \mbox{for} \en X \in TM \\
&& [\nabla, \tau] = 0
\eea 
Then $L_2((M,E),g,h) = L_2(M,E^+) \oplus L_2(M,E^-)$
\[ D = \left( \begin{array}{cc} 0 & D^- \\ D^+ & 0 \end{array} \right) \]
and $D^- = (D^+)^*$. If $M^n$ is compact then as usual
\be
\ind D := \ind D^+ := \dimker D^+ - \dimker D^- \equiv \tr (\tau e^{-tD^2}),
\end{equation}
where we understand $\tau$ as
\[ \tau = \left( \begin{array}{cc} I & 0 \\ 0 & -I \end{array} \right) . \]
For open $M^n$ $\ind D$ in general is not defined since $\tau e^{-tD^2}$ is not
of trace class. The appropriate approach on open manifolds is relative index 
theory for pairs of operators $D,D'$. If $D,D'$ are selfadjoint in the same
Hilbert space and $e^{tD^2}-e^{-t{D'}^2}$ would be of trace class then
\be
\ind (D,D') := \tr (\tau (e^{-tD^2} - e^{-t{D'}^2}))
\end{equation}
makes sense, but at the first glance (9.7) should depend on $t$.

If we restrict to Clifford bundles 
$E \in \cl {\cal B}^{N,n}(I,B_k)$
with involution $\tau$ then we assume that the maps entering in the definition
of $\comp^{1,r+1}_{L,diff,F}(E)$ or $\gencomp^{1,r+1}_{L,diff,rel}(E)$ 
are $\tau$--compatible, i.e. after identification of $E|_{M \setminus K}$ 
and $f^*_E E'|_{M' \setminus K}$ holds
\be
[f^*_E \nabla^{h'}, \tau] = 0, \quad [f^* \cdot', \tau]_+ = 0.
\end{equation}
We call $E|_{M \setminus K}$ and $E'|_{M' \setminus K'}$ $\tau$--compatible.
Then, according to section 8, 
\be
\tr (\tau (e^{-tD^2} P - e^{-t(U^*i^*D'iU)^2} P' ))
\end{equation}
makes sense.

\begin{theorem}
Let
$((E,h,\nabla^h) \lra (M^n,g), \tau) \in \cl {\cal B}^{N,n}(I,B_k)$
be a graded Clifford bundle, $k \ge r>n+2$.

{\bf a)} If
${\nabla'}^h \in \comp^{1,r}(\nabla) \subset C^{1,r}_{E}(B_k)$, 
$\nabla'$ $\tau$--compatible, i.e. $[\nabla',\tau]=0$ then
\[ \tr (\tau (e^{-tD^2} - e^{-t{D'}^2})) \]
is independent of $t$.

{\bf b)} If
$E' \in \gencomp^{1,r+1}_{L,diff,rel}(E)$ is $\tau$--compatible with $E$, i.e.
$[\tau, X \cdot']_+=0$ for $X \in TM$ and $[\nabla',\tau]=0$, then
\[ \tr (\tau (e^{-tD^2} P - e^{-t(U^*i^*D'iU)^2} P')) \]
is independent of $t$.
\end{theorem}

\begin{proof}
a) follows from our 6.1 and 5.1 in [2]. b) follows from our 9.2.
\qed 
\end{proof}

\begin{prop}
If $E' \in \gencomp^{1,r+1}_{L,diff,rel}(E)$ and 
\bea
&& \tau (e^{-tD^2} P - e^{-t(U^*i^*D'iU)^2} P') \non \\
&& \tau (e^{-tD^2} D - e^{-t(U^*i^*D'iU)^2} (U^*i^*D'iU)) \non 
\eea
are for $t>0$ of trace class and the trace norm of 
$\tau (e^{-tD^2} D - e^{-t(U^*i^*D'iU)^2} (U^*i^*D'iU))$
is uniformly bounded on compact $t$--intervals $[a_0,a_1]$, $a_0>0$, then
\[ \tr(\tau (e^{-tD^2} p - e^{-t(U^*i^*D'iU)^2} (U^*i^*D'iU))) \]
is independent of $t$.
\end{prop}

\begin{proof}
The assumptions of 5.1 in [2] are much more restrictive then ours here.
Nevertheless, the main idea of the proof carries over to here. We sketch the
proof. Let $(\phi_i)_i$ be a sequence of smooth functions
$\in C^\infty_c (\overline{M \setminus K})$, 
satisfying
$\sup |d \phi_i| \underset{i \ra \infty}{\lra} 0$,
$0 \le \phi_i \le \phi_{i+1}$ and
$\phi_i \underset{i \ra \infty}{\lra} 1$.
Denote by $M_i$ the multiplication operator with $\phi_i$ on
$L_2((M \setminus K, E|_{M \setminus K}),g,h)$. 
We extend $M_i$ by 1 to the complement of 
$L_2((M \setminus K, E),g,h)$ in $H$.
We have to show
\[ \frac{d}{dt} \tr \tau (e^{-tD^2} P - e^{-t(U^*i^*D'iU)^2} P') = 0. \]
$e^{-tD^2} P - e^{-t(U^*i^*D'iU)^2} P'$
is of trace class, hence
\[ \tr \tau (e^{-tD^2} P - e^{-t(U^*i^*D'iU)^2} P')= 
\lim_{j \ra \infty} \tr \tau M_j (e^{-tD^2} P - e^{-t(U^*i^*D'iU)^2} P')M_j. \]
$M_j$ restricts to compact sets and we can differentiate under the trace and
we obtain 
\bea
&& \frac{d}{dt} \tr \tau M_j (e^{-tD^2} P - e^{-t(U^*i^*D'iU)^2} P') M_j = \non \\
&& = \frac{d}{dt} \tr \tau (M_j U^* (e^{-t(UDU^*)^2} P - e^{-t(i^*D'i)^2} P') U M_j = \non \\
&& =  - \tr \tau (U^* M_j (e^{-t(UDU^*)^2} (UDU^*)^2 - e^{-t(i^*D'i)^2} (i^*D'i)^2) M_j U) \non
\eea
Consider
$\tr \tau (U^* M_j (e^{-t(UDU^*)^2} (UDU^*)^2 M_j U) = \tr \tau M_j e^{-tD^2} D^2 M_j$.
According to [2], 
$\tr \tau (M_j e^{-tD^2} D^2 M_j) = \tr M_j \grad \phi_i \cdot \tau D e^{-tD^2}$.
Quite similar
\bea
&& \tr \tau (M_j (e^{-t(i^*D'i)^2} (i^*D'i)^2) M_j) = \non \\
&& = \tr \tau \phi_j e^{-\frac{t}{2}(i^*D'i)^2} (i^*Di) (i^*D'i)
e^{-\frac{t}{2}(i^*D'i)^2} \phi_j =  \non \\
&& = \tr (i^*Di) e^{-\frac{t}{2}(i^*D'i)^2} \phi_j \tau \phi_j
e^{-\frac{t}{2}(i^*D'i)^2} (i^*D'i) =  \non \\
&& = \tr e^{-\frac{t}{2}(i^*D'i)^2} (i^*D'i) \phi^2_j \tau 
e^{-\frac{t}{2}(i^*D'i)^2} (i^*D'i) =  \non \\
&& = \tr e^{-\frac{t}{2}(i^*D'i)^2} i^* (2 \phi_j \grad' \phi_j \cdot' + \phi^2_j D') i \tau 
e^{-\frac{t}{2}(i^*D'i)^2} (i^*D'i) =   \non \\
&& = \tr 2 i^* M_j \grad' \phi_j \cdot' i \tau (i^*D'i) e^{-t(i^*D'i)^2} -
\tr \tau M_j e^{-t(i^*D'i)^2} (i^*D'i)^2 M_j , \non 
\eea
hence
\[ \tr \tau (M_j e^{-t(i^*D'i)^2} (i^*D'i)^2 M_j) = \tr M_j i^* \grad' \phi_j
\cdot' i \tau (i^*D'i) e^{-t(i^*D'i)^2} \]
and finally
\bea
&& \frac{d}{dt} \tr \tau M_j (e^{-tD^2} P - e^{-t(U^*i^*D'iU)^2} P') M_j =  \non \\
&& = \tr \tau M_j [ \grad \phi_j \cdot e^{-tD^2} D - \grad' \phi_j \cdot'
e^{-t(U^*i^*D'iU)^2} (U^*i^*D'iU) ] =  \non \\
&& = \tr \tau M_j [ (\grad \phi_j - \grad' \phi_j) \cdot e^{-tD^2} + \grad' \phi_j 
(\cdot - \cdot') e^{-tD^2} + \non \\
&& + \grad' \phi_j \cdot' (e^{-tD^2}  
- e^{-t(U^*i^*D'iU)^2} (U^*i^*D'iU)) ] . \non 
\eea
But this tends to zero uniformly for $t$ in compact intervals since $\grad \phi_j$, 
$\grad' \phi_j$ do so.
\qed
\end{proof}

We denote $Q^\pm = D^\pm$
\be
Q = \left( \begin{array}{cc} 0 & Q^+ \\ Q^- & 0  \end{array} \right), \quad 
H = \left( \begin{array}{cc} H^+ & 0 \\ 0 & H^-  \end{array} \right) = 
\left( \begin{array}{cc} Q^- Q^+ & 0  \\ 0 & Q^+ Q^- \end{array} \right) = Q^2, 
\end{equation}
${Q'}^\pm = U^*i^*{D'}^\pm iU = (U^*i^*D'iU)^\pm$, $Q'$, $H'$ analogous, assuming
(9.3) -- (9.5) as before and $\cdot', \nabla'$ $\tau$--compatible. $H,H'$ form 
by definition a supersymmetric scattering system if the wave operators
\be
W^\mp (H,H') := \lim\limits_{t \ra \mp \infty} e^{itH} e^{-tH'} \cdot P_{ac} (H')
\en \mbox{exist and are complete }
\end{equation}
and
\be
QW^\mp (H,H') = W^\mp (H,H') H' \mbox{ on } {\cal D}_{H'} \cap {\cal H}'_{ac} (H').
\end{equation}
Here $P_{ac}(H')$ denotes the projection on the absolutely continuous subspace
${\cal H}'_{ac} (H') \subset {\cal H}$ of $H'$.

A well known sufficient criterion for forming a supersymmetric scattering system
is given by

\begin{prop}
Assume for the operators graded $Q,Q'$ (= supercharges)
\[ e^{-tH} - e^{-tH'} \]
and 
\[ e^{-tH} Q - e^{-tH'} Q \]
are for $t>0$ of trace class. Then they form a supersymmetric scattering system.
\qed
\end{prop}

\begin{coro}
Assume the hypothesises of 9.1.
Then $D,D'$ or $D,U^*i^*D'iU$ form a supersymmetric scattering system, respectively.
In particular, the restriction of $D,D'$ or $D, U^*i^*D'iU$ to their absolutely
continuous spectral subspaces are unitarily equivalent, respectively.
\qed
\end{coro}

Until now we have seen that under the hypothesises of 9.1
\be
\ind (D,\tilde{D'}) = \tr \tau (e^{-tD^2} P - e^{-t{\tilde{D'}}^2} P'),
\end{equation}
$\tilde{D'} = D'$ or $\tilde{D'} = U^*i^*D'iU$, is a well defined number,
independent of $t>0$ and hence yields an invariant of the pair $(E,E')$, still
depending on $K,K'$. Hence we should sometimes better write
\be
\ind (D,\tilde{D'},K,K').
\end{equation}
We want to express in some good cases $\ind (D,\tilde{D'},K,K')$ by other relevant
numbers. Consider the abstract setting (9.10). If $\inf \sigma_e (H) > 0$ then
$\ind D := \ind D^+$ is well defined.

\begin{lemma}
If $e^{-tH}p-e^{-tH'}P'$ is of trace class for all $t>0$ and
$\inf \sigma_e(H), \inf \sigma_e(H')>0$ then
\be
\lim\limits_{t \ra \infty} \tr \tau (e^{-tH}P-e^{-tH'}P') = \ind Q^+ - \ind Q^-.
\end{equation}
\end{lemma}

This is just theorem 5.2 in [2].
\qed

We infer from this

\begin{theorem}
Assume the hypothesises of 9.1 and $\inf \sigma_e (D^2) > 0$. 
Then $\inf \sigma_e({D'}^2), \inf \sigma_e(U^*i^*D'iU)^2 > 0$ 
and for each $t>0$
\be
\tr \tau (e^{-tD^2} - e^{-t{\tilde{D'}}^2}) = \ind D^+ - \ind {D'}^+.
\end{equation}
\end{theorem}

\begin{proof}
In the case 9.1.a, $\inf \sigma_e({D'}^2)>0$ follows from 
3.15 and (9.16) then follows from 9.5.
Consider the case 9.1.b. We can replace the comparison of 
$\sigma_e(D^2)$ and $\sigma_e((U^*i^*D'iU)^2)$
by that of 
$\sigma_e(UD^2U^*)$ and $\sigma_e((i^*D'i)^2)$.
Moreover, for self adjoint $A$, $0 \notin \sigma_e(A)$
if and only if $\inf \sigma_e(A^2)>0$.
Assume $0 \notin \sigma_e(UDU^*)$ and $0 \in \sigma_e(i^*D'i)$.
We must derive a contradiction. 
Let $(\Phi_\nu)_\nu$ be a Weyl sequence for
$0 \in \sigma_e(i^*D'i)$
satisfying additionally
$|\Phi_\nu|_{L_2}=1$, 
$\supp \Phi_\nu \subseteq M \setminus K = M' \setminus K'$
and for any compact
$L \subset M \setminus K = M' \setminus K'$
\be
|\Phi_\nu|_{L_2(M \setminus L} \underset{\nu \ra \infty} {\lra} 0 .
\end{equation}
We have $\lim\limits_{\nu \ra \infty} i^*D'i \Phi_\nu = 0$.
Then also 
$\lim\limits_{\nu \ra \infty} D'\Phi_\nu = 0$.
We use in the sequel the following simple fact.
If $\beta$ is an $L_2$--function, in particular if $\beta$ is even Sobolev, then
\be
|\beta \cdot \Phi_\nu|_{L_2} \lra 0 .
\end{equation}
Now
$(UDU^*) \Phi_\nu = (UDU^*-D') \Phi_\nu + D' \Phi_\nu$.
Here
$D' \Phi_\nu \underset{\nu \ra \infty} {\lra} 0$.
Consider
$(UDU^*-D') \Phi_\nu = (\alpha^{\frac{1}{2}} D \alpha^{-\frac{1}{2}} - D') \Phi_\nu
= \left( -\frac{\grad \alpha}{2\alpha} \cdot + D - D' \right) \Phi_\nu$.
Assume $\alpha \not\equiv 1$.
Then
$\beta = \left| \frac{\grad \alpha}{2\alpha} \right| \in 
\Omega^{0,2,\frac{r}{2}} (T (M \setminus K))$
satisfies the assumptions above and
\be
\lim\limits_{\nu \ra \infty} \left| - \frac{\grad \alpha}{2\alpha} \cdot \Phi_\nu 
\right|_{L_2} = 0 .
\end{equation}
If $\alpha \equiv 1$ this term does not appear. Write, according to (
7.9) -- (7.12), 
\be
(D-D') \Phi_\nu = \eta^{op}_1 \Phi_\nu + \eta^{op}_2 \Phi_\nu + \eta^{op}_3 \Phi_\nu .
\end{equation}
$\left| {g'}^{ik} \frac{\partial}{\partial x^k} \cdot \right|$
is bounded (we use a uniformly locally finite cover by normal charts, an
associated bounded decomposition of unity etc.).
$\beta = |\nabla - \nabla'|$
is Sobolev hence $L_2$ and by (9.18)
\be
|\eta^{op}_2 \Phi_\nu|_{L_2} \underset{\nu \ra \infty} {\lra} 0.
\end{equation}
Now
$|\nabla \Phi_\nu|_{L_2} \le C_1 (|\Phi_\nu|_{L_2} + |D \Phi_\nu|_{L_2} ) \le
C_2 (|\Phi_\nu|_{L_2} +|D' \Phi_\nu|_{L_2})$.
$g-g'$
is Sobolev, hence, according to (9.18) with
$\beta=|g-g'|$, 
$||g-g'| \cdot \Phi_\nu|_{L_2} \underset{\nu \ra \infty} {\lra} 0$
and finally
$||g-g'| \cdot D' \Phi_\nu|_{L_2} \underset{\nu \ra \infty} {\lra} 0$.
This yields
\be
|\eta^{op}_1 \Phi_\nu|_{L_2} \underset{\nu \ra \infty} {\lra} 0.
\end{equation}
We conclude in the same manner from $\cdot - \cdot'$ Sobolev and
$|\nabla' \Phi_\nu|_{L_2} \le C_3 (|\Phi_\nu|_{L_2} + |D' \Phi_\nu|_{L_2} )$
that
\be
|\eta^{op}_3 \Phi_\nu|_{L_2} \underset{\nu \ra \infty} {\lra} 0.
\end{equation}
(9.19) -- (9.23) yield
$(UDU^*) \Phi_\nu \lra 0$,
$0 \in \sigma_e (UDU^*)$,
$\inf \sigma_e (D^2) = 0$,
a contradiction, hence
$\inf \sigma_e ((U^*i^*D'iU)^2) > 0$,
$\inf \sigma_e ((i^*D'i)^2) > 0$,
$0 \notin \sigma_e (i^*D'i)$,
$0 \notin \sigma_e (D')$,
$\inf \sigma_e ({D'}^2) > 0$.
We infer from 9.2 and 9.5 that for $t>0$
\be
\tr \tau e^{-tD^2} - e^{-t(U^*i^*D'iU)^2} = \ind D^+ - \ind (U^*i^*D'iU)^+ .
\end{equation}
We are done if we can show
\be
\ind (U^*i^*D'iU)^+ = \ind {D'}^+ .
\end{equation}
$\Phi \in \ker (U^*i^*D'iU)^+$
means
$(U^*i^*D'iU)^+ \Phi = (U^*i^*D'iU) (\chi_{K'} \Phi + U^* i^{-1} 
\chi_{M' \setminus K'} \Phi) = \chi_{K'} {D'}^+ \Phi + U^* i^* 
\chi_{M' \setminus K'} {D'}^+ \Phi = 0$.
But this is equivalent to ${D'}^+ \Phi = 0$.
Similar for ${D'}^-$. (9.25) holds and hence (9.16).
\qed
\end{proof}

It would be desirable to express 
$\ind (D,\tilde{D'},K,K')$
by geometric topological terms. 
In particular, this would be nice in the case
$\inf \sigma_e(D^2) > 0$.
In the compact case, one sets
$\ind_a D := \ind_a D^+ = \dim \ker D^+ - \dim \ker (D^+)^* = 
\dim \ker D^+ - \dim \ker D^- = \lim\limits_{t \ra \infty} \tr \tau e^{-tD^2}$.
On the other hand, for $t \ra 0^+$ there exists the well known asymptotic
expansion for the kernel of 
$\tau e^{-tD^2}$.
Its integral at the diagonal yields the trace. If 
$\tr \tau e^{-tD^2}$
is independent of $t$ (as in the compact case), we get the index theorem where
the integrand appearing in the $L_2$--trace consists only of the $t$--free
term of the asymptotic expansion. Here one would like to express things in the
asymptotic expansion of the heat kernel of
$e^{-t{D'}^2}$ instead of $e^{-t(U^*i^*D'iU)^2}$.
For this reason we restrict in the definition of the topological index to the
case 
$E' \in \comp^{1,r+1}_{L,diff,F}(E)$
or
$E' \in \comp^{1,r+1}_{L,diff,F,rel}(E)$,
i.e. we admit Sobolev perturbation of $g,\nabla^h,\cdot$
but the fibre metric $h$ should remain fixed.
Then for
$D' = D(g',h,{\nabla'}^h,\cdot')$
in $L_2((M,E),g,h)$
the heat kernel of 
$e^{-t(U^*D'U)^2} = U^*e^{-t{D'}^2}U$
equals to 
$\alpha(q)^{-\frac{1}{2}} W'(t,q,p) \alpha(p)^{\frac{1}{2}}$.
At the diagonal this equals to $W'(t,m,m)$,
i.e. the asymptotic expansion at the diagonal of the original
$e^{-t{D'}^2}$ and the transformed to $L_2((M,E),g,h)$ coincide.

Consider
\be
\tr \tau W(t,m,m) \underset{t \ra 0^+} {\sim} t^{-\frac{n}{2}} b_{-\frac{n}{2}}
(D,m) + \cdots + b_0 (D,m) + \cdots
\end{equation}
and
\be
\tr \tau W'(t,m,m) \underset{t \ra 0^+} {\sim} t^{-\frac{n}{2}} b_{-\frac{n}{2}}
(D',m) + \cdots + b_0 (D',m) + \cdots .
\end{equation}
We show in the next section that
\be
b_i (D,m) - b_i (D',m) \in L_1 , \quad  - \frac{n}{2} \le i \le 1.
\end{equation}
Define for 
$E' \in \gencomp^{1,r+1}_{L,diff,F}(E)$
\be
\ind_{top} (D,D') := \int\limits_M b_0 (D,m) - b_0 (D',m) .
\end{equation}
According to (9.28), $\ind_{top} (D,D')$ is well defined.

\begin{theorem}
Assume $E' \in \gencomp^{1,r+1}_{L,diff,F,rel} (E)$

{\bf a)} Then
\bea
\ind (D,D',K,K') &=& \int\limits_K b_0 (D,m) - \int\limits_{K'} b_0 (D',m) + \\
&+& \int\limits_{M \setminus K = M' \setminus K'} b_0 (D,m) - b_0 (D',m). 
\eea
{\bf b)} If $E' \in \gencomp^{1,r+1}_{L,diff,F} (E)$ then
\be
\ind (D,D') = \ind_{top} (D,D').
\end{equation}
{\bf c)} If $E' \in \gencomp^{1,r+1}_{L,diff,F} (E)$ and 
$\inf \sigma_e(D^2) > 0$ then
\be
\ind_{top} (D,D') = \ind_a D - \ind_a D'.
\end{equation}
\end{theorem}

\begin{proof}
All this follows from 9.1, the asymptotic expansion, (9.28) and the fact that
the $L_2$--trace of a trace class integral operator equals to the integral
over the trace of the kernel.
\qed
\end{proof}

\begin{remarks}
{\rm 
{\bf 1.} 
If $E' \in \gencomp^{1,r+1}_{L,diff,rel} (E)$, $g$ and $g'$,
$\nabla^h$ and ${\nabla'}^h$, $\cdot$ and $\cdot'$ coincide in 
$V = M \setminus L = M' \setminus L'$, 
$L \supseteq K$, $L' \supseteq K'$, 
then in (7.4) -- (7.53) $\alpha-1$ and the $\eta$'s
have compact support and we conclude form (8.38), (8.39) and the
heat kernel estimates in section 5 that
\be
\int\limits_V |W(t,m,m) - W'(t,m,m)| \,\, dm \le C \cdot e^{-\frac{d}{t}}
\end{equation}
and obtain
\be
\ind (D,D',L,L') = \int\limits_L b_o (D,m) - \int\limits_{L'} b_0 (D',m) .
\end{equation}
This follows immediately from 9.7.a. and (9.35) and contains corollary 5.2
in [2]. It is one of the main results in [17].

{\bf 2.} 
The point here is that we admit much more general perturbations than in 
preceding approaches to prove relative index theorems.

{\bf 3.}
$\inf \sigma_e (D^2) > 0$ 
is an ivariant of 
$\gencomp^{1,r+1}_{L,diff,F} (E)$. 
If we fix $E$, $D$ as reference point in 
$\gencomp^{1,r+1}_{L,diff,F} (E)$
then 9.7.c enables us to calculate the analytical index for all other
$D$'s in the component from $\ind D$ and a pure integration. 

{\bf 4.}
$\inf \sigma_e (D^2) > 0$ 
is satisfied e.g. if in 
$D^2 = \nabla^* \nabla + {\cal R}$
the operator ${\cal R}$
satisfies outside a compact $K$ the condition
\be
{\cal R} \ge \kappa_0 \cdot id, \kappa_0 > 0.
\end{equation}
(9.36) is an invariant of 
$\gencomp^{1,r+1}_{L,diff,F} (E)$
(with possibly different $K$, $\kappa_0$).
}

\qed
\end{remarks}

It is possible that $\ind D$, $\ind D'$ are defined even if $0 \in \sigma_e$.
For the corresponding relative index theorem we need the scattering index.

To define the scattering index and in the next section relative $\zeta$--functions, 
we must introduce the spectral shift function of Birman/Krein/Yafaev.
Let $A$, $A'$ be bounded self adjoint operators, $V=A-A'$ of trace class, 
$R'(z) = (A'-z)^{-1}$.
Then the spectral shift function
\be
\xi(\lambda) = \xi(\lambda,A,A') := \pi^{-1} \lim\limits_{\epsilon \ra 0}
\arg \det (1+VR'(\lambda+i\epsilon))
\end{equation}
exists for a.e. $\lambda \in \R$.
$\xi(\lambda)$ is real valued, $\in L_1(\R)$ and
\be 
\tr (A-A') = \int\limits_\R \xi (\lambda) \,\, d \lambda , \quad
|\xi|_{L_1} \le |A-A'|_1 . 
\end{equation}
If $I(A,A')$ is the smallest interval containing
$\sigma(A) \cup \sigma(A')$
then $\xi(\lambda)=0$ for $\lambda \notin I(A,A')$.

Let 
\[ {\cal G} = \{ f: \R \lra \R \en | \en f \in L_1 \quad \mbox{and} \quad 
\int\limits_\R |\widehat{f}(p)| (1+|p|) \,\, dp < \infty \} . \]
Then for $\phi \in {\cal G}$, $\phi(A)-\phi(A')$
is of trace class and
\be
\tr (\phi(A)-\phi(A')) = \int\limits_\R \phi'(\lambda) \xi(\lambda) \,\, d \lambda .
\end{equation}
We recall proposition 2.1 from [19].

\begin{lemma}
Let $H, H' \ge 0$, selfadjoint in ${\cal H}$, $e^{-tH}-e^{-tH'}$
for $t>0$ of trace class. Then there exist a unique function
$\xi = \xi(\lambda) = \xi(\lambda,H,H') \in L_{1,loc} (\R)$
such that for $>0$, 
$e^{-t\lambda} \xi (\lambda) \in L_1(\R)$
and the following holds.

{\bf a)}
$\tr (e^{-tH}-e^{-tH'}) = -t \int\limits^\infty_0 e^{-t\lambda} \xi(\lambda) \,\, d\lambda$.

{\bf b)} 
For every $\phi \in {\cal G}$, $\phi(H)-\phi(H')$ is of trace class and
\[ \tr (\phi(H)-\phi(H')) = \int\limits_\R \phi'(\lambda) \xi(\lambda) \,\, d\lambda . \]
{\bf c)}
$\xi(\lambda) = 0$ for $\lambda<0$.
\qed
\end{lemma}

We apply this to our case
$E' \in \gencomp^{1,r+1}_{L,diff,rel} (E)$.
According to 9.4, $D$ and $U^*i^*D'iU$ form a supersymmetric scattering system,
$H=D^2$, $H'=(U^*i^*D'iU)^2$.
In this case
\[ e^{2\pi i \xi(\lambda,H,H')} = \det S(\lambda) , \]
where $S = (W^+)^* W^- = \int S(\lambda) \,\, dE'(\lambda)$ 
and $H'_{ac} = \int \lambda \,\, d E'(\lambda)$.

Let $P_d(D)$, $P_d(U^*i^*D'iU)$
be the projector on the discrete subspace in ${\cal H}$, respectively and
$P_c=1-P_d$ the projector onto the continuous subspace. Moreover we write
\be
D^2 = \left( \begin{array}{cc} H^+ & 0 \\ 0 & H^- \end{array} \right), \quad
(U^*i^*D'iU)^2 = \left( \begin{array}{cc} {H'}^+ & 0 \\ 0 & {H'}^- \end{array} \right) .
\end{equation}
We make the following additional assumption.
\be
e^{-tD^2} P_d(D), e^{-t(U^*i^*D'iU)^2} P_d(U^*i^*D'iU) \en
\mbox{ are for } t>0 \mbox{ of trace class.}
\end{equation}
Then for $t>0$
\[ e^{-tD^2} P_c(D) - e^{-t(U^*i^*D'iU)^2} P_c(U^*i^*D'iU) \]
is of trace class and we can in complete analogy to (9.37) define
\bea
\xi^c (\lambda,H^\pm,{H'}^\pm) & := & -\pi \lim\limits_{\epsilon \ra 0^+} \arg \det
[ 1+(e^{-tH^\pm} P_c(H^\pm) - e^{-t{H'}^\pm} P_c({H'}^\pm)) \non \\
&& ( e^{-t{H'}^\pm} P_c({H'}^\pm) - e^{-\lambda t} - i \epsilon )^{-1} ]
\eea
According to (9.38), 
\be
\tr (e^{-tH^\pm} P_c (H^\pm) - e^{-t{H'}^\pm} P_c ({H'}^\pm)) 
= -t \int\limits^\infty_0 \xi^c ( \lambda, H^\pm, {H'}^\pm ) 
e^{-t \lambda} \,\, d \lambda.
\end{equation}
We denote as after (9.13) $\tilde{D'} = D'$ in the case
$\nabla' \in \comp^{1,r} (\nabla)$ and $\tilde{D'} = U^*i^*D'iU$
in the case $E' \in \gencomp^{1,r+1}_{L,diff,rel} (E)$.
The assumption (9.41) in particular implies that for the restriction
of $D$ and $\tilde{D'}$ to their discrete subspace the analytical
index is well defined and we write
$\ind_{a,d} (D,\tilde{D'}) = \ind_{a,d} (D) - \ind_{a,d} (\tilde{D'})$
for it. Set
\be
n^c (\lambda, D, \tilde{D'}) := - \xi^c (\lambda, H^+, {H'}^+) + 
\xi^c (\lambda,H^-,{H'}^-).
\end{equation}

\begin{theorem}
Assume the hypothesises of 9.1 and (9.41). \\
Then 
$n^c (\lambda, D, \tilde{D'}) = n^c (D, \tilde{D'})$
is constant and
\be
\ind (D,\tilde{D'}) - \ind_{a,d} (D,\tilde{D'}) = n^c (D,\tilde{D'}).
\end{equation}
\end{theorem}

\begin{proof}
\bea
\ind(D,\tilde{D'}) &=& \tr \tau (e^{-tD^2} P - e^{-t{\tilde{D'}}^2} P') = \non \\
&=& \tr \tau e^{-tD^2} P_d (D) P - \tr \tau e^{-t{\tilde{D'}}^2} P_d (\tilde{D'}) P' + \non \\
&& + \tr \tau (e^{-t D^2} P_c (D) - e^{-t{\tilde{D'}}^2} P_c (\tilde{D'})) = \non \\
&=& \ind_{a,d} (D,\tilde{D'}) + t \int\limits^\infty_0 e^{-t \lambda} n^c
(\lambda, D, \tilde{D'}) \,\, d \lambda. \non
\eea
According to 9.1, 
$\ind(D,\tilde{D'})$ 
is independent of $t$. The same holds for
$\ind_{a,d} (D,\tilde{D'})$. 
Hence
$t \int\limits^\infty_0 e^{-t \lambda} n^c (\lambda, D, \tilde{D'}) \,\, d \lambda$
is independent of $t$. This is possible only if
$\int\limits^\infty_0 e^{-t \lambda} n^c (\lambda, D, D') \,\, d \lambda = \frac{1}{t}$
or $n^c(\lambda, D, \tilde{D'})$ is independent of $\lambda$.
\qed
\end{proof}

\begin{coro}
Assume the hypothesises of 9.10 and additionally 
$\inf \sigma_e (D^2|_{(\ker D^2)^\perp}) > 0$.
Then $n^c (D, \tilde{D'}) = 0$.
\end{coro}

\begin{proof}
In this case 
$\ind_{a,d} (D, \tilde{D'}) = \ind D - \ind \tilde{D'} = \ind(D, \tilde{D'})$,
hence $n^c=0$.
\qed
\end{proof}

\section{Relative $\zeta$--functions, $\eta$--functions, determinants and torsion}

\setcounter{equation}{0}

Assume
$E' \in \gencomp^{1,r+1}_{L,diff,F} (E)$.
Then we have in $L_2 ((E,M),g,h)$ the asymptotic expansion
\be
  \tr \, W(t,m,m) 
  \underset{t \rightarrow 0^+} {\sim}    
  t^{-\frac{n}{2}} b_{-\frac{n}{2}} (m) + 
  t^{-\frac{n}{2}+1} b_{-\frac{n}{2}+1} + \dots
\end{equation}
and analogously for 
$\tr \alpha^{-\frac{1}{2}}(m) W'(t,m,m) \alpha^{\frac{1}{2}}(m) = \tr \, W'(t,m,m)$ 
with 
\[ b_{-\frac{n}{2}+1}(m) =  b_{-\frac{n}{2}+l} (D(g,h,\nabla),m), \qu 
   b'_{-\frac{n}{2}+l} (m) =  b_{-\frac{n}{2}+l} (D(g',h,\nabla'),m) . \]
Here we use that the odd coefficients vanish, i.e. terms with
$t^{-\frac{n}{2} + \frac{1}{2}}$, $t^{-\frac{n}{2} + \frac{3}{2}}$
etc. do not appear.
The heat kernel coefficients have for $l \ge 1$ a representation
\be
  b_{-\frac{n}{2}+l} = \sum^l_{k=1} \sum^k_{q=0} 
  \sum_{i_1 + i_2 + \cdots + i_k = 2 (l-k)} \nabla^{i_1} R^g \dots \nabla^{i_q} R^g
  \tr \, (\nabla^{i_{q+1}} R^E \dots \nabla^{i_k} R^E ) C^{i_1, \dots , i_k} ,
\end{equation}                                          
where $C^{i_1, \dots , i_k}$ stands for a contraction with respect to $g$, 
i.e. it is built up by linear combination of products of the 
$g^{ij}$, $g_{ij}$.

\begin{lemma} 
$b_{-\frac{n}{2}+l} - b'_{-\frac{n}{2}+l} \in L_1(M,g)$, 
$0 \le l \le \frac{n+3}{2}$. 
\end{lemma}

\begin{proof} 
First we fix $g$. Forming the difference 
$b_{-\frac{n}{2}+l} - b'_{-\frac{n}{2}+l}$, we obtain a sum of terms of
the kind 
\be
   \nabla^{i_1} R^g \dots \nabla^{i_q} R^g \, \tr \, 
   [\nabla^{i_{q+1}} R^E \dots \nabla^{i_k} R^E - 
   {\nabla'}^{i_{q+1}} {R'}^E \dots {\nabla'}^{i_k} {R'}^E]
   C^{i_1,\dots,i_k} .
\end{equation}
The highest derivative of $R^q$ with respect to $\nabla^g$
occurs if $q=k, i_1= \dots = i_{q-1} = 0$. Then we have
\be
   (\nabla^g)^{2l-2k} R^g . 
\end{equation}
By assumption, we have bounded geometry of order $\ge r > n+2$,
i. e. of order $\ge n+3$. Hence $(\nabla^g)^i R^g$ is bounded for
$i \le n+1$. To obtain bounded $\nabla^j R^g$--coefficients of $[\dots]$ in 
(10.3), we must assume
\be
  2l-2 \le n+1, \qu l \le \frac{n+3}{2} .
\end{equation}
Similarly we see that the highest occuring derivatives of $R^E$, 
${R'}^E$ in $[\dots]$ are of order $2l-2$. The corresponding expression is
\bea
   R^E \nabla^{2l-2} R^E - {R'}^E {\nabla'}^{2l-2} {R'}^E & = &
   (R^E-{R'}^E)(\nabla^{2l-2} R^E) \non \\
   && + {R'}^E (\nabla^{2l-2} R^E - \nabla^{2l-2} {R'}^E) . 
\eea
We want to apply the module structure theorem.
$\nabla - \nabla' \in \Omega^{1,1,r} ({\cal G}^{Cl}_E, \nabla) =
 \Omega^{1,1,r} ({\cal G}^{Cl}_E, \nabla')$ 
implies $R^E-{R'}^E \in \Omega^{2,1,r-1}$.
We can apply the module structure theorem (and conclude that all norm
products of derivatives of order $\le 2l-2$ are absolutely
integrable) if $2l-2 \le r-1$, $2l-2 \le n+1$, $l \le \frac{n+3}{2}$.
Hence, (10.5) $\in L_1$ since $R^E$, ${R'}^E$, $C^{i_1 \dots i_k}$ bounded. 
It is now a very
simple combinatorial matter to write $[\dots]$ in (10.3) as a sum of terms
each of them is a product of differences $(\nabla^i R^E - {\nabla'}^i {R'}^E)$
with bounded functions $\nabla^j R^E$, ${\nabla'}^{j'} {R'}^E$. Remember
$\nabla, \nabla' \in {\cal C}_E(B_k)$. 
Admit now change of $g$.
We write
$\nabla^{i_1} R^g \cdots \nabla^{i_q} R^g \tr (\nabla^{i_{q+1}} R^E \cdots
\nabla^{i_k} R^E) C^{i_1,\dots,i_k}$ 
as
${\cal R}^1(g) \tr {\cal R}^2 (h,\nabla) C(g)$,
similarly
${\cal R}^1(g') \tr {\cal R}^2 (h,\nabla') C(g')$.
Then we have to consider expressions
\bea
&& {\cal R}^1(g) \tr {\cal R}^2 (h,\nabla) C(g)
- {\cal R}^1(g') \tr' {\cal R}^2 (h,\nabla') C(g') = \non \\
&& = [ {\cal R}^1(g) - {\cal R}^1(g') ] \tr {\cal R}^2 (h,\nabla) C(g) + \non \\
&& + {\cal R}^1(g') \tr {\cal R}^2 (h,\nabla) [ C(g) - C(g') ] \non \\
&& + {\cal R}^1(g) \tr [ {\cal R}^2 (h,\nabla) - {\cal R}^2 (h,\nabla') ] C(g'). \non
\eea
But each term $[ \dots ] \in L_1 (g)$ and the others are bounded
what we infer as above. This proves 10.1.
\qed
\end{proof}

\begin{lemma} 
There is an expansion
\be
  \tr \, (e^{-tD^2}-e^{-t(U^*D'U)^2})= t^{-\frac{n}{2}} a_{-\frac{n}{2}}+
  \dots + t^{-\frac{n}{2}+[\frac{n+3}{2}]} a_{-\frac{n}{2}+[\frac{n+3}{2}]}
  + O(t^{-\frac{n}{2}+[\frac{n+3}{2}]+1}) .
\end{equation}
\end{lemma}

\begin{proof} 
Set
\be
  a_{-\frac{n}{2}+i} = \int \Big( b_{-\frac{n}{2}+i}(m) -
  b'_{-\frac{n}{2}+i}(m)\Big) \,\, dm
\end{equation}
and use
\bea
  \tr \, W(t,m,m) &=&  t^{-\frac{n}{2}} b_{-\frac{n}{2}}+
  \dots + t^{-\frac{n}{2}+[\frac{n+3}{2}]} b_{-\frac{n}{2}+[\frac{n+3}{2}]} + \non \\
  && + O(m, t^{-\frac{n}{2}+[\frac{n+3}{2}]+1}) , \\           
   \tr \, W'(t,m,m) &=&  t^{-\frac{n}{2}} b'_{-\frac{n}{2}} + \dots 
  + O'(m, t^{-\frac{n}{2}+[\frac{n+3}{2}]+1}) \\
  \tr \, (e^{-tD^2}-e^{-t(U^*D'U)^2}) &=& \int \Big( 
  \tr \, W(t,m,m) - \tr \, W'(t,m,m) \Big) \,\, dm . \nonumber
\eea                 
Using lemma 10.1, the only critical point is 
\be
  \int\limits_M O(m, t^{-\frac{n}{2}+[\frac{n+3}{2}]+1}) -
  O'(m, t^{-\frac{n}{2}+[\frac{n+3}{2}]+1}) \,\, dm = 
  O(t^{-\frac{n}{2}+[\frac{n+3}{2}]+1}) .  
\end{equation}                      
(10.11) requires a very careful investigation of the concrete
representatives for $O(m, t^{-\frac{n}{2}+[\frac{n+3}{2}]})$.
We did this step by step, following [16], p. 21/22, 66 -- 69. Very
roughly spoken, the $m$--dependence of $O(m,\cdot)$ is given by 
the parametwise construction, i. e. by differences
of corresponding derivatives of the $\Gamma^\beta_{i \alpha}$,
${\Gamma'}^\beta_{i \alpha}$, which are integrable by assumption.
\qed
\end{proof}

\begin{defi} 
Assume $E' \in \gencomp^{1,r+1}_{L,diff,F} (E)$. 
Set
\be
  \zeta_1 (s,D^2,(U^*D'U)^2) := \frac{1}{\Gamma(s)} \int\limits^1_0 t^{s-1}
  \tr \, (e^{-tD^2}-e^{-t(U^*D'U)^2}) \,\, dt.
\end{equation}                     
\end{defi}

Using 10.7, 
\be
  \int\limits^1_0 t^{s-1} t^{-\frac{n}{2}+[\frac{n+3}{2}]} \,\, dt =
  \frac{1}{s+\frac{n}{2}+[\frac{n+3}{2}]} ,
\end{equation}
\be
  \frac{1}{\Gamma(s)} \int\limits^1_0 t^{s-1} 
  O(t^{-\frac{n}{2}+[\frac{n+3}{2}]+1}) \,\, dt \en \mbox{holomorphic for} \en
  \re(s)+(-\frac{n}{2})+[\frac{n+3}{2}]+1>0
\end{equation}
and $[\frac{n+3}{2}] \ge \frac{n}{2}+1$, we
obtain a function meromorphic in $\re(s) > -1$, holomorphic in 
$s=0$ with simple poles at $s=\frac{n}{2}-l$, $l \le [\frac{n+3}{2}]$.

Much more troubles causes the integral $\int\limits^\infty_1$.
Here we must additionally assume
\be
\mu_e(D) =  \inf \sigma_e (D^2|_{(\ker D^2)^\perp}) > 0 . 
\end{equation}
(10.15) implies 
$\mu_e ({D'}^2) = \inf \sigma_e ((U^*D'U)^2 |_{(\ker (U^*D'U)^2)^\perp}) > 0$.
Denote by $\mu_0 (D^2)$, $\mu_0 ({D'}^2) = \mu_0 ((U^*D'U)^2)$ 
the smallest positive eigenvalue of $D^2$, ${D'}^2$, respectively and set
\bea
\mu(D^2) & = & \min \{ \mu_e(D^2), \mu_0(D^2) \}, \non \\
\mu({D'}^2) & = & \min \{ \mu_e({D'}^2), \mu_0({D'}^2) \}, \non \\
\mu(D^2, {D'}^2) & := & \min \{ \mu(D^2), \mu({D'}^2) \} >0 .
\eea
If there is no such eigenvalue for $D^2$ then set 
$\mu(D^2) = \mu_e(D^2)$,
analogous for ${D'}^2$.
$D^2$, ${D'}^2$, $(U^*D'U)^2$
have in $] 0, \mu(D^2, {D'}^2) [$
no further spectral values.

We assert that the spectral function 
$\xi (\lambda) = \xi (\lambda, D^2, (U^*D'U)^2)$
is constant in the interval
$[ 0, \mu(D^2, {D'}^2) / 2 [$.

Consider the function
$\omega_\epsilon (x) = \left\{ \begin{array}{cl} c_\epsilon 
e^{- \frac{\epsilon^2}{\epsilon^2-x^2}} & |x| \le \epsilon \\ 
0 & |x| > \epsilon \end{array} \right. $
and choose $c_\epsilon$ s. t.
$\int \omega_\epsilon (x) \,\, dx = 1$.
Let 
$0 < 3 \epsilon < \frac{\delta}{2}$
and 
$\raisebox{0.3ex}{$\chi$}_{[-\delta-2\epsilon, \delta+2\epsilon]}$
the characteristic function of 
$[-\delta-2\epsilon, \delta+2\epsilon]$.
Then 
$\phi_{\delta,\epsilon} := 
\raisebox{0.3ex}{$\chi$}_{[-\delta-2\epsilon, \delta+2\epsilon]} * \omega_\epsilon$
satifies
$0 \le \phi_{\delta,\epsilon} \le 1$,
$\phi_{\delta,\epsilon} (x) = 1$
on 
$[-\delta-\epsilon, \delta+\epsilon]$,
$\phi_{\delta,\epsilon} (x) = 0$
for 
$x \notin ]-\delta-3\epsilon, \delta+3\epsilon[$,
$\phi'_{\delta,\epsilon} \le K \cdot \epsilon^{-1}$
and 
$\lim\limits_{\epsilon \ra 0} \lim\limits_{\delta \ra 0} \phi_{\delta,\epsilon}
= \delta -$ distribution. Assume
$\delta + 3 \epsilon < \frac{\mu}{2}$.
A regular distribution
$f \in D' (] - \frac{\mu}{2}, \frac{\mu}{2} [)$
equals to zero if and only if
$\sk{f(\lambda)}{w_\epsilon (\lambda-a) \sin (k(\lambda-a))} = 0$, 
$\sk{f(\lambda)}{w_\epsilon (\lambda-a) \cos (k(\lambda-a))} = 0$
for all sufficiently small $\epsilon$ and all $a$ (s.t. 
$|a| + \epsilon < \frac{\mu}{2}$ ) and for all $k$.
([22], p. 95). This is equivalent to 
$\sk{f}{w_\epsilon (\lambda-a)} = 0$
for all sufficiently small $\epsilon$ and $a$ and the latter is equivalent to
$\sk{f(\lambda)}{\phi_{\delta,\epsilon} - \phi_{\delta',\epsilon}} = 0$
for all $\delta, \delta'$ and all $\epsilon$
(s.t. $\delta + 3 \epsilon, \delta' + 3 \epsilon < \frac{\mu}{2}$).
We get for 
$\mu = \mu (D^2, {D'}^2)$,
$0 < 3 \epsilon < \frac{\delta}{2}$,
$\delta + 3 \epsilon < \frac{\mu}{2}$
that
$\phi_{\delta,\epsilon} (D^2) - \phi_{\delta,\epsilon} ((U^*D'U)^2)$
is independent of $\delta, \epsilon$,
$\tr (\phi_{\delta,\epsilon} (D^2) - \phi_{\delta,\epsilon} ((U^*D'U)^2))$
is independent of $\delta, \epsilon$,
$0 = \tr (\phi_{\delta,\epsilon} (D^2) - \phi_{\delta,\epsilon} ((U^*D'U)^2))
- \tr (\phi_{\delta',\epsilon} (D^2) - \phi_{\delta',\epsilon} ((U^*D'U)^2))
= \int\limits^{\frac{\mu}{2}}_0 (\phi_{\delta,\epsilon} - \phi_{\delta',\epsilon})'
(\lambda) \xi(\lambda) \,\, d \lambda$, 
i.e. the distributional derivative of $\xi$ equals to zero, 
$\xi(\lambda)$ is a constant regular distribution. We write
$\xi(\lambda) |_{[ 0, \frac{\mu}{2} [} = 
\tr (\phi_{\delta,\epsilon} (D^2) - \phi_{\delta,\epsilon} ((U^*D'U)^2)) \equiv - h$.
Set quite parallel to [19]
$\tilde{\xi}(\lambda) := \xi(\lambda) + h$
which yields
$\tilde{\xi}(\lambda) = 0$
for 
$\lambda < \frac{\mu}{2}$
and
$-t \int\limits^\infty_0 e^{-t\lambda} \xi (\lambda) \,\, d \lambda
= h - \int\limits^\infty_{\frac{\mu}{2}} e^{-t\lambda} \tilde{\xi} (\lambda) \,\, d \lambda$.
The latter integral converges for $t>0$ and can for $t \ge 1$ estimated by
\[ e^{-t \frac{\mu}{4}} \int\limits^\infty_{\frac{\mu}{2}} |\tilde{\xi} (\lambda)|
   e^{-t \frac{\lambda}{4}} \,\, d \lambda \le C e^{-t \frac{\mu}{4}}. \]

Hence we proved
\begin{prop}
Assume $E' \in \comp^{1,r+1}_{L,diff,F} (E)$,
$\inf \sigma_e (D^2|_{(\ker D^2)^\perp}) > 0$
and set 
$h = \tr (\phi_{\delta,\epsilon} (D^2) - \phi_{\delta,\epsilon} ((U^*D'U)^2))$
as above. Then there exist $c>0$ s.t.
\be
\tr (e^{-tD^2} - e^{-t(U^*D'U)^2}) = h + O(e^{-ct}) .
\end{equation}
\qed
\end{prop}

Define for $\re s <0$
\be
\zeta_2 (s,D^2,{D'}^2) := \frac{1}{\Gamma(s)} \int\limits^\infty_1
t^{s-1} [ \tr (e^{-tD^2} - e^{-t(U^*D'U)^2}) - h] \,\, ds.
\end{equation}
Then $\zeta_2(s,D^2,{D'}^2)$ is holomorphic in $\re s <0$ and admits a
meromorphic extension to $\C$ which is holomorphic in $s=0$.

Define finally
\bea
\zeta (s,D^2,{D'}^2) & := & \frac{1}{\Gamma(s)} \int\limits^\infty_0
t^{s-1} [ \tr (e^{-tD^2} - e^{-t(U^*D'U)^2}) - h ] \,\, dt \non \\
& = & \zeta_1 (s,D^2,{D'}^2) + \zeta_2 (s,D^2,{D'}^2) - \frac{1}{\Gamma (s)}
\int\limits^1_0 t^{s-1} h \,\, dt \non \\
& = & \zeta_1 (s,D^2,{D'}^2) + \zeta_2 (s,D^2,{D'}^2) - \frac{h}{\Gamma (s+1)}
\eea

We proved

\begin{theorem}
Suppose $E' \in \gencomp^{1,r+1}_{L,diff,F} (E)$, 
$\inf \sigma_e (D^2 |_{(\ker D^2)^\perp}) > 0$
and set $h$ as above. Then $\zeta (s,D^2,{D'}^2)$
is after meromorphic extension well defined in $\re s > -1$ and
holomorphic in $s=0$.
\qed
\end{theorem}

We know from QFT that functional integrals of the type
\be
\int\limits_A e^{-\sk{H \phi}{\phi}} \,\, d \phi
\end{equation}
play a decisive role. It is very difficult to give a reasonable sense
to (10.20). An oftenly used kind to do this is to define -- in analogy
to the Gau{\ss} integral -- (10.20) by the regularized determinant
$\det H = e^{-\zeta'(0,H)}$, where $\zeta(s,H)$ is the zeta function
of $H$. This makes sense if $\zeta(s,H)$ is defined and holomorphic
at $s=0$. Exactly spoken, this is for open manifolds very rarely the case and definitely
wrong for underlying manifolds satisfying $(I)$ and $(B_k)$.
We are now able to rescue this situation by considering relative 
determinants, 
\be
\det (H,H') := e^{-\zeta'(0,H,H')} .
\end{equation}
If 
$E', E'' \in \gencomp^{1,r+1}_{L,diff,F} (E)$ 
then we denote as above
$\tilde{D'}^2 = (U^*D'U)^2$, 
$\tilde{D''}^2 = (V^*D''V)^2$
for the transformed operators acting in 
$L_2((M,E),g,h)$.

\begin{theorem}
Suppose $E', E'' \in \gencomp^{1,r+1}_{L,diff,F} (E)$ and
$\inf \sigma (D^2 |_{(\ker D^2)^\perp})$ $>0$.

{\bf a)} Then
$\zeta(s,D^2,\tilde{D'}^2)$, 
$\zeta(s,D^2,\tilde{D''}^2)$,
$\zeta(s,\tilde{D'}^2,\tilde{D''}^2)$
are after meromorphic extension in $\re s > -1$ well defined and
holomorphic in $s=0$.
In particular
$\det (D^2,{D'}^2) = e^{-\zeta'(0,D^2,{\tilde{D'}}^2)}$, 
$\det (D^2,{D''}^2) = e^{-\zeta'(0,D^2,{\tilde{D''}}^2)}$, 
$\det (\tilde{D'}^2,{D''}^2) = e^{-\zeta'(0,{\tilde{D'}}^2,{\tilde{D''}}^2)}$
are well defined.

{\bf b)} There holds
\be
\det (\tilde{D'}^2, D^2) = \det (D^2, \tilde{D'}^2)^{-1}
\end{equation}
etc. and
\be
\det (D^2, \tilde{D''}^2) = \det (D^2, \tilde{D'}^2) \cdot 
\det (\tilde{D'}^2, \tilde{D''}^2).
\end{equation}
\end{theorem}

\begin{proof}
a) follows from 10.4 and the fact that
$E, E'' \in \gencomp (E)$ 
implies
$E'' \in \gencomp (E') (= \gencomp(E))$.

b) immediately follows from the definitions and
$\tr (e^{-tD^2}-e^{-t{\tilde{D''}}^2}) = \tr (e^{-tD^2}-e^{-t{\tilde{D'}}^2})
+ \tr (e^{-t{\tilde{D'}}^2}-e^{-t{\tilde{D''}}^2})$.
\qed
\end{proof}

If we now restrict to the case 
$E = (\Lambda^* T^* M \otimes \C, g_\Lambda)$
then, as we have seen already in section 7, 
$g' \in \gencomp^{1,r+1}(g)$
does not imply
$E' = (\Lambda^* T^* M \otimes \C, g'_\Lambda) \in \gencomp^{1,r+1}_{L,diff,F}(E)$
since the fibre metric changes, 
$g_\Lambda \lra g'_\Lambda$.
Hence the above considerations for constructing the relative $\zeta$--function 
are not immediately applicable. Fortunately we can define relative $\zeta$--functions
also in this case. We recall from [16], p. 65 -- 74 the following well known fact
which we used in (10.1), (10.2) already. Let $P$ be a self adjoint elliptic partial
differential operator of order 2 such that the leading symbol of $P$ is positive
definite, acting on sections of a vector bundle
$(V,h) \lra (M^n,g)$.
Let
$W_P(t,p,m)$
be the heat kernel of
$e^{-tP}$, $t>0$.
Then for
$t \lra 0^+$
\[ \tr W_P(t,m,m) \sim t^{-\frac{n}{2}} b_{-\frac{n}{2}} (m) +
   t^{-\frac{n}{2}+1} b_{-\frac{n}{2}+1} (m) + \cdots \]
and the $b_\nu (m)$ can be locally calculated as certain derivatives of the 
symbol of $P$ according to fixed rules. As established by Gilkey, for 
$P=\Delta$ or $P=D^2$ the $b$'s 
can be expressed by curvature expressions (including derivatives). 
This is (10.1), (10.2), (10.3). We apply this to
$e^{-t\Delta}$ and $e^{-t(U^*i^*\Delta'iU)}$ 
but we want to compare the asymptotic expansions of 
$W_\Delta (t,m,m)$ and $W_{\Delta'}(t,m,m)$.
The expansion of
\be
W_{i^*\Delta'i} (t,m,m) \mbox{ in } L_2(g',g_{\Lambda^*})
\end{equation}
and
\be
W_{U^*i^*\Delta'iU} (t,m,m) \mbox{ in } L_2(g,g_{\Lambda^*})
\end{equation}
coincide since
\be
W_{U^*i^*\Delta'iU} (t,m,m) = \alpha^{-\frac{1}{2}}(m) 
W_{i^*\Delta'i} (t,m,m) \alpha^{\frac{1}{2}} (m).
\end{equation}
The point is to compare the expansions of 
\be
W_{i^*\Delta'i} (t,m,m) \mbox{ in } L_2(g',g_{\Lambda^*})
\end{equation}
and 
\be
W_{\Delta'} (t,m,m) \mbox{ in } L_2(g',g'_{\Lambda^*}), 
\end{equation}
i.e. we have to compare the symbol of 
$i^*\Delta'i = i^*\Delta'$ and $\Delta'$.
For $q=0$ they coincide.
Let $q=1$, $m \in M$, $\omega_1, \dots, \omega_n$
a basis in $T^*_m M$, 
$\Phi \in \Omega^1(M)$, 
$\Delta' \Phi|_m = \xi^1 \omega_1 + \cdots + \xi^n \omega_n$.
Then, according to (7.55)
$i^* (\Delta' \Phi_i) |_m = g^{kl} g'_{ik} \xi^i \omega_l$,
i.e.
\be
((i^*-1) \Delta'_1 \Phi) |_m = (g^{kl} g'_{ik} \xi^i - \xi^l) \omega_l .
\end{equation}
Hence for the (local) coefficients of 
$i^*\Delta'i$
as differential operator holds
\be
\mbox{coeff of } (i^*\Delta'_1i) = (g^{kl} g'_{ik}) 
\mbox{ coeff of } (\Delta'_1) .
\end{equation}
Quite similar for 
$0<q<n$, e.g.
\[ \mbox{coeff of } (i^*\Delta'_2i) = (g^{k_1l_1} g^{k_2l_2} 
g'_{i_1k_1} g'_{i_2k_2}) \mbox{ coeff of } (\Delta'_2) . \]

\begin{prop}
Let $r>n+2$, $g' \in \comp^{1,r+1}(g)$, $l \le \frac{n+3}{2}$,
$b_{-\frac{n}{2}+l} (\Delta(g,g_{\Lambda^*}),g,g_{\Lambda^*},m)$ and
$b_{-\frac{n}{2}+l} (U^*i^*\Delta'(g',g'_{\Lambda^*})iU,g,g_{\Lambda^*})$
the coefficients of the asymptotic expansion of
$\tr_{g_{\Lambda^*}} W_\Delta (t,m,m)$ and
$\tr_{g_{\Lambda^*}} W_{U^*i^*\Delta'iU} (t,m,m)$
in $L_2(M,g)$, repectively. Then
\be
b_{-\frac{n}{2}+l} (\Delta,g,g_{\Lambda^*},m) - b_{-\frac{n}{2}+l}
(U^*i^*\Delta'iU,g,g_{\Lambda^*},m) \in L_1 (M,g).
\end{equation}
\end{prop}

\begin{proof}
Write
\bea
&& b_{-\frac{n}{2}+l} (\Delta,g,g_{\Lambda^*},m) - b_{-\frac{n}{2}+l}
(U^*i^*\Delta'iU,g,g_{\Lambda^*},m) \non \\
&& = b_{-\frac{n}{2}+l} (\Delta,g,g_{\Lambda^*},m) - b_{-\frac{n}{2}+l}
(\Delta',g',g'_{\Lambda^*},m) + \\
&& + b_{-\frac{n}{2}+l} (\Delta',g',g'_{\Lambda^*},m) - b_{-\frac{n}{2}+l}
(i^*\Delta'i,g',g_{\Lambda^*},m) + \\
&& + b_{-\frac{n}{2}+l} (i^*\Delta'i,g',g_{\Lambda^*},m) - b_{-\frac{n}{2}+l}
(U^*i^*\Delta'iU,g,g_{\Lambda^*},m)
\eea
where
$b_{-\frac{n}{2}+l} (\Delta',g',g'_{\Lambda^*},m)$,
$b_{-\frac{n}{2}+l} (i^*\Delta'i,g',g_{\Lambda^*})$
are explained in (10.28), (10.27), respectively. (10.34) vanishes according to
(10.26). (10.32) $\in L_1(M,g)$
according to the expressions (10.2) and 
$g' \in \comp^{1,r}(g)$.
We conclude this as in the proof of 10.1. Finally (10.33) $\in L_1$
according to 
$i^*\Delta'i = (i^*-1) \Delta' + \Delta'$,
$\mbox{coeff } (i^*\Delta'i) = \mbox{coeff } (i^*-1) + \mbox{coeff } (\Delta')$,
$(i^*-1) \in \Omega^{0,r} (\mbox{End } (\Lambda^*))$, 
the rules for calculating the heat kernel expansion and according to the
module structure theorem.
\qed
\end{proof}

\begin{theorem}
Let $(M^n,g)$ be open, satisfying 
$(I)$, $(B_k)$, $k \ge r+1 > n+3$,
$g' \in \comp^{1,r+1} (g)$, $\Delta = \Delta (g,g_{\Lambda^*})$, 
$\Delta' = \Delta (g',g'_{\Lambda^*})$ 
the graded Laplace operators, 
$U$, $i$ as in (7.54), and assume 
$\inf \sigma_e (\Delta |_{(ker \Delta)^\perp}) > 0$.

{\bf a)} Then for $t>0$
\[ e^{-t\Delta} - e^{-tU^*i^*\Delta'iU} \]
is of trace class.

{\bf b)} Denote 
$h = \tr (\phi_{\delta,\epsilon} (\Delta) - 
\phi_{\delta,\epsilon} (U^*i^*\Delta'iU))$
for 
$0 < 3\epsilon < \frac{\delta}{2}$, $\delta+3\delta < \frac{\mu}{2}$, 
$\mu = \inf \{ \mbox{nonzero spectrum of } \Delta, i^*\Delta'i \}$.
Then
\[ \zeta_q (s,\Delta,\Delta') := \frac{1}{\Gamma(s)} \int\limits^\infty_0
t^{s-1} [ \tr (e^{-t\Delta_q} - e^{-t(U^*i^*\Delta'_qiU)}) - h] \,\, dt \]
has a well defined meromorphic extension to $\re (s) > -1$ which is 
holomorphic in $s=0$.

{\bf c)} The relative analytic torsion
$\tau^a(M^n,g,g')$, 
\be
\log \tau^a (M^n,g,g') := \sum\limits^n_{q=0} (-1)^q q \cdot 
\zeta'_q (0,\Delta,\Delta')
\end{equation}
is well defined.
\end{theorem}

\begin{proof}
a) is just theorem 7.9. b) immediately follows from 10.6 and the proof of
10.4. c) is a consequence of b).
\qed
\end{proof}

\newpage
It is very easy to provide many classes of {\bf Examples}.

{\bf 1)} Let
$(H^{2k}_{-1}, g_{H^{2k}_{-1}}) = (\R^{2k}, dr^2 + (\sinh r)^2 
d \sigma^2_{S^{2k-1}})$
be the $2k$--dimensional real hyperbolic space and
$(\R^{2k},g')$
such that
$g \in \comp^{1,r+1}(g_{H^{2k}_{-1}})$.
Then
$\tau^a (H^{2k}_{-1}, g_{H^{2k}_{-1}},g')$
is a well defined important non local relative invariant for the pair
$(H^{2k}_{-1}, (\R^{2k}, g'))$.

{\bf 2)} This can be e.g. generalized to finite connected sums
$(M^{2k},g) = \left( \underset{1} {\overset {m} {\#}} H^{2k}_{-1}, 
g \underset{1} {\overset {m} {\#}} H^{2k}_{-1} \right)$
and
$g' \in \comp^{1,r+1} \left( \underset{1} {\overset {m} {\#}} 
g_{H^{2k}_{-1}} \right)$.

{\bf 3)} A further generalization is 
$(M^{2k},g) =$ (compact topological and metrical perturbation of
$\underset{1} {\overset {m} {\#}} H^{2k}_{-1}, g_{M^{2k}}$), 
$g'_M \in \comp^{1,r+1} (g_M)$.

{\bf 4)} Still more general manifolds $(M^{2k},g)$ which fall into the
domain of 10.7 are those with warped product metric at infinity satisfying
the assumptions of 10.7 which can in this case easily be controlled.

{\bf 5)} Multiplication by $S^1$ yields odd dimensional examples.
\qed

\bi
Finally we turn to the relative $\eta$--invariant. On a closed manifold
$(M^n,g)$ and for a generalized Dirac operator the $\eta$--function is
defined as 
\be
\eta_D(s) := \sum\limits_{{\scriptsize \begin{array}{c} \lambda \in \sigma(D) \\ 
\lambda \neq 0 \end{array}}} \frac{\mbox{sign } \lambda}{|\lambda|^s} =
\frac{1}{\Gamma \left( \frac{s+1}{2} \right)} \int\limits^\infty_0
t^{\frac{s-1}{2}} \,\, \tr ( D e^{-tD^2} ) \,\, dt.
\end{equation}
$\eta_D(s)$ is defined for $\re(s) > n$, 
it has a meromorphic extension to $\C$ with isolated simple 
poles and the residues at all poles are locally computable. 
$\Gamma \left( \frac{s+1}{2} \right) \cdot \eta_D(s)$
has its poles at
$\frac{n+1-\nu}{2}$ 
for $\nu \in \N$.
One cannot conclude directly $\eta$ is regular at $s=0$ since
$\Gamma (u)$ is regular at $u=\frac{1}{2}$, i.e. 
$\Gamma \left( \frac{s+1}{2} \right)$
is regular at $s=0$. But one can show in fact using methods of algebraic
topology that $\eta(s)$ is regular at $s=0$. A purely analytical proof
for this is presently not known (cf. [16], p. 114/115).

(10.36) does not make sense on open manifolds. But we are able to define
a relative $\eta$--function and under an additional assumption the relative
$\eta$--invariant. $D e^{-tD^2}$ is an integral operator with heat kernel
$D_p W_D (t,m,p)$
which has at the diagonal a well defined asymptotic expansion 
(cf. [16], p. 75, lemma 1.9.1 for the compact case)
\be
\tr D_p W(t,m,p) |_{p=m} \underset{t \ra 0^+} \sim \sum\limits_{l \ge 0}
b_{-\frac{n+l}{2}} (D^2,D,m) t^{\frac{-n+l}{2}}. 
\end{equation}
In [3] has been proved that the heat kernel expansion on closed manifolds
also holds on open manifolds with the same coefficients (it is a local
matter) independent of the trace class property. The (simple) proof there
is carried out for 
$e^{-tD^2}$, $\tr W(t,m,m)$, 
but can be word by word repeated for
$D e^{-tD^2}$, $D W(t,m,m)$. 
The rules for calculating the
$b_{\frac{-n+l}{2}} (D^2,D,m)$
are quite si\-mi\-lar to them for
$b_{\frac{-n+l}{2}} (D,m)$
(cf. [16], Lemma 1.9.1). 
We sum up these considerations in

\begin{prop}
Let $E' \in \gencomp^{1,r+1}_{L,diff}(E)$, $r+1 > n+3$.
Then for $t>0$
\be
e^{-tD^2} D - e^{-t(U^*i^*D'iU)^2} (U^*i^*D'iU)
\end{equation}
is of trace class, for $t \ra 0^+$ there exists an asymptotic expansion
\be
\tr (e^{-tD^2} D - e^{-t(U^*i^*D'iU)^2} (U^*i^*D'iU)) =
\sum\limits^{n+3}_{l=0} \int\limits_M a_{\frac{-n+l}{2}} (m) \,\,
dvol_m (g) t^{\frac{-n+l}{2}} + o(t)
\end{equation}
\end{prop}

\begin{proof}
The first assertion is just theorem 7.10. (10.39) can just be derived as
proposition 10.6.
\qed
\end{proof}

We recall from [19] the following

\begin{prop}
Assume that $D$ and $\tilde{D'} = U^*i^*D'iU$ satisfy (10.37), (10.38)
and that the spectra of $D$ and $D'$ have a common gap $[a,b]$, 
$\left( \sigma(D) \cup \sigma(\tilde{D'}) \right) \cap [a,b] = \emptyset$.
Then there exists a spectral shift function
$\xi(\lambda) = \xi \left( \lambda,D,\tilde{D'} \right)$
having the following properties.

1. $\xi \in L_{1,loc} (\R)$ and $\xi(\lambda)=0$ for $\lambda \in [a,b]$.
\hfill {\rm (10.40)} 

2. For all $\phi \in C^\infty_c (\R)$, $\phi(D)-\phi ( \tilde{D'} )$ 
is a trace class operator and

\qquad $\tr \left( \phi(D)-\phi(\tilde{D'}) \right) = \int\limits_\R \phi'(\lambda)
\xi(\lambda) \,\, d \lambda$.
\hfill {\rm (10.41)}

3. $\tr \left( e^{-tD^2} D - e^{-t {\tilde{D'}}^2} \tilde{D'} \right) = 
\int\limits_\R \frac{d}{d \lambda} (\lambda e^{-t\lambda^2})
\xi(\lambda) \,\, d \lambda$.
\hfill {\rm (10.42)}

\qed
\end{prop}

\setcounter{equation}{42}

\begin{prop}
Assume $E' \in \gencomp^{1,r+1}_{L,diff} (E)$ and
$\inf \sigma_e (D^2 |_{(ker D^2)^\perp})$ $ > 0$.
Then there exists $c>0$ s. t. 
\be
\tr (e^{-t D^2} D - e^{-t(U^*i^*D'iU)^2} (U^*i^*D'iU)) = O(e^{-ct}).
\end{equation}
\end{prop}

\begin{proof}
We conclude as in the proof of 10.3 that there exists $\mu > 0$ s. t.
$\sigma(D) \cap \sigma(i^*D'i) \cap ([-\mu, -\frac{1}{\nu}] \cup 
[\frac{1}{\nu}, \mu]) = \emptyset$ 
for all $\nu \ge \nu_0$.
Hence, according to (10.40), 
$\int\limits^\mu_{-\mu} \frac{d}{dt} (\lambda e^{-t\lambda^2}) 
\xi(\lambda) \, d \lambda = 0$ 
and
\bea
&& | \tr (e^{-tD^2} D - e^{-t(U^*i^*D'iU)^2} ) | \non \\
&& = \left| 
\int\limits^\mu_{-\infty} \frac{d}{d \lambda} 
(\lambda e^{-t \lambda^2}) \xi (\lambda) \,\, d \lambda 
+ \int\limits^\infty_\mu \frac{d}{d \lambda} 
(\lambda e^{-t \lambda^2}) \xi (\lambda) \,\, d \lambda 
\right| \non \\
&& \le e^{-t\frac{\lambda^2}{2}} \left[
\int\limits^\mu_{-\infty} e^{-t \frac{\lambda^2}{2}} 
|1-2t\lambda| |\xi(\lambda)| \,\, d \lambda
+ \int\limits^\infty_\mu e^{-t \frac{\lambda^2}{2}} 
|1-2t\lambda| |\xi(\lambda)| \,\, d \lambda
\right] \non \\
&& = C \cdot e^{-t \frac{\mu^2}{2}} . \hspace{10cm} \Box \non
\eea
\end{proof}

\begin{theorem}
Assume $E' \in \gencomp^{1,r+1}(E)$, $k \ge n+1 > n+3$
and $\inf \sigma_e (D^2 |_{(ker D^2)^\perp}) > 0$.
Then there is a well defined relative $\eta$--function
\be
\eta(s,D,D') := \frac{1}{\Gamma \left( \frac{s+1}{2} \right)}
\int\limits^\infty_0 t^{\frac{s-1}{2}} \,\, \tr ( D e^{-tD^2}
- U^*i^*D'iU e^{-t(U^*i^*D'iU)^2} ) \,\, dt
\end{equation}
which is defined for $\re s > \frac{n}{2}$ and admits a 
meromorphic extension to $\re s > -5$. It is holomorphic 
at $s=0$ if the coefficient
$\int a_{- \frac{1}{2}} (m) \, d vol_m (g)$
of $t^{-\frac{1}{2}}$ equals to zero. Then there is a well 
defined relative $\eta$--invariant of the pair $(E, E')$.
\end{theorem}

\begin{proof}
We write again $U^*i^*D'iU = \tilde{D'}$.
Then according to (10.38),
\bea
\eta(s,D,D') &=& \frac{1}{\Gamma \left( \frac{s+1}{2} \right)}
\int\limits^\infty_0 t^{\frac{s-1}{2}} \left[ 
\sum\limits^{n+3}_{l=0} \int\limits_M a_{\frac{-n+l}{2}} \,\, d vol_m (g)
t^{\frac{n+l}{2}} + O(t^{\frac{4}{2}}) \,\, dt
\right] \non \\
&=& \frac{1}{\Gamma \left( \frac{s+1}{2} \right)}
\left[ \int\limits^1_0 \cdots dt + \int\limits^\infty_1 \cdots dt \right] \non \\
&=& \frac{1}{\Gamma \left( \frac{s+1}{2} \right)}
\sum\limits^{n+3}_{l=0} \frac{1}{\frac{s}{2}-\frac{n}{2}+\frac{l}{2}+\frac{1}{2}}
\int\limits_M a_{\frac{-n+l}{2}} \,\, d vol_m (g) +  \\
&+& \frac{1}{\Gamma \left( \frac{s+1}{2} \right)}
\int\limits^t_0 t^{\frac{s-1}{2}} O(t^{\frac{4}{2}}) \,\, dt + \\
&+& \frac{1}{\Gamma \left( \frac{s+1}{2} \right)}
\int\limits^\infty_1 t^{\frac{s-1}{2}} \,\, \tr (D e^{tD^2} - 
\tilde{D'} e^{-t {\tilde{D'}}^2} ) \,\, dt
\eea
We infer from (10.43) that (10.47) is holomorphic in $\C$. (10.46)
is holomorphic in $\re s > -5$. (10.45) admits a meromorphic extension
to $\C$. $\eta(s,D,D')$ is holomorphic at $s=0$ if the coefficient
$\int\limits_M a_{-\frac{1}{2}} (m) \, d vol_m (g)$
equals to zero.
\qed
\end{proof}

This finishes our analytical approach to relative Dirichlet series for pairs
of open manifolds arising from geometry. In a forthcoming paper a combinatorial
approach will be presented.

\bi
\bi
J\"urgen Eichhhorn

Institut f\"ur Mathematik und Informatik

Friedrich--Ludwig--Jahn--Strasse 15a

D-17487 Greifwald

Germany

eichhorn@uni-greifswald.de


\begin{thebibliography}{99} 

\bibitem{} {\sc N. V. Borisov, W. M\"uller} and {\sc R. Schrader}, 
{\em Relative index theorems and sypersymmetric scattering theory}, 
Comm. Math. Phys. 114 (1983), 475--513.

\bibitem{} {\sc U. Bunke},
{\em Relative index theory}, 
J. Funct. Analysis 105 (1992), 63--67.
            
\bibitem{} {\sc U. Bunke}, 
{\em Spektraltheorie von Diracoperatoren auf offenen Mannigfaltigkeiten}, 
thesis Greifswald University 1991.
            
\bibitem{} {\sc J. Cheeger, M. Gromov} and {\sc M. Taylor}, 
{\em Finite propagation speed, kernel estimates for functions of the
Laplace operator, and the geometry of complete manifolds}, 
J. Diff. Geom. 17 (1992), 15--35.

\bibitem{} {\sc H. Donnelly}, 
{\em Essential spectrum and heat kernel}, 
J. Functional Ana\-ly\-sis 75 (1987), 326--381.

\bibitem{} {\sc J. Eichhorn}, 
{\em Uniform Structures of Metric Spaces and Open Manifols}, 
Results in Mathematics 40 (2001), 144--191.

\bibitem{} {\sc J. Eichhorn}, 
{\em Invariants for proper metric spaces and open Riemannian manifolds}, 
to appear in Math. Nachr.

\bibitem{} {\sc J. Eichhorn}, 
{\em Bordism theory for open manifolds}, 
preprint, Greifswald 2001.

\bibitem{} {\sc J. Eichhorn}, 
{\em The manifold structure of maps between open manifolds}, 
Annals of Global Analysis and Geometry 11 (1993), 253--300.
            
\bibitem{} {\sc J. Eichhorn},
{\em Gauge theory on open manifolds of bounded geometry}, 
Int. Journ. Mod. Physics 7 (1992) 3927--3977.

\bibitem{} {\sc J. Eichhorn}, 
{\em Elliptic operators on noncompact manifolds}, 
Teubner Texte zur Mathematik 106 (1988), 4--169.
            
\bibitem{} {\sc J. Eichhorn}, 
{\em Spaces of Riemannian metrics on open manifolds}, 
Results in Mathematics 27 (1995), 256--283.
            
\bibitem{} {\sc J. Eichhhorn, J. Fricke},  
{\em The module structure theorem for Sobolev spaces on open manifolds}, 
Math. Nachr. 194 (1998), 35--47.

\bibitem{} {\sc J. Eichhorn, Yu. Kordjukov},  
{\em Differential operators with Sobolev coefficients}, 
in preparation.
             
\bibitem{} {\sc J. Eichhorn}, 
{\em The boundedness of connection coefficients and their derivatives}, 
Math. Nachr. 152 (1991), 145--158.

\bibitem{} {\sc P. Gilkey}, 
{\em Invariance theory, the heat equation and the Atiyah--Singer index theorem}, 
Studies in Advanced Mathematics, Boca Raton 1995.

\bibitem{} {\sc M. Gromov} and {\sc H. B. Lawson}, 
{\em Positive scalar curvature and the Dirac operator on complete
Riemannian manifolds},
Publ. Math. IHES 58 (1983), 295--408.

\bibitem{} {\sc E. Hebey},   
{\em Sobolev spaces on Riemannian manifolds},
Lecture Notes in Mathematics 1635, Berlin 1996.

\bibitem{} {\sc W. M\"uller}, 
{\em Relative zeta functions, relative determinants and scattering theory}, 
Comm. Math. Phys. 192 (1998), 309--347.
            
\bibitem{} {\sc G. Salomonsen},  
{\em Equivalence of Sobolev Spaces}, 
preprint Arhus 1999.

\bibitem{} {\sc W. S. Wladimirov}, 
{\em Equations of Mathematical Physics}, 
Moscow 1971.



\end{thebibliography}
\end{document}